\documentclass[10pt]{amsart}

\usepackage{amsfonts}
\usepackage{amscd}
\usepackage{amsmath}
\usepackage{amsthm}
\usepackage{amssymb}
\usepackage{amsxtra}
\usepackage{latexsym}

\input xy
\xyoption{all}

\newfont{\geschwollen}{cmsy10 at 10pt}

\newcommand{\cA}{{\mathbb A}}

\newcommand{\cC}{{\mathbb C}}
\newcommand{\cZ}{{\mathbb Z}}

\newcommand{\cN}{{\mathbb N}}

\newcommand{\cT}{{\mathbb T}}
\newcommand{\cS}{{\mathbb S}}

\newcommand{\caC}{{\mathcal C}}
\newcommand{\caD}{{\mathcal D}}
\newcommand{\caB}{{\mathcal B}}

\newcommand{\caE}{{\mathcal E}}

\newcommand{\caN}{{\mathcal N}}
\newcommand{\caS}{{\mathcal S}}
\newcommand{\caO}{{\mathcal O}}

\newcommand{\caT}{{\mathcal T}}
\newcommand{\caP}{{\mathcal P}}
\newcommand{\caL}{{\mathcal L}}

\newcommand{\comp}{\mathit{Comp}}

\newcommand{\inj}{\text{-}\mathrm{inj}}
\newcommand{\cof}{\text{-}\mathrm{cof}}
\newcommand{\cell}{\text{-}\mathrm{cell}}

\newcommand{\sset}{\mathrm{\bf SSet}}
\newcommand{\set}{\mathrm{\bf Set}}
\newcommand{\ho}{\mathrm{Ho}\,}

\newcommand{\Hom}{\mathrm{Hom}}
\newcommand{\uhom}{\underline{\mathrm{Hom}}}
\newcommand{\Ab}{\mathrm{\bf Ab}}

\newcommand{\op}{\mathrm{op}}
\newcommand{\Op}{\mathrm{Op}}

\newcommand{\mmod}{\text{--}\mathit{Mod}}

\newcommand{\eins}{1 \hspace{-2.3pt} \mathrm{l}}

\newcommand{\id}{\mathrm{Id}}

\newcommand{\colim}{\mathrm{colim}}

\newcommand{\Alg}{{\mathrm{Alg}}}

\newcommand{\val}{\mathrm{val}}
\newcommand{\dc}{\mathrm{dc}}

\newcommand{\old}{\mathrm{old}}
\newcommand{\new}{\mathrm{new}}
\newcommand{\ch}{\mathrm{ch}}
\newcommand{\conc}{\mathrm{conc}}
\newcommand{\Aut}{\mathrm{Aut}}
\newcommand{\am}{\mathrm{am}}
\newcommand{\Ass}{\mathrm{Ass}}
\newcommand{\Mor}{\mathrm{Mor}}
\newcommand{\Cat}{\mathrm{\bf Cat}}
\newcommand{\End}{\mathrm{End}}
\newcommand{\Comm}{\mathrm{Comm}}
\newcommand{\boxstar}{\Box_*}
\newcommand{\Gpd}{\mathrm{\bf Gpd}}

\newcounter{appnum}

\newtheorem{rem}{Remark}

\newtheorem{prop}{Proposition}
\newtheorem{lem}{Lemma}
\newtheorem{thm}{Theorem}
\newtheorem{cor}{Corollary}
\newtheorem{defi}{Definition}
\newtheorem{ex}{Example}
\newtheorem{assump}{Assumption}

\begin{document}

\title{Operads, algebras and modules in general model categories}
\author{Markus Spitzweck}
\date{January 2001}
\address{Department of Mathematics \\ Bonn University \\
Bonn, Beringstr. 3 \\ Germany}
\email{msp@math.uni-bonn.de}

\begin{abstract}

In this paper we develop the theory of operads, algebras
and modules in cofibrantly generated symmetric monoidal model
categories. We give $J$-semi model structures,
which are a slightly weaker version of model structures, for operads and
algebras and model  structures for modules. We prove homotopy invariance
properties for the categories of algebras and modules.
In a second part we develop the theory of $\cS$-modules and algebras of
\cite{ekmm} and \cite{km}, which allows a general homotopy theory for
commutative algebras and pseudo unital symmetric monoidal categories
of modules over them. Finally we prove a base change and projection formula.

\end{abstract}

\maketitle

\tableofcontents

\section{Introduction}

Recently important new applications of model categories appeared,
the most notable one maybe in the work of Voevodsky and others
on the $\cA^1$-local stable homotopy category of schemes.
But also for certain questions in homological algebra model
categories are quite useful, for example when one deals with
unbounded complexes in abelian categories.
In topology, mainly in the stable homotopy category, 
one is used to deal with objects having additional structures, 
for example modules over ring spectra.
The work of \cite{ekmm} made it possible to handle
commutativity appropriately, namely
the special properties of the linear isometries operad lead
to a strictly associative and commutative tensor product for
modules over $E_\infty$-ring spectra. As a consequence
many constructions in topology became more elegant or even
possible at all (see \cite{modern}).
Moreover the category of $E_\infty$-algebras
could be examined with homotopical methods because this
category carries a model structure. In \cite{km} a parallel theory
in algebra was developed (see \cite{may}).

Parallel to the achievements in topology the abstract model
category theory was further developed (see \cite{hov-book}
for a good introduction to model categories, see also \cite{dhk}). 
Categories of algebras and of modules over algebras in monoidal
model categories have been considered (\cite{ss}, \cite{hov-mon}).
Also localization techniques for model categories have become
important, because they yield many new useful model structures
(for example the categories of spectra of \cite{hov-stab}).
The most general statement 
for the existence of localizations is given in \cite{hir}.

In all these situations it is as in topology desireable 
to be able to work in
the commutative world, i.e. with commutative algebras and modules over them.
Since a reasonable model structure for commutative
algebras in a given symmetric monoidal model category is quite
unlikely to exist the need for a theory of $E_\infty$-algebras
arises. Also for the category of modules over an $E_\infty$-algebra
a symmetric monoidal structure is important. 
One of the aims of this paper is to give adequate answers
to these requirements.

$E_\infty$-algebras are algebras over particular operads.
Many other interesting operads appeared in various
areas of mathematics, starting from the early application
for recognition principles of iterated loop spaces
(which was the reason to introduce operads), 
later for example to handle homotopy Lie algebras which are
necessary for general deformation theory, the operads appearing
in two dimensional conformal quantum field theory or the
operad of moduli spaces of stable curves in
algebraic geometry.
In many cases the necessary operads are only well
defined up to quasi isomorphism or another sort of
weak equivalence (as for example is the case for $E_\infty$-agebras),
therefore a good homotopy theory of operads is desireable.
A related question is then the invariance (up to homotopy)
of the categories of algebras over weakly equivalent operads
and also of modules over weakly equivalent algebras.
We will also give adequate solutions to these questions.
This part of the paper was motivated by and owes many ideas from
\cite{hin1} and \cite{hin2}.

So in the first part we will develop the theory of operads, 
algebras and modules in the general situation
of a cofibrantly generated symmetric monoidal model category
satisfying some technical conditions which are usually fulfilled.
Our first aim is to provide these categories with model structures.
It turns out that in general we cannot quite get model structures
in the case of operads and algebras,
but a slightly weaker structure which we call a $J$-semi model structure.
A version of this structure already appeared in \cite{hov-mon}.
To the knowledge of the author no restrictions arise 
in the applications when using
$J$-semi model structures instead of model structures.
The $J$-semi model structures are necessary since
the free operad and algebra functors are not linear
(even not polynomially).
These structures appear in two versions, 
an absolute one and a version
relative to a base category.

We have two possible conditions for an operad or an algebra
to give model structures on the associated categories of algebras
or modules, the first one is being cofibrant (which is in some
sense the best condition), and the second one being cofibrant
in an underlying model category. 

In the second part of the paper
we demonstrate that the theory of $\cS$-modules of \cite{ekmm} and \cite{km}
can also be developed in our context if the given symmetric monoidal
model category $\caC$ either 
receives a symmetric monoidal left Quillen functor from 
$\sset$ (i.e. is simplicial) or from $\comp_{\ge 0}(\Ab)$.
The linear isometries operad $\caL$ gives via one of these functors
an $E_\infty$-operad in $\caC$ with the same special
properties responsible for the good behavior
of the theories of \cite{ekmm} and \cite{km}.
These theories do not yield honest
units for the symmetric monoidal category of modules over 
$\caL$-algebras, and we have to deal with the same problem.
In the topological theory of \cite{ekmm} it is possible
to get rid of this problem, in the algebraic or simplicial one
it is not.
Nevertheless it turns out that the properties the unit
satisfies are good enough to deal with operads,
algebras and modules in the category of modules over a
cofibrant $\caL$-algebra.
This seems to be a little contraproductive, but we need this
to prove quite strong results on the behavior of algebras and
modules with respect to base change and projection morphisms.
These results are even new
for the cases treated in \cite{ekmm} and \cite{km}.

In an appendix we show that one can always
define a product
on the homotopy category of modules over
an $\caO$-algebra for an arbitrary $E_\infty$-operad $\caO$
without relying on the special properties of the
linear isometries operad, but we do not construct associativity
and commutativity isomorphisms in this situation!
In the case when $\cS$-modules are available
this product structure is naturally isomorphic to the one
defined using $\cS$-modules.

%The basic example of a symmetric monoidal category which
%is not simplicial is the category of complexes of
%modules over a commutative ring $R$ with unit,
%but it is ``equivalent'' (in a zick zack sense)
%to a simplicial one, namely to the stabilization via
%symmetric spectra of the simplicial symmetric monoidal model category
%of simplicial $R$-modules. This can be done with many other
%symmetric monoidal categories structured on this example,
%e.g. complexes of sheaves on a small site, complexes
%of Nisnevich sheaves with transfers,
%$\cA^1$-localizations thereof, etc.
%So the hypothesis of being simplicial is not really
%a restriction.

Our constructions have explicit applications, for example
for the $\cA^1$-homotopy categories of Voevodsky, for
triangulated categories of motives over a general base,
for the ``tangential base point'' constructions 
$\grave{\text{a}}$ la Grothendieck, Deligne and others
and its ``motivic'' interpretation, which we demonstrate
in a forthcoming paper,
to develop the theory of schemes in 
symmetric monoidal cofibrantly generated model
categories (see \cite{tv}), etc.

I would like to thank Bertrand Toen for many useful
discussions on the subject.

\section{Preliminaries}

We first review some standard arguments from model category theory
which we will use throughout the paper (see for the first part
e.g. the introduction to \cite{hov-mon}).

Let $\caC$ be a cocomplete category.
For a pushout diagram in $\caC$
$$\xymatrix{A \ar[r]^f & B \\ K \ar[u]_\varphi \ar[r]^g & L \ar[u]}$$
we call $f$ the pushout of $g$ by $\varphi$,
and we call $B$ the pushout of $A$ by $g$ with attaching map $\varphi$.
If we say that $B$ is a pushout of $A$ by $g$ and $g$ is
an object of $\caC$ then we mean that $B=g$ and $A$ need not
be defined in this case (the sense of this statement will become
clear in the statements describing pushouts
of operads and algebras over operads in a model category $\caC$).

Let $I$ a set of maps in $\caC$.
Let $I\inj$ denote the class of maps in $\caC$ which have the right
lifting property with respect to $I$, $I\cof$ the class of maps in $\caC$
which have the left lifting property with respect to $I\inj$ and
$I\cell$ the class of maps which are transfinite
compositions of pushouts of maps from $I$.
Note that $I\cell \subset I\cof$ and that 
$I\inj$ and $I\cof$ are closed under retracts.

Let us suppose now that the domains of the maps in $I$ are small
relative to $I\cell$. Then by the small object argument there
exists a functorial factorization of every map in $\caC$ into
a map from $I\cell$ followed by a map from $I\inj$. Moreover
every map in $I\cof$ is a retract of a map in $I\cell$ such that
the retract induces an isomorphism on the domains of the two maps.
Also the domains of the maps in $I$ are small relative to $I\cof$.

Now let $\caC$ be equipped with a symmetric monoidal structure
such that the product $\otimes:\; \caC \times \caC \rightarrow \caC$
preserves colimits (e.g. if the monoidal structure is closed).
We denote the pushout product of maps $f:\; A \rightarrow B$
and $g:\; C \rightarrow D$,
$$A \otimes D \sqcup_{A \otimes C} B \otimes C \rightarrow B \otimes D
\;\text{,}$$
by $f \Box g$.

For ordinals $\nu$ and $\lambda$ we use the convention that
the well-ordering on the product ordinal $\nu\times \lambda$
is such that the elements in $\nu$ have higher significance.
We will need the

\begin{lem} \label{cell-pushout-str}
Let $f: \; K_0 \rightarrow K_\mu=\colim_{i < \mu} K_i$ and
$g: \; L_0 \rightarrow L_\lambda=\colim_{i < \lambda} L_i$ be transfinite
compositions with transition maps $f_i:\; K_i \rightarrow K_{i+1}$
and $g_i:\; L_i \rightarrow L_{i+1}$.
Then the pushout product $f \Box g$ is a transfinite composition
$M_0 \rightarrow M_{\mu\times \lambda}=\colim_{i < \mu\times \lambda}
M_i$ over the product ordinal $\mu \times \lambda$ where the transition maps
$M_{(i,j)} \rightarrow M_{(i,j+1)}$ are pushouts by the maps
$f_i \Box g_j$.
\end{lem}
\begin{proof}
For any $(i,j) \le \mu\times \lambda$ define $M_{(i,j)}$ to be
the colimit of the diagram
$$\xymatrix{A_\mu \otimes B_0 & & A_{i+1} \otimes B_j & & A_i 
\otimes B_\lambda \\ & A_{i+1} \otimes B_0 \ar[lu] \ar[ru] & &
A_i \otimes B_j \ar[lu] \ar[ru]}\;\text{.}$$
Clearly $M_{(0,0)}=A_\mu \otimes B_0 \sqcup_{A_0 \otimes B_0}
A_0 \otimes B_\lambda$ is the domain and
$M_{\mu \times \lambda}=A_\mu \otimes B_\lambda$
the codomain of $f \Box g$. 
Moreover it is easy to see that the pushout of $M_{(i,j)}$
by $f_i \Box g_j$ with the obvious attaching map is canonically
isomorphic to $M_{(i,j+1)}$. Since $\otimes$ preserves colimits the
assignment $(i,j) \mapsto M_{(i,j)}$ is a transfinite composition.
\end{proof}

The pushout product is associative. For maps 
$f_i:\; A_i \rightarrow B_i$, $i=1,\ldots,n$, in
$\caC$ giving a map from the domain of $g:=\text{\LARGE $\Box$}_{i=1}^n
f_i$ to an object $X \in \caC$ is the same as to give maps $\varphi_j$
from the $$S_j:=(\bigotimes_{i=1}^{j-1} B_i) \otimes A_j \otimes
\bigotimes_{i=j+1}^n B_i$$ to $X$ for $j=1,\ldots,n$
such that $\varphi_j$ and $\varphi_{j'}$ ($j'>j$) coincide
on $$I_{j,j'}:=(\bigotimes_{i=1}^{j-1} B_i) \otimes A_j \otimes
(\bigotimes_{i=j+1}^{j'-1} B_i) \otimes A_{j'} \otimes
\bigotimes_{i=j'+1}^n B_i$$ after the obvious compositions.
We call the $S_j$ the {\em summands} of the domain of $g$
and the $I_{j,j'}$ the {\em intersections} of these summands.
Sometimes some of the $f_i$ will coincide. Then there is an
action of a product of symmetric groups on $g$, and the quotient
of a summand with respect to the induced action of the stabilizer
of this summand will also be called a summand (and similarly for
the intersections).

\vskip.2cm

For the rest of the paper we fix a cofibrantly generated
symmetric monoidal model category $\caC$ with generating cofibrations
$I$ and generating trivial cofibrations $J$.
For simplicity we assume that the domains of $I$ and $J$
are small relative to the whole category $\caC$.
The interested reader may weaken this hypothesis appropriately
in the statements below.

For a monad $\cT$ in $\caC$ we write
$\caC[\cT]$ for the category of $\cT$-algebras in $\caC$.
The following theorem summarizes the general method to equip
categories of objects in $\caC$ with ``additional structure''
with model structures (e.g. as in \cite{hov-mon}[Theorem 2.1]).

\begin{thm} \label{mon-mod-str}
Let $\cT$ be a monad in $\caC$, assume that $\caC[\cT]$
has coequalizers and suppose 
that every map in $\cT J \cell$,
where the cell complex is built in $\cC[\cT]$, is a weak equivalence
in $\caC$.
Then there is a cofibrantly generated model structure on
$\caC[\cT]$, where a map is a weak equivalence or fibration if
and only if it is a weak equivalence or fibration in $\caC$.
\end{thm}
\begin{proof}
We apply \cite{hov-book}[Theorem 2.1.19] with generating cofibrations
$\cT I$, generating trivial cofibrations $\cT J$
and weak equivalences the maps which are weak equivalences in $\caC$.

By \cite{mclane}[VI.2, Ex 2], $\caC[\cT]$ is complete and
by \cite{ttt}[9.3 Theorem 2] cocomplete.
Property 1 of \cite{hov-book}[Theorem 2.1.19] is clear,
properties 2 and 3 follow by adjunction from our smallness
assumptions on the domains of $I$ and $J$.
Since each element of $J$ is in $I\cof$,
hence a retract of a map in $I\cell$, each element
of $\cT J$ is in $\cT I \cof$, hence together with our
assumption we see that property 4 is fulfilled.
By adjunction $\cT I \inj$ (resp. $\cT J \inj$)
is the class of maps in $\caC[\cT]$ which are trivial fibrations
(resp. fibrations) in $\caC$.
Hence property 5 and the second alternative of 6 are fulfilled.
\end{proof}

In most of the cases we are interested in the hypothesis of
this theorem that every map in $\cT J \cell$ is a weak equivalence
won't be fulfilled. The reason is that we are considering
monads which are not linear. The method to circumvent this problem 
was found by Hovey in \cite{hov-mon}[Theorem 3.3]. He considers
categories which are not quite model categories. We will call them
semi model categories.

\begin{defi} \label{semi-mod-rel-def}
(I) A {\em $J$-semi model category over $\caC$} is a
left adjunction $F:\;\caC \rightarrow \caD$ and
subcategories of weak equivalences, fibrations and
cofibrations in $\caD$ such that the following axioms are
fulfilled:
\begin{enumerate}
\item The adjoint of $F$ preserves fibrations
and trivial fibrations.
\item $\caD$ is bicomplete and the two out of three and
retract axioms hold in $\caD$.
\item Cofibrations in $\caD$ have the left lifting property with
respect to trivial fibrations, and trivial cofibrations whose
domain becomes cofibrant in $\caC$ have the left lifting
property with respect to fibrations.
\item Every map in $\caD$ can be functorially factored into a cofibration
followed by a trivial fibration, and every map in $\caD$
whose domain becomes cofibrant in
$\caC$ can be functorially factored into a trivial
cofibration followed by a fibration.
\item Cofibrations in $\caD$ whose domain becomes cofibrant
in $\caC$ become cofibrations in $\caC$, and the initial
object in $\caD$ is mapped to a cofibrant object in $\caC$.
\item Fibrations and trivial fibrations are closed under pullback.
\end{enumerate}
We say that $\caD$ is {\em cofibrantly generated} if there are
sets of morphisms $I$ and $J$ in $\caD$ such that
$I\inj$ is the class of trivial fibrations and $J\inj$ the
class of fibrations in $\caD$ and if the domains of $I$ are
small relative to $I\cell$ and the domains of $J$ are small
relative to maps from $J\cell$ whose domain becomes
cofibrant in $\caC$.

$\caD$ is called {\em left proper (relative
to $\caC$)} if pushouts by cofibrations
preserve weak equivalences whose domain and codomain become
cofibrant in $\caC$ (hence all objects which appear
become cofibrant in $\caC$).
$\caD$ is called {\em right proper}
if pullbacks by fibrations preserve weak equivalences.

(II) A category $\caD$ is called a $J$-semi model category
if conditions (2) to (6) of Definition \ref{semi-mod-rel-def}
are fulfilled where the condition of becoming cofibrant in
$\caC$ is replaced by the condition of being cofibrant.

The same is valid for the definition of being cofibrantly generated
and of being right proper.
\end{defi}
(Note that the only reasonable property to require in a
definition for a $J$-semi model category
to be left proper, namely that weak equivalences between cofibrant
objects are preserved
by pushouts by cofibrations, is automatically fulfilled
as is explained below when we consider homotopy pushouts.)

{\bf Alternative: }
One can weaken the definition of a $J$-semi model category
(resp. of a $J$-semi model category over $\caC$) slightly 
by only requiring that a factorization of a map in $\caD$
into a cofibration followed by a trivial fibration should exist 
if the domain of this map is cofibrant (resp. becomes
cofibrant in $\caD$). We then include
into the definition of cofibrant generation that the cofibrations
are all of $I \cof$. Using this definition all
statements from section \ref{operads} on remain true
if one does not impose any further smallness assumptions
on the domains of $I$ and $J$. This follows in each of the cases
from the fact that the domains of $I$ and $J$ are small relative
to $I\cof$.

\vskip.2cm

Of course a $J$-semi model category over $\caC$ is a $J$-semi
model category.
There is also the notion of an $I$-semi (and also
$(I,J)$-semi) model category (over $\caC$), where the parts
of properties 3 and 4 concerning cofibrations
are restricted to maps whose domain is cofibrant
(becomes cofibrant in $\caC$).

\vskip.2cm

We summarize the main properties of a $J$-semi model category $\caD$
(relative to $\caC$) (compare also \cite{hov-mon}[p. 14]):

By the factorization property and the retract argument it follows that
a map is a cofibration if and only if it has the left lifting
property with respect to the trivial fibrations.
Similarly a map is a trivial fibration if and only if it has
the right lifting property with respect to the cofibrations.
These two statements remain true under the alternative definition
if $\caD$ is cofibrantly generated.

A map in $\caD$ whose domain is cofibrant
(becomes cofibrant in $\caC$)
is a trivial cofibration if and only if it has the left lifting
property with respect to the fibrations, and a map 
whose domain is cofibrant (becomes cofibrant in $\caC$)
is a fibration if and only if it has the right lifting property
with respect to the trivial cofibrations.

Pushouts preserve cofibrations (also under the alternative definition
if $\caD$ is cofibrantly generated).
Trivial cofibrations with cofibrant domain (whose domain becomes
cofibrant in $\caC$) are preserved under pushouts by maps with
cofibrant codomain (whose codomain becomes cofibrant in $\caC$).

In the relative case
the functor $F$ preserves cofibrations
(also in the alternative definition
if $\caD$ is cofibrantly generated),
and trivial cofibrations with cofibrant domain.

Ken Brown's Lemma (\cite{hov-book}[lemma 1.1.12]) remains true,
and its dual version has to be modified to the following
statement: Let $\caD$ be a $J$-semi model category (over $\caC$) and
$\caD'$ be a category with a subcategory of weak equivalences which
satisfies the two out of three property. Suppose $F:\; \caD \rightarrow
\caD'$ is a functor which takes trivial fibrations between fibrant objects
with cofibrant domain (whose domain becomes cofibrant in $\caC$)
to weak equivalences. Then $F$ takes all weak equivalences
 between fibrant objects
with cofibrant domain (whose domain becomes cofibrant in $\caC$)
to weak equivalences.

\vskip.2cm

We define cylinder and path objects and the various versions
of homotopy as in \cite{hov-book}[Definition 1.2.4].
Cylinder and path objects exist for cofibrant objects
(for objects which become cofibrant in $\caC$).

We give the $J$-semi version of \cite{hov-book}[Proposition 1.2.5]:

\begin{prop}
Let $\caD$ be a $J$-semi model category (over $\caC$) and let
$f,g:\; B \rightarrow X$ be two maps in $\caD$.
\begin{enumerate}
\item If $f \overset{l}{\sim} g$ and $h: \; X 
\rightarrow Y$, then $hf \overset{l}{\sim} hg$. Dually, if
$f \overset{r}{\sim} g$ and $H:\; A \rightarrow B$, then
$fh \overset{r}{\sim} gh$.
\item Let $h:\; A \rightarrow B$ and suppose $A$ and $B$ are cofibrant
(become cofibrant in $\caC$) and $X$ is fibrant. Then
$f \overset{l}{\sim} g$ implies $fh \overset{l}{\sim} gh$.
Dually, let $h:\; X \rightarrow Y$. Suppse $X$ and $Y$ are cofibrant
(become cofibrant in $\caC$) and $B$ is cofibrant. Then
$f \overset{r}{\sim} g$ implies $hf \overset{r}{\sim} hg$.
\item If $B$ is cofibrant, then left homotopy is an
equivalence relation on $\Hom(B,X)$.
\item If $B$ is cofibrant and $X$ is cofibrant (becomes cofibrant
in $\caC$), then $f \overset{l}{\sim} g$ implies
$f \overset{r}{\sim} g$. Dually, if $X$ is fibrant and $B$ is cofibrant
(becomes cofibrant in $\caC$), then
$f \overset{r}{\sim} g$ implies $f \overset{l}{\sim} g$.
\item If $B$ is cofibrant and $h:\; X \rightarrow Y$ is a trivial fibration
or weak equivalence between fibrant objects with $X$ cofibrant
(such that $X$ becomes cofibrant in $\caC$), then $h$ induces an isomorphism
$$\Hom(B,X)/\overset{l}{\sim} \;\; \overset{\cong}{\longrightarrow} \;\;
\Hom(B,Y)/\overset{l}{\sim}\;\text{.}$$
Dually, suppose $X$ is fibrant and cofibrant (becomes cofibrant
in $\caC$) and $h:\; A \rightarrow B$ is a trivial cofibration
with $A$ cofibrant (such that $A$ becomes cofibrant in $\caC$) or
a weak equivalence between cofibrant objects, then $h$ induces
an isomorphism
$$\Hom(B,X)/\overset{r}{\sim} \;\; \overset{\cong}{\longrightarrow} \;\;
\Hom(A,X)/\overset{r}{\sim}\;\text{.}$$
\end{enumerate}
\end{prop}
This Proposition is also true for the alternative definition of a
$J$-semi model category (over $\caC$).
We changed the order between 4 and 5, because it is a priori not
clear that right homotopy is an equivalence relation (under suitable
condition), this follows only after comparison with the left homotopy
relation.

As in \cite{hov-book}[Corollary 1.2.6 and 1.2.7]
it follows that if $B$ is cofibrant and $X$ is fibrant and cofibrant
(becomes cofibrant in $\caC$), then left and right homotopy coincide
and are equivalence relations on $\Hom(B,X)$ and the homotopy
relation on $\caD_{cf}$ is an equivalence relation and compatible with
composition.
The statement of \cite{hov-book}[Proposition 1.2.8] that a map
in $\caD_{cf}$ is a weak equivalence if and only if it is a
homotopy equicalence is proved exactly in the same way. The same holds for
the fact that $\ho \caD_{cf}$ is naturally isomorphic to
$\caD_{cf}/\sim$ (\cite{hov-book}[Corollary 1.2.9]).
Finally the existence of the cofibrant and fibrant replacement functor
$RQ$ implies that the map $\ho \caD_{cf} \rightarrow
\ho \caD$ is an equivalence.

\begin{defi}
A functor $L: \; \caD \rightarrow \caD'$ between $J$-semi model
categories is a {\em left Quillen functor} if it has a right adjoint
and if the right adjoint preserves fibrations and trivial fibrations.
\end{defi}
Of course in the relative situation $F$ is a left Quillen functor.
We show that
a left Quillen functor induces an adjunction between the homotopy categories
(also when we use the alternative definition).
$L$ preserves (trivial) cofibrations between cofibrant objects,
hence by Ken Brown's Lemma it preserves weak equivalences between
cofibrant objects. This induces a functor $\ho \caD \rightarrow
\ho \caD'$. By the dual version of Ken Brown's Lemma the adjoint
of $L$ preserves weak equivalences between fibrant and cofibrant objects
which gives a functor $\ho \caD' \rightarrow \ho \caD$.
One easily checks that $L$ preserves cylinder objects on
cofibrant objects and that the adjoint of $L$ preserves path objects
on fibrant objects. As in Lemma \cite{hov-book}[Lemma 1.3.10] it follows
that on the derived functors between $\ho \caD$ and $\ho \caD'$ there
is induced a natural derived adjunction.

\vskip.2cm

Next we are going to consider Reedy model structures and homotopy
function complexes.
We have the analogue of \cite{hov-book}[Theorem 5.1.3]:

\begin{prop}
Let $\caD$ be a $J$-semi model category and $\caB$ be a direct
category. Then the diagram category $\caD^\caB$ is a $J$-semi
model category with objectwise weak equivalences and fibrations
and where a map $A \rightarrow B$ is a cofibrations if and only if
the maps $A_i \sqcup_{L_i A} L_i B \rightarrow B_i$ are
cofibrations for all $i \in \caB$.
\end{prop}
\begin{proof}
As in \cite{hov-book}[Proposition 5.1.4] one shows that
cofibrations have the left lifting property with respect to trivial
fibrations. Then it follows that if $A \rightarrow B$ is a map
in $\caD^\caB$ with $A$ cofibrant such that the maps
$A_i \sqcup_{L_i A} L_i B \rightarrow B_i$ are (trivial) cofibrations
then the map $\colim A \rightarrow \colim B$ is a (trivial)
cofibration in $\caD$. So a good trivial cofibration (definition as in
the proof of \cite{hov-book}[Theorem 5.1.3]) with cofibrant domain is
a trivial cofibration and trivial cofibrations with cofibrant domain
have the left lifting property with respect to fibrations.
We then can construct functorial factorizations into a good trivial
cofibration followed by a fibration for maps with cofibrant domain
as in the proof of \cite{hov-book}[Theorem 5.1.3]) and also the
factorization into a cofibration followed by a trivial fibration
(for the alternative definition for maps with cofibrant domain).
It follows that a trivial cofibration with cofibrant domain is
a good trivial cofibration.
All other properties are immediate.
\end{proof}

Similarly but easier we have that for an inverse category $\caB$
the diagram category $\caD^\caB$ is a $J$-semi model category.

We can combine both results as in \cite{hov-book}[Theorem 5.2.5]
to get

\begin{prop} Let $\caD$ be a $J$-semi model category and
$\caB$ a Reedy category. Then $\caD^\caB$ is a $J$-semi model category
where a map $f:\; A \rightarrow B$ is a weak equivalence if and only if
it is objectwise a weak equivalence, a cofibration if and only if
the maps $A_i \sqcup_{L_i A} L_i B \rightarrow B_i$ are cofibrations
and a fibration if and only if the maps
$A_i \rightarrow B_i \times_{M_i B} M_i A$ are fibrations.
\end{prop}

It is easily checked that cosimplicial and simplicial frames
(see \cite{hov-book}[Definition 5.2.7])
exist on cofibrant objects. In the following we denote by
$A^\bullet$ and $A_\bullet$ functorial cosimplicial and simplicial
frames on cofibrant $A \in \caD$.
We are going to equip the category $\caD_{cf}$ with a strict
$2$-category structure $\caD_{cf}^{\le 2}$
with underlying $1$-category $\caD_{cf}$
and with associated homotopy category $\ho \caD_{cf}$.
Let $A,B \in \caD_{cf}$. As in \cite{hov-book}[Proposition 5.4.7]
there are weak equivalences
$$\Hom_\caD(A^\bullet,B) \rightarrow 
\text{diag}(\Hom_\caD(A^\bullet, B_\bullet))
\leftarrow \Hom_\caD(A, B_\bullet)$$
in $\sset$ which are isomorphisms in degree $0$,
and we define the morphism category
$\Hom_{\caD_{cf}^{\le 2}}(A,B)$ to be the groupoid associated to
one of these simplicial sets. By the groupoid associated to a $K \in \sset$
we mean the groupoid with set of objects $K[0]$ and set of morphisms
$\Hom(x,y)$ for $x,y \in K[0]$ the homotopy classes of paths from
$x$ to $y$ in the topological realization of $K$.
We have to give composition functors
$$\Hom_{\caD_{cf}^{\le 2}}(A,B) \times \Hom_{\caD_{cf}^{\le 2}}(B,C)
\rightarrow Hom_{\caD_{cf}^{\le 2}}(A,C)\;\text{.}$$
These are the normal composition on objects and are induced on
the morphisms by the map of simplicial sets
$$\Hom_\caD(A^\bullet,B) \times \Hom_\caD(B,C_\bullet)
\rightarrow \text{diag}(\Hom_\caD(A^\bullet,C_\bullet))\;\text{.}$$
In the following we write $\underset{0}{\circ}$ for the composition
of $2$-morphisms over objects and $\underset{1}{\circ}$
for the composition of $2$-morphisms over $1$-morphisms.
We claim that for $A,B,C \in \caD_{cf}$, morphisms
$f,g:\; A \rightarrow B$, $f',g':\; B \rightarrow C$ and $2$-morphisms
$\varphi:\; f \rightarrow g$, $\psi:\; f' \rightarrow g'$ we
have $$\psi \underset{0}{\circ} \varphi =
(\id_{f'} \underset{0}{\circ} \varphi) \underset{1}{\circ}
(\psi \underset{0}{\circ} \id_g) = (\psi \underset{0}{\circ} \id_f)
\underset{1}{\circ} (\id_{g'} \underset{0}{\circ} \varphi)\;\text{.}$$
This follows from the corresponding equation of homotopy classes
of paths in $\Hom_\caD(A^\bullet,B) \times \Hom_\caD(B,C_\bullet)$.
Moreover for a $1$-morphism $f'':\; C \rightarrow D$ we have
$(\id_{f''}\underset{0}{\circ} \psi) \underset{0}{\circ} \id_f=
\id_{f''}\underset{0}{\circ} (\psi \underset{0}{\circ} \id_f)$,
and the assignments $\Hom_{\caD_{cf}^{\le 2}}(B,C) \rightarrow
\Hom_{\caD_{cf}^{\le 2}}(B,D)$, $a \mapsto \id_{f''} \underset{0}{\circ} a$,
and $\Hom_{\caD_{cf}^{\le 2}}(B,C) \rightarrow
\Hom_{\caD_{cf}^{\le 2}}(A,C)$, $a \mapsto a \underset{0}{\circ} \id_f$,
are functors.
From these three properties it follows that $\underset{0}{\circ}$
is associative and that $\underset{0}{\circ}$ and
$\underset{1}{\circ}$ are compatible.
Hence $\caD_{cf}^{\le 2}$ is a strict $2$-category. We set
$\ho^{\le 2} \caD := \caD_{cf}^{\le 2}$.
One can show that this $2$-category is
weakly equivalent to the $2$-truncation of the $1$-Segal category 
(see \cite{si-hi}) associated to $\caD$.

\vskip.2cm

Let $\llcorner$ be the category whose diagrams (i.e. functors
into another category) are the ``lower left triangles'',
and $\square$ the category whose diagrams are the commutative squares
like the square at the beginning of this section.
There is an obvious inclusion functor $\llcorner \rightarrow \square$.
For a category $\caD$
denote by $\caD^\llcorner$ (resp. $\caD^\square$) the category of
$\llcorner$-diagrams (resp. of $\square$-diagrams) in $\caD$.
There is a restriction functor $r:\;\caD^\square \rightarrow \caD^\llcorner$.

Let $\caD$ be a $J$-semi model category.
Then there is a canonical way to define a {\em homotopy pushout functor}
$$\hbox{$\sqcup$}:\; (\ho \caD)^\llcorner \rightarrow (\ho \caD)^\square$$
which sends a triangle
$\xymatrix{B & \\ A \ar[u] \ar[r] & C}$ to the square
$\xymatrix{B \ar[r] & B \hbox{$\sqcup$}_A C \\ 
A \ar[u] \ar[r] & C \ar[u]}$,
together with a natural isomorphism from $r \circ \sqcup$ to the identity.
This is done by lifting a triangle to a triangle in $\caD$ where
all objects are cofibrant and at least one map is a cofibration.
Then by the cube lemma (\cite{hov-book}[Lemma 5.2.6]),
which is also valid for $J$-semi model categories, the pushout does
not depend on the choices and indeed yields a well-defined square
in $\ho \caD$.
We call a square in $\ho \caD$ a {\em homotopy pushout square} if
it is in the essential image of the functor $\sqcup$.
This is by definition 
the same as to say that it is the image of a homotopy pushout
square in $\caD$, which is defined to be any commutative square
weakly equivalent to a pushout square
$$\xymatrix{B \ar[r] & D \\ A \ar[u]_f \ar[r]^g & C \ar[u]}$$
where all objects are cofibrant and $f$ or $g$ is a cofibration.

Taking $A$ to be an initial object in $\ho \caD$
(i.e. the image of an initial object in $\caD$) the product
$\sqcup_A$ gives the categorical coproduct on $\ho \caD$.
For general $A$ the homotopy pushout need not be a categorical
pushout in $\ho \caD$.

We show that the homotopy pushout {\em has} a categorical interpretation
in the $2$-category $\ho^{\le 2} \caD$: 
%In the following we use the term $2$-category
%and $2$-functor always in the weak sense.
Let for the moment $\caD$ be an arbitrary $2$-category.
A commutative square
$$\xymatrix{B \ar[r]^{g'} \ar@{=>}[dr]^\varphi & D \\ 
A \ar[u]_f \ar[r]^g & C \ar[u]_{f'}}$$
in $\caD$ is called a {\em homotopy pushout}, if for any object
$T \in \caD$ the square
$$\xymatrix{\Hom(D,T) \ar[r] \ar[d] & \Hom(B,T) \ar[d] \ar@{=>}[dl] \\
\Hom(C,T) \ar[r] & \Hom(A,T)}$$
is a homotopy pullback in the $2$-category $\Gpd$ of small groupoids.
We recall the definition of a homotopy pullback in $\Gpd$:
For a triangle $K \overset{f}{\rightarrow} G 
\overset{g}{\leftarrow} L$ in $\Gpd$ we define the 
{\em homotopy fibre product} $K\times_G^h L$ to be the groupoid with
objects triples $(x,y,\varphi)$, where $x \in K$, $y \in L$ and
$\varphi:\; f(x) \overset{\cong}{\rightarrow} g(y)$ an isomorphism,
and morphisms $(x,y,\varphi) \rightarrow (x',y',\varphi')$
pairs of morphisms $x \rightarrow x'$, $y \rightarrow y'$ making
the obvious diagram in $G$ commutative.
Now for a commutative square
$$\xymatrix{M \ar[r] \ar[d] & L \ar[d] \ar@{=>}[dl] \\
K \ar[r] & G}$$ in $\Gpd$
there is a canonical functor $M \rightarrow K\times_G^h L$,
and we say that the square is a homotopy pullback if this
functor is an equivalence.

Let $\caD$ be again a $J$-semi model category.
We claim now that the image of a homotopy pushout square
$$\xymatrix{B \ar[r] & D \\ A \ar[u]_f \ar[r]^g & C \ar[u]}$$
in $\caD_{cf}$ in $\ho^{\le 2} \caD$ is a homotopy pushout square
in the sense just defined: 
So let $T \in \caD_{cf}$. Then $\Hom(\text{square},T_\bullet)$ is a homotopy
pullback square in $\sset$, since, if $f$ is a cofibration with $A$
cofibrant, the map $\Hom(f,T_\bullet)$ is a fibration in $\sset$.
As is easily verified the functor $\sset \rightarrow \Gpd$ preserves
homotopy pullbacks, hence our claim follows.

So we have shown the following:
$\ho^{\le 2} \caD$ has categorical homotopy pushouts,
every homotopy pushout square in $\ho \caD$ comes
from one in $\ho^{\le 2} \caD$,
every homotopy pushout square in $\ho^{\le 2} \caD$ is
equivalent to the image of a homotopy pushout square in
$\caD$ and all such images are homotopy pushout squares in
$\ho^{\le 2} \caD$.

Note that it follows that for any $T \in \ho \caD$
and homotopy pushout square as above the map
$$\Hom(B \sqcup_A C,T) \rightarrow \Hom(B,T)\times_{\Hom(A,T)}
\Hom(C,T)\;\text{,}$$ where all homomorphism sets are in $\ho \caD$,
is always surjective.

There is a dual {\em homotopy pullback functor} $\times$
and the dual notion of a homotopy pullback square
in both $\ho \caD$ and $\ho^{\le 2} \caD$.

For a cofibrant object $A \in \caD$ the category $A \downarrow \caD$
of objects under $A$ is again a $J$-semi
model category. The $2$-functor $$\caD \rightarrow \Cat\;\text{,}$$
$$A \mapsto \ho((QA) \downarrow \caD)$$
where $QA \rightarrow A$ is a cofibrant replacement,
descents to a $2$-functor $$\ho^{\le 2} \caD \rightarrow \Cat\;\text{,}$$
$$A \mapsto D(A \downarrow \caD)$$
such that the image functors $f_*$ of all maps $f$ in $\ho^{\le 2} \caD$
have right adjoints $f^*$.
The functor $f_*$ preserves homotopy pushout squares,
and the functor $f^*$ preserves homotopy pullback and homotopy pushout
squares.
For $f:\; 0 \rightarrow A$ the map from an initial object
to an object in $\ho^{\le 2} \caD$ the functor $f^*:\; D(A \downarrow \caD)
\rightarrow \ho \caD$ factors through $A \downarrow \ho \caD$
and the map from $A$ to the image of the initial object in
$D(A \downarrow \caD)$ is an isomorphism.

Consider a commutative square
$$\xymatrix{B \ar[r]^{g'} \ar@{=>}[dr]^\varphi & D \\ 
A \ar[u]_f \ar[r]^g & C \ar[u]_{f'}}$$
in $\ho^{\le 2} \caD$. Let $E \in D(B \downarrow \caD)$.
There is a base change morphism 
$$g_*f^* E \rightarrow {f'}^*g'_* E$$
adjoint to the natural map
$f^*E \rightarrow f^*{g'}^*g'_* M\overset{\varphi}{\cong} g^*{f'}^* g'_* M$.
This base change morphism applied to diagrams
$$\xymatrix{B \ar[r] \ar@{=>}[dr] & C \\ 
A \ar[u] \ar[r]^\id & A \ar[u]}$$
enables one to construct a $2$-functor
$$(A \downarrow \ho^{\le 2} \caD) \rightarrow D(A \downarrow \caD)$$
which gives an equivalence after $1$-truncation of the left hand side.

\begin{rem}
The above construction should generalize
to give functors between (weak) $(n+1)$-categories
$$\ho^{\le n+1} \caD \rightarrow n-\Cat$$
$$A \mapsto D^{\le n}(A \downarrow \caD)\;\text{,}$$
where $\ho^{\le n+1}$ is the $(n+1)$-truncation of the $1$-Segal
category associated to $\caD$,
$n-\Cat$ is the $(n+1)$-category of $n$-categories
and $D^{\le n}(A \downarrow \caD):=\ho^{\le n} (QA \downarrow \caD)$
for $QA \rightarrow A$ a cofibrant replacement.
\end{rem}

There are dual constructions for objects over an object in $\caD$.

The following theorem is the main source to obtain $J$-semi
model categories.

\begin{thm} \label{mon-semi-mod-str}
Let $\cT$ be a monad in $\caC$ and assume that $\caC[\cT]$
has coequalizers. Suppose 
that every map in $\cT J \cell$ whose domain
is cofibrant in $\caC$
is a weak equivalence in $\caC$ and every map in
$\cT I \cell$ whose domain is cofibrant
in $\caC$ is a cofibration in $\caC$ (here in both cases
the cell complexes are built in $\caC[\cT])$.
Assume furthermore that the initial object in $\caC[\cT]$
is cofibrant in $\caC$.
Then there is a cofibrantly generated $J$-semi model structure on
$\caC[\cT]$ over $\caC$, where a map is a weak equivalence or fibration if
and only if it is a weak equivalence or fibration in $\caC$.
\end{thm}
\begin{proof}
We define the weak equivalences (resp. fibrations) as the maps in
$\caC[\cT]$ which are weak equivalences (resp. fibrations) as maps
in $\caC$. By adjointness the fibrations are $\cT J \inj$ and the
trivial fibrations are $\cT I \inj$.
We define the class of cofibrations to be $\cT I \cof$. Since
the adjoint of $\cT$ is the forgetful functor property 1 of
Definition \ref{semi-mod-rel-def} is clear. 

The bicompleteness of $\caC[\cT]$ follows as in the proof of
Theorem \ref{mon-mod-str}. The 2-out-of-3 and retract axioms
for the weak equivalences and the fibrations hold in $\caC[\cT]$
since they hold in $\caC$, the retract axiom for the cofibrations
holds because $\cT I \cof$ is closed under retracts.
So property 2 is fulfilled.

The first half of property 3 is true by the definition of the cofibrations.
By our smallness assumptions we have functorial factorizations of
maps into a cofibration followed by a trivial fibration and into a
map from $\cT J \cell$ followed by a fibration. We claim that a map $f$
in $\cT J \cell$ whose domain is cofibrant in $\caC$ is a trivial
cofibration. $f$ is a weak equivalence by assumption. 
Factor $f$ as $p\circ i$ into a cofibration followed by a trivial fibration.
Since $f$ has the left lifting property with respect to $p$,
$f$ is a retract of $i$ by the retract argument, hence also a cofibration.
Hence we have shown property 4. 

Now let $f$ be a trivial cofibration
whose domain is cofibrant in $\caC$.
We can factor $f$ as $p\circ i$ with $i \in \cT J \cell$ and $p$
a fibration. $p$ is a trivial fibration by the 2-out-of-3 property,
hence $f$ has the left lifitng property with respect to $p$, so $f$
is a retract of $i$ and has therefore the left lifitng property
with respect to fibrations. This is the second half of
property 3. Property 5 immediately follows from the assumptions,
and property 6 is true since limits in $\caC[\cT]$ are computed
in $\caC$.
\end{proof}
{\bf Alternative:}
Assume that $\caC[\cT]$ has coequalizers, that sequential
colimits in $\caC[\cT]$ are computed in $\caC$ and that the pushout
of an object in $\caC[\cT]$ which is cofibrant in $\caC$ by a
map from $\cT I$ (resp. from $\cT J$) is a cofibration
(resp. weak equivalence) as a map in $\caC$. Then the same
conclusion holds as in the Theorem above.
Moreover the conclusion also holds for the
alternative definition of $J$-semi model category without 
the smallness assumptions on the domains of $I$ and $J$ which we made
at the beginning of this section.

\vskip.2cm

\begin{ex}
Let $\Ass(\caC)$ be the category of associative unital algebras
in $\caC$. Then $\Ass(\caC)$ is a $J$-semi model category over $\caC$
(see \cite[Theorem 3.3]{hov-mon}).
\end{ex}

Will will need the
\begin{lem} \label{rmod-tensor}
Let $R$ be a ring with unit in $\caC$, $i$ a map in $(I\otimes R)\cof$
(taken in $R\mmod_r$)
and $j$ a map in $R\mmod$ which is a (trivial) cofibration in $\caC$.
Then $i\Box_R j$ is a (trivial) cofibration in $\caC$.
If $i$ is in $(J\otimes R) \cof$, then $i\Box_R j$ is a trivial
cofibration in $\caC$.
\end{lem}
\begin{proof}
This follows either by \cite{hov-book}[Lemma 4.2.4] applied
to the adjunction of two variables
$R\mmod_r \times R\mmod \rightarrow \caC$, $(M,N) \mapsto M\otimes_R N$,
or by Lemma \ref{cell-pushout-str}.
\end{proof}

\section{Operads} \label{operads}

For a group $G$ write $\caC[G]$ for the category of objects
in $\caC$ together with a right $G$-action. This is the same
as $\eins[G] \mmod_r$, where $\eins[G]$ is the group ring
of $G$ in $\caC$. Let $\caC^\cN$ be the category of sequences
in $\caC$ and $\caC^\Sigma$ the category of symmetric
sequences, i.e. $\caC^\Sigma=\coprod_{n \in \cN} \caC[\Sigma_n]$.
Finally let $\caC^{\cN,\bullet}$
(resp. $\caC^{\Sigma,\bullet}$) be the category of
objects $X$ from $\caC^\cN$ (resp. from $\caC^\Sigma$) together with
a map $\eins \rightarrow X(1)$.

\begin{prop}
For any group $G$ the category $\caC[G]$ has a natural structure
of cofibrantly generated model category with generating
cofibrations $I[G]$ and generating trivial cofibrations $J[G]$.
\end{prop}
\begin{proof}
Easy from \cite{hov-book}[Theorem 2.1.19].
\end{proof}

Hence there are also canonical model structures
on $\caC^\cN$, $\caC^{\cN,\bullet}$,
$\caC^\Sigma$ and $\caC^{\Sigma,\bullet}$.

Note that a map of groups $\varphi:\; H \rightarrow G$
induces a left Quillen functor $\caC[H] \rightarrow \caC[G]$.
If $\varphi$ is injective the right adjoint to this functor
preserves (trivial) cofibrations.

Let $\Op(\caC)$ be the category of operads in $\caC$, where
an operad in $\caC$ is defined as in \cite{km}[Definition 1.1].
Let $F: \caC^\cN \rightarrow \Op(\caC)$ be the functor which
assigns to a sequence $X$ the free operad $FX$ on $X$.
This functor naturally factors through $\caC^{\cN,\bullet}$,
$\caC^\Sigma$ and $\caC^{\Sigma,\bullet}$, and the functors
starting from one of these categories going to $\Op(\caC)$
are also denoted by $F$. The right adjoints of $F$,
i.e. the forgetful functors, map $\caO$ to $\caO^\sharp$.

For any object $A \in \caC$ there is the endomorphism operad
$\End^\Op (A)$ given by $\End^\Op(A)(n) = \uhom(A^{\otimes n},A)$.

We come to the main result of this section:
\begin{thm} \label{op-mod-str}
The category $\Op(\caC)$ is a cofibrantly generated
$J$-semi model category over 
$\caC^{\Sigma,\bullet}$
with generating cofibrations $FI$ and generating trivial cofibrations
$FJ$.
If $\caC$ is left proper (resp. right proper), then $\Op(\caC)$
is left proper relative to $\caC^{\Sigma,\bullet}$ (resp. right proper).
\end{thm}

We first give an explicit description of free operads and
pushouts by free operad maps, which will be needed for the proof
of this Theorem.

%Let $\caT$ denote a set of representatives of isomorphism classes
%of finite connected directed graphs $T$ such that any vertex
%of $T$ has $\le 1$ ingoing arrows. For the graphs we consider
%we do not necessarily require that every arrow has an endpoint.
%Arrows without endpoints are called tails.
%We call the elements of $\caT$ {\em trees}.
%Choose once and for all total orders on all sets of outgoing
%arrows of all vertices of all trees. This corresponds to
%a fixed embedding into the plane of each tree.
%For $i \in \caT$ write $f(t)$ for the number of tails of $t$.
\begin{defi}
\begin{enumerate}
\item An {\em $n$-tree} is a finite connected directed graph $T$
such that any vertex of $T$ has $\le 1$ ingoing arrows,
the outgoing arrows of each vertex $v$ of $T$ are
numbered by $1,\ldots,\val(v)$, where $\val(v)$ is the
number of these arrows, and there are $n$ arrows which do not end
at any vertex, which are called {\em tails} and which are numbered by
$1,\ldots,n$.
By definition the empty tree has one tail, so it is
a $1$-tree.
\item A {\em doubly colored} 
$n$-tree is an $n$-tree together with a decomposition
of the set of vertices into {\em old} and {\em new} vertices.
\item A {\em proper} doubly colored $n$-tree is a doubly colored $n$-tree
such that every arrow starting from an old vertex is either a tail or goes to
a new vertex.
\end{enumerate}
\end{defi}
We denote the set of $n$-trees by $\caT(n)$,
the set of doubly colored $n$-trees by $\caT_\dc(n)$
and the set of proper doubly colored $n$-trees by $\caT_\dc^p(n)$.
Set $\caT:=\coprod_{n \in \cN} \caT(n)$ and
$\caT_\dc^{(p)}:=\coprod_{n \in \cN} \caT_\dc^{(p)}(n)$.

The $n$-trees will describe the $n$-ary operations of free operads,
and indeed $\caT(\bullet)$ is endowed with a natural operad structure
in $\set$.
Let $n,m_1,\ldots,m_n \in \caN$, $m:=\sum_{i=1}^n m_i$ and $T \in \caT(n)$, 
$T_i \in \caT(m_i)$, $i=1,\ldots,n$.
Then the corresponding structure map $\gamma$ of this operad sends
$(T,T_1,\ldots,T_n)$ to the tree which one obtains
from $T$ by glueing the root of $T_i$ to the $i$-th tail of $T$
for every $i=1,\ldots,n$. The previously $j$-th tail of
$T_i$ gets the label $j+ \sum_{k=1}^{i-1} m_k$.
The free right action of $\Sigma_n$ on $\caT(n)$
(which is also defined on $\caT_\dc^{(p)}(n)$)
is such that $\sigma \in \Sigma$ sends a tree $T \in \caT(n)$
to the tree obtained from $T$ by changing
the label $i$ of a tail of $T$ 
into $\sigma^{-1}(i)$. So
$\gamma(T,T_1,\ldots,T_n)^{\sigma(m_{\sigma(1)},\ldots,m_{\sigma(n)})}=
\gamma(T^\sigma,T_{\sigma(1)},\ldots,T_{\sigma(n)})$, where
$\sigma(m_1,\ldots,m_n)$ permutes blocks of lenth $m_i$ in
$1,\ldots,m$ as $\sigma$ permutes $1,\ldots,n$.

Note that an $n$-tree has a natural embedding into the plane and
this embedding is equivalent to the numbering of the arrows.
It follows that there exists a {\em canonical} labelling of the tails of
an $n$-tree, namely the one which labels the tails succesively from
the left to the right in the planar embedding of the tree.

For $T$ an element of $\caT$ or $\caT_\dc^{(p)}$ let $V(T)$ denote the set
of vertices of $T$ (this is defined up to unique isomorphism,
since our trees do not have automorphisms) and let $u(T)$ be
the number of vertices of $T$ of valency $1$ and $U(T)$ be
the set of vertices of $T$ of valency $1$.
For $T \in \caT_\dc^{(p)}$ write
$V_\old(T)$ (resp. $V_\new(T)$) for the set of old (resp. new) vertices
of $T$ and $U_\old(T)$ (resp. $U_\new(T)$) for the set of
old (resp. new) vertices in $U(T)$ and $u_\old(T)$ (resp.
$u_\new(T)$) for their number.

\begin{prop} \label{free-op}
\begin{enumerate}
\item The free operad $FX$ on $X \in \caC^\cN$ is given
by $$(FX)(n)=\coprod_{T \in \caT(n)} \;\;\; \bigotimes_{v \in V(T)}
X(\val(v)) \;\text{.}$$
\item The free operad $FX$ on $X \in \caC^{\cN,\bullet}$
is given by a $\omega$-sequence
$$FX=\colim_{i < \omega} F_iX$$ in $\caC^\cN$, where $(F_iX)_n$
is a pushout of $(F_{i-1}X)_n$ by the map
$$\coprod_{\begin{array}{c}
T \in \caT(n) \\ u(T)=i \end{array}} \left(
\bigotimes_{v \in V(T) \setminus U(T)}
X(\val(v))\right) \otimes e^{\Box (U(T))} \;\text{,}$$
where $e$ is the unit map $\eins \rightarrow X(1)$.
\item The free operad on $X \in \caC^\Sigma$ is given by
$$(FX)(n)=\left(\coprod_{T \in \caT(n)} \;\;\; \bigotimes_{v \in V(T)}
X(\val(v))\right) / \sim \;\text{,}$$
where the equivalence relation $\sim$ identifies for every
isomorphism of directed graphs $\varphi:\; T \rightarrow T'$,
$T,T' \in \caT(n)$, which respects the numbering of the tails
but not necessarily of the arrows, the summands $\bigotimes_{v \in V(T)}
X(\val(v))$ and $\bigotimes_{v \in V(T')}
X(\val(v))$ by the map $\bigotimes_{v \in V(T)} \sigma_v$,
where $\sigma_v: X(\val(v)) \rightarrow X(\val(\varphi(v)))=X(\val(v))$ 
is the action of the element $\sigma_v \in \Sigma_{\val(v)}$ 
such that $\varphi$ maps the $i$-th arrow of $v$ to
the $\sigma_v(i)$-th arrow of $\varphi(v)$.
\item The free operad $FX$ on $X \in \caC^{\Sigma,\bullet}$
is given by a $\omega$-sequence
$$FX=\colim_{i < \omega} F_iX$$ in $\caC^\cN$,
where $(F_iX)_n$
is a pushout of $(F_{i-1}X)_n$ by the map
$$\left(\coprod_{T \in \caT(n), u(T)=i} \left(
\bigotimes_{v \in V(T) \setminus U(T)}
X(\val(v))\right) \otimes e^{\Box (U(T))}\right)
/ \sim \;\text{,}$$
where $e$ is as in 2 and the equivalence relation $\sim$ is like
in 3.
\end{enumerate}
In cases 2 and 4 the attaching map is induced from the operation
of removing a vertex of valency $1$ from a tree. Note that the
morphism in 4 and the attaching morphism respects the
equivalence relation. 
The $\Sigma_n$-actions are induced from the $\Sigma_n$-action
on $\caT(n)$.
\end{prop}
\begin{proof}
We claim that in all four cases the functors $F$ define
a monad the algebras of which are the operads in $\caC$.
So we have to define in all four cases maps $m:\; FFX \rightarrow FX$
and $e:\;X \rightarrow FX$ satisfying the axioms for a monad.
We will restrict ourselves to case i) and leave the other cases
to the interested reader.

The domain of the map $m(n)$ is a coproduct over all
$T \in \caT(n)$, $T_v \in \caT(\val(v))$ for all $v \in V(T)$
over the $$\bigotimes_{v \in V(T), w \in V(T_v)} 
\!\!\!\!\!\!\!\!\! X(\val(w))\;\text{,}$$
and the map $m$ sends such an entry via the identity to the entry
associated to the tree in $\caT(n)$ obtained by replacing every
vertex $v$ of $T$ by the tree $T_v$ in such a way that the
numbering of the arrows starting at $v$ and the numbering of the tails
of $T_v$ correspond. The map $e$ sends $X(n)$ to the summand
$X(n)$ in $FX$ which belongs to the tree with one vertex and
$n$ tails such that the labelling of the arrows coincides with
the labelling of the tails (which are of course all arrows in this case)
(i.e. the labelling of the tails is the canonical one).
It is clear that $m$ is associative and $e$ is a 
two-sided unit. To see that an $F$-algebra is the same as an operad 
one proceeds as follows:
Let $X$ be an $F$-algebra. Let $\caO(n):=X(n)$.
The structure maps of the operad structure we will define on $\caO$ 
are obtained from the algebra map by restricting it to the
summands belonging to trees where every arrow starting at the root
goes to a vertex which has only tails as outgoing vertices
and where the labelling of the tails is the canonical one.
The unit in $\caO(1)$ corresponds to the empty tree.
The right action of a $\sigma \in \Sigma_n$ on $\caO(n)$ is
given by the algebra map restricted to the tree with one vertex and
$n$ tails such that the $i$-th arrow simultaneously is 
the $\sigma^{-1}(i)$-th tail. That $1$ acts as the identity is
the unit property of $X$, and the associativity of the action follows from
the associativity of $X$.
It is easy to see that 
the associativity and symmetry properties of $\caO$ also follow
from the associativity of $X$. The unit properties
follow from the behaviour of the empty tree.

On the other hand let $\caO$ be an operad.
We define an $F$-algebra structure
on $X:=\caO^\sharp$:
Let $T \in \caT(n)$ be a tree with canonical labelling of the
tails. Then it is clear how to define a map from the summand
in $FX$ corresponding to $T$ to $X(n)$ by iterated
application of the structure maps of $\caO$ (the unit of $\caO$
is needed to get the map for the empty tree).
The map on the summand corresponding to $T^\sigma$ for $\sigma \in \Sigma_n$
is the map for $T$ followed by the action of $\sigma$ on $X(n)=\caO(n)$.
One then can check that the associativity, symmetry and unit properties
of the structure maps of $\caO$ imply that we get indeed an $F$-algebra
with structure map $FX \rightarrow X$ just described.
\end{proof}
For describing pushouts by free operad maps we
need an operation which changes a new vertex in a tree
in $\caT_\dc^p(n)$ into an old vertex
and gives again a tree in $\caT_\dc^p(n)$.
This is given by first making the new vertex into
an old vertex to get an element of $\caT_\dc(n)$ and
then removing all arrows joining only old vertices and
identifying the old vertices which have been joined.
The numbering of the arrows of the new tree is most easily described
by noting that this numbering corresponds to a planar embedding of the tree
and the operation of removing the arrows and identifying the vertices
can canonically be done in the plane.
For $T \in \caT_\dc^p$ and $v \in V_\new(T)$ denote
by $\ch_T(v) \in \caT_\dc^p$ the tree obtained
by changing the new vertex $v$ in $T$ into an old vertex.
Note that for $\caO \in \Op(\caC)$ there is a concatenation map
$$\conc_T^\caO(v):\;\;
\caO(\val(v))\otimes \!\!\!\! \bigotimes_{v' \in V_\old(V)} 
\!\!\!\! \caO(\val(v'))
\;\longrightarrow 
\bigotimes_{v' \in V_\old(\ch_T(v))} \!\!\!\!\!\! \caO(\val(v'))$$
induced by applying the operad maps of $\caO$.

\begin{prop} \label{op-pushout}
Let $\caO \in \Op(\caC)$ and $f:\; A \rightarrow B$ 
and $\varphi:\; A \rightarrow \caO^\sharp$ 
be maps in $\caC^\cN$.
Then the pushout $\caO'$ of $\caO$ by $Ff$ with attaching map the
adjoint of $\varphi$ is given by
a $\omega\times (\omega+1)$-sequence 
$\caO'=\colim_{(i,j) < \omega\times(\omega+1)} 
\caO_{(i,j)}$ in $\caC^\cN$, where for $j < \omega$
$\caO_{(i,j)}(n)$ is a pushout of $\caO_{(i,j)-1}(n)$ in $\caC$ 
by the quotient of the map
$$\coprod \left( \bigotimes_{v \in V_\old(T)
\setminus U_\old(T)} \!\!\!\!\!\!\! \caO(\val(v))\right)\otimes
e^{\Box(U_\old(T))}\Box \: 
\underset{v \in V_\new(T)}{\text{\huge $\Box$}}
f(\val(v)) \;\text{,}$$
where the coproduct is over all $T \in \caT_\dc^p(n)$
with $\sharp V_\new(T)=i$ and $u_\old(T)=j$,
with respect to the equivalence relation which identifies
for every
isomorphism of doubly colored
directed graphs $\varphi:\; T \rightarrow T'$,
$T,T' \in \caT_\dc^p$, which respects the labeling of the tails
and of the arrows starting at new vertices, 
the summands corresponding to $T$ and $T'$
via a map analoguous to the map in
Proposition \ref{free-op}.3.
Here $e$ is the unit $\eins \rightarrow \caO(1)$ and
the attaching map is the following:
The domain of the above map is obtained by glueing $i+j$ objects
together, hence we have to give $i+j$ maps compatible with glueing.
The first $i$ maps are induced by removing one of the vertices
in $U_\old(T)$ from $T$, and the other $j$ maps are induced by
changing one of the vertices in $V_\new(T)$ into an old vertex and applying
the maps $\conc_v^\caO(T)$, $v \in V_\new(T)$.
(Note that for $n=1$ the operad $\caO$ appears in the second step
of the limit, in all other cases in the first.)
The $\Sigma_n$-actions are induced from the ones on $\caT_\dc^p(n)$.

There are similar descriptions of pushouts of $\caO$
by free operad maps on maps from $\caC^{\cN,\bullet}$,
$\caC^\Sigma$ and $\caC^{\Sigma,\bullet}$.
\end{prop}
\begin{proof}
Let $\widetilde{\caO}(n)$ be the colimits described in the Proposition.
First of all we check that this is well defined, i.e.
that firstly the $i+j$ maps we have described glue together.
This is the case because the processes of removing old
vertices of valency $1$ and/or changing a new vertex into an
old one and concatenating commute with each other.
Secondly this map factors through the
quotient described in the Proposition because of the symmetry
properties of $\caO$ and because of the fact that in previous steps
quotients with respect to analoguous equivalence relations have been 
taken.

Next we have to equip $\widetilde{\caO} \in \caC^\cN$
with an operad structure.
The unit is the one coming from $\caO$. 
%For $X \in \caC^\cN$ let $X(n,m_1,\ldots,m_n):=X(n)\otimes X(m_1)
%\otimes \cdots \otimes X(m_n)$.
%The map $\caO(n,m_1,\ldots,m_n) \rightarrow \widetilde{\caO}(n,m_1,\ldots,
%m_n)$ is a transfinite composition over $(\omega\times (\omega+1))^{n+1}$.
%The structure map $\gamma:\; \widetilde{\caO}(n,m_1,\ldots,m_n)
%\rightarrow \widetilde{\caO}(m)$ ($m=\sum_{i=1}^n m_i$)
%is def
We define the structure map 
$\gamma:\; \widetilde{\caO}(n)\otimes \widetilde{\caO}(m_1)
\otimes \cdots \otimes \widetilde{\caO}(m_n) \rightarrow
\widetilde{\caO}(m)$ ($m =\sum_{i=1}^n m_i$)
in the following way: 
For $T \in \caT_\dc^p(n)$ let 
$$S(T):=\left(\bigotimes_{v \in V_\old(T)} \caO(\val(v)) \right) \otimes
\!\!\! \bigotimes_{v \in V_\new(T)} \!\!\! B(\val(v)) \;\text{.}$$
First one defines for trees $T \in \caT_\dc^p(n)$,
$T_i \in \caT_\dc^p(m_i)$, $i=1,\ldots,n$,
a map $$\xi_{(T,T_1,\ldots,T_n)}:\;
S(T)\otimes S(T_1)\otimes \cdots \otimes S(T_n)
\rightarrow \widetilde{\caO}(m)\;\text{.}$$ Therefore one glues
the tree $T_i$ to the tail of $T$ with label $i$ and concatenates
such that one gets a tree $\widetilde{T} \in \caT_\dc^p(m)$.
Then by applying structure maps of $\caO$ one gets a map
$S(T)\otimes S(T_1)\otimes \cdots \otimes S(T_n)
\rightarrow S(\widetilde{T})$ and composes this with the canonical map
$S(\widetilde{T}) \rightarrow \widetilde{\caO}(m)$.

Let $m_0:=n$. Suppose we have already defined for a $0\le k \le n$
and for all trees $T_i \in \caT_\dc^p(m_i)$, $i=k,\ldots,n$,
a map $(\bigotimes_{i=0}^{k-1}\widetilde{\caO}(m_i)) \otimes
S(T_k) \otimes \cdots \otimes S(T_n) \rightarrow \widetilde{\caO}(m)$.
From this data one then obtains the same data for $k+1$ instead of $k$
as follows: 
Let $T_i \in \caT_\dc^p(m_i)$, $i=k+1,\ldots,n$.
One defines the map 
$\varphi:\; (\bigotimes_{i=0}^k\widetilde{\caO}(m_i)) \otimes
S(T_{k+1}) \otimes \cdots \otimes S(T_n) \rightarrow \widetilde{\caO}(m)$
by transfinite induction on the terms
of the $\omega \times (\omega+1)$-sequence defining 
$\widetilde{\caO}(m_k)$: 
%by using the data for $k$ to define the
%map on the codomains of the maps along which the pushouts are taken.
%We have to check that this definition is compatible with
%the attaching maps. 
So let $(\bigotimes_{i=0}^{k-1}\widetilde{\caO}(m_i)) \otimes
\caO_{(i,j)}(m_k) \otimes
S(T_k) \otimes \cdots \otimes S(T_n) \rightarrow \widetilde{\caO}(m)$
be already defined and let $\psi$ be the map by which
$\caO_{(i,j+1)}(m_k)$ is a pushout of $\caO_{(i,j)}(m_k)$.
We define $\varphi$ on 
$(\bigotimes_{i=0}^{k-1}\widetilde{\caO}(m_i)) \otimes
\caO_{(i,j+1)}(m_k) \otimes S(T_k) \otimes \cdots \otimes S(T_n)$
by using the data described above for $k$ to get the map 
after taking the appropriate quotient on the
codomain of $\psi$.
One has to check the compatibility of this map with the given map
via the attaching map. To do this for one of
the $i+j$ summands of the domain of $\psi$ one uses the fact
that the same kind of
compatibility is valid in $\widetilde{\caO}(m)$.
Finally when arriving at $k=n$ we get the desired structure map.

The associativity of the structure maps follows by proving the
corresponding statement for the $\xi_{(T,T_1,\ldots,T_n)}$.
This one gets by first glueing trees without concatenating and
then observing that the concatenation processes at different places
commute. The symmetry properties follow in the same way as for free
operads, the unit properties are forced by the fact that in the
$\psi$'s the pushout product over the unit maps is taken.
Hence $\widetilde{\caO}$ is an operad. It receives canonical 
compatible maps in $\Op(\caC)$ from $\caO$ and $FB$.

In the end we have to show that our operad $\widetilde{\caO}$
indeed satisfies the universal property of the pushout by $Ff$.
We need to show that a map $g:\;\caO \rightarrow \caO''$ in
$\Op(\caC)$ together with a map $h:\;B \rightarrow (\caO'')^\sharp$
compatible with the attaching map is the same as a map
$g':\; \widetilde{\caO} \rightarrow \caO''$.
To get $g'$ from $g$ and $h$ one first defines for any $T \in \caT_\dc^p(n)$
a map $S(T) \rightarrow \caO''$ using the structure maps
of $\caO''$. Then one checks that these maps indeed glue
together to a $g'$.
To get $g$ and $h$ from $g'$ one composes $g'$ with
$\caO \rightarrow \widetilde{\caO}$ and $B \rightarrow FB
\rightarrow \widetilde{\caO}$.
These processes are invers to each other.
\end{proof}
\begin{lem} \label{op-pushout-cof}
Let $\caO$, $f$, $\varphi$ and $\caO'$ be as in Proposition \ref{op-pushout},
assume that $\caO \in \Op(\caC)$ 
is cofibrant as object in $\caC^{\Sigma,\bullet}$ and that
$f$ is a (trivial) cofibration. Then the pushout
$\caO \rightarrow \caO'$ is a (trivial) cofibration
in $\caC^{\Sigma,\bullet}$.
There is an analoguous statement for $f$ a (trivial) cofibration
in $\caC^{\cN,\bullet}$, $\caC^\Sigma$ and $\caC^{\Sigma,\bullet}$.
\end{lem}
In the following Lemma we use the fact that if we have a $G$-action
on an object $L$ and a $\Sigma_n$-action on $M$, then
there is a canonical action of the wreath product
$\Sigma_n \ltimes G^n$ on $M\otimes L^{\otimes n}$.
\begin{lem} \label{semi-dir-cof}
Let $n_1,\ldots,n_k \in \cN_{>0}$,
let $G_1,\ldots,G_k$ be groups and
$g_i$ be a cofibration in $\caC[G_i]$, $i=1,\ldots,k$.
Let $f$ be a cofibration in $\caC[\prod_{i=1}^k \Sigma_{n_i}]$.
Then the map $$h:=f \Box \: \text{\LARGE $\Box$}_{i=1}^k 
\: g_i^{\Box n_i}$$ is a cofibration
in $\caC[(\prod_{i=1}^k \Sigma_{n_i}) \ltimes 
(\prod_{i=1}^k G^{n_i})]$. If $f$ or one of the $g_i$ is trivial,
so is $h$.
\end{lem}
\begin{proof}
We restrict to the case $k=1$, the general case is done in the same way.
Set $n:=n_1$, $G:=G_1$ and $g:=g_1$.
We can assume that $g \in I[G]\cell$ and $f \in I[\Sigma_n]\cell$
(or $f \in J[G]\cell$ or $g \in J[\Sigma_n]\cell$).
Let $g:\; L_0 \rightarrow \colim_{i < \mu} L_i$
and $f:\; M_0 \rightarrow \colim_{i < \lambda} M_i$
such that $L_i \rightarrow L_{i+1}$ is a pushout by
$\psi_i \in I[G]$ and $M_i \rightarrow M_{i+1}$ is
a pushout by $\varphi_i \in I[\Sigma_n]$.
Then by Lemma \ref{cell-pushout-str}
$f \Box g^{\Box n}$ is a $\lambda\times \mu^n$-sequence,
and the transition maps are pushouts by
the $\varphi_i\Box \psi_{i_1}\Box
\cdots \Box \psi_{i_n}$, $i<\lambda$; $i_1,\ldots,i_n< \mu$.
We can modify this sequence to make it invariant under
the $\Sigma_n$-action: Let $S$ be the set of unordered sequences
of length $n$ with entries in $\mu$, and for $s \in S$ let $j_s$
be the set of ordered sequences of length $n$ with entries
in $\mu$ which map to $s$. 
Let $s,s' \in S$. In the following let us view
$s$ and $s'$ as monotonly increasing sequences of length $n$.
We say that $s<s'$ if there is a $1\le i < n$ such that
$s(j)=s'(j)$ for $i<j$ and $s(i)<s'(i)$. With this order $S$ is
well-ordered. Now $g^{\Box n}$ is an $S$-sequence
with $s$-th transition map $\psi_s':= \coprod_{w \in j_s} \psi_{w(1)}
\Box \cdots \Box \psi_{w(n)}$, so $f \Box g^{\Box n}$
is the corresponding $\lambda\times S$-sequence with
transition maps the $\varphi_i \Box \psi_s'$, $i < \lambda$, $s \in S$.
Note that on these maps there is a $\Sigma_n \ltimes G^n$-action.
Now to prove our claim it suffices to show that every
$\varphi_i \Box \psi_s'$ is a (trivial) cofibration
in $\caC[\Sigma_n \ltimes G^n]$, which can easily be seen
by noting that every $\varphi_i$ and $\psi_i$ is of the form
$h[G]$ for $h \in I$ (or $h \in J$).
\end{proof}
\begin{proof}[Proof of Lemma \ref{op-pushout-cof}]
Let $\sim$ be the equivalence relation on $\caT_\dc^p$ which
identifies $T$ and $T'$ in $\caT_\dc^p$ if there is an
isomorphism of directed graphs $T \rightarrow T'$ 
which respects the labeling of the arrows starting at new vertices.
Let $C$ be an equivalence class of $\sim$ in $\caT_\dc^p(n)$.
The $\Sigma_n$-action on $\caT_\dc^p(n)$ restricts to a
$\Sigma_n$-action on $C$. We have to show that
the part of the map in Proposition \ref{op-pushout}
given as the appropriate quotient of
$$\coprod_{T \in C} \left( \bigotimes_{v \in V_\old(T)
\setminus U_\old(T)} \caO(\val(v))\right)\otimes
e^{\Box(U_\old(T))}\Box \: 
\underset{v \in V_\new(T)}{\text{\huge $\Box$}}
f(\val(v))\;\;\; \text{(*)}$$
is a (trivial) cofibration in $\caC[\Sigma_n]$.
Let $\Gamma$ be a doubly colored directed graph,
where the arrows starting at new vertices are labelled,
isomorphic to the objects of the same type
underlying the objects from $C$.
Set $$\varphi :=\left( \bigotimes_{v \in V_\old(\Gamma)
\setminus U_\old(\Gamma)} \caO(\val(v))\right)\otimes
e^{\Box(U_\old(\Gamma))}\Box \: 
\underset{v \in V_\new(\Gamma)}{\text{\huge $\Box$}}
f(\val(v))\;\text{.}$$
On $\varphi$ there is an action of $\Aut(\Gamma)$.
Let $t$ be the set of tails of $\Gamma$. There is an action of
$\Aut(\Gamma)$ on $t$.
It is easily seen that the quotient of the map (*) we are considering
is isomorphic to $\varphi\times_{\Aut(\Gamma)} \Sigma_t$.
Hence we are finished if we show that $\varphi$ is a
(trivial) cofibration in $\caC[\Aut(\Gamma)]$.
This is done by induction on the depth of $\Gamma$.
Let $\Gamma_1,\ldots,\Gamma_k$ be the different isomorphism types
of doubly colored directed graphs, such that the arrows starting
at new vertices are labelled,
sitting at the initial vertex
of $\Gamma$ with multiplicities $n_1,\ldots,n_k$ and set
$G_i:=\Aut(\Gamma_i)$, $i=1,\ldots,k$.
Then, if the initial vertex of $\Gamma$ is old,
$\Aut(\Gamma)=(\prod_{i=1}^k \Sigma_{n_i})\ltimes
(\prod_{i=1}^k G_i^{n_i})$, otherwise
$\Aut(\Gamma)=\prod_{i=1}^k G_i^{n_i}$, and the map $\varphi$ is
given like the map $h$ in Lemma \ref{semi-dir-cof}. Now
the claim follows from Lemma \ref{semi-dir-cof} and the induction
hypothesis.
\end{proof}

\begin{proof}[Proof of Theorem \ref{op-mod-str}]
We apply Theorem \ref{mon-semi-mod-str} to the monad $\cT$ which maps
$X$ to $(FX)^\sharp$.
It is known that $\Op(\caC)$
is cocomplete.
Since filtered colimits in $\Op(\caC)$ are computed in $\caC^\cN$,
it follows from Lemma \ref{op-pushout-cof}
that those maps from $FI\cell$ (resp. $FJ\cell)$ whose domain
is cofibrant in $\caC^{\Sigma,\bullet}$ are cofibrations
(resp. trivial cofibrations) in $\caC^{\Sigma,\bullet}$.

It is clear that $\Op(\caC)$ is right proper if $\caC$ is.
If $\caC$ is left proper, then $\caC^{\Sigma,\bullet}$ is left proper,
and the pushout in $\Op(\caC)$ by a cofibration whose domain
is cofibrant in $\caC^{\Sigma,\bullet}$ is a retract of a transfinite
composition of pushouts by cofibrations in $\caC^{\Sigma,\bullet}$,
hence weak equivalences are preserved by these pushouts.
\end{proof}

\begin{rem}
Let $R$ be a commutative ring with unit and $\caC$ be
the symmetric monoidal model category of unbounded
chain complexes of $R$-modules with the projective
model structure. Here the generating trivial cofibrations
are all maps $0 \rightarrow D^n R$, $n \in \cZ$.
The $D^n R$ are clearly null-homotopic.
From this it follows that for a generating trivial cofibration $f$
in $\caC^\cN$
the codomain of the maps
in Proposition \ref{op-pushout} along which the pushouts are taken
(these maps have domain $0$, so the pushouts are trival)
are also null-homotopic by the homotopy which is
on the summand corresponding to a tree $T \in \caT_\dc^p$ the sum over the
homotopies from above over all new vertices of $T$
(this homotopy factors through the quotient
which is taken).
Hence the conditions of Theorem \ref{mon-mod-str} are fulfilled,
so we get a model structure on $\Op(\caC)$ which is the
same as the one provided by \cite[Theorem 6.1.1]{hin1}.
\end{rem}

\begin{rem}
One can use exactly the same methods as above to give the
category of colored operads in $\caC$ for any set of labels the structure of
a $J$-semi model category. In the case of unbounded complexes
over a commutative unital ring as above this $J$-semi model structure
is again a model structure.
\end{rem}

\section{Algebras} \label{algebras}

For an operad $\caO \in \Op(\caC)$ let us denote
by $\Alg(\caO)$ the category of algebras over $\caO$.
Let $F_\caO:\; \caC \rightarrow \Alg(\caO)$ be the
free algebra functor which is given by
$$F_\caO(X)= \coprod_{n \ge 0} \caO(n)\otimes_{\Sigma_n} X^{\otimes n}
\;\text{.}$$
The right adjoint of $F_\caO$ maps $A$ to $A^\sharp$.

\begin{rem} \label{alg-rem}
An $\caO$-algebra structure on an object $A \in \caC$ is the same
as to give a map of operads $\caO \rightarrow \End^\Op(A)$.
\end{rem}

\begin{lem} \label{alg-mor}
Let $I$ be a small category and let $D:\; I \rightarrow \Op(\caC)$,
$i \mapsto \caO_i$, be a functor. Set $\caO:=\colim_{i \in I} \caO_i$
and let $A,B \in \caC$.
Then the following is valid.
\begin{enumerate}
\item To give an $\caO$-algebra structure on $A$ is the same
as to give $\caO_i$-algebra structures on $A$ compatible 
with all transition maps in $D$.
\item Assume that $A$ and $B$ have $\caO$-algebra structures
and let $f:\; A \rightarrow B$ be a map in $\caC$.
Then $f$ is a map of $\caO$-algebras if and only if it is a
map of $\caO_i$-algebras for all $i \in D$.
\end{enumerate}
\end{lem}
\begin{proof}
The first part follows from the Remark above.

Let $f$ be compatible with all $\caO_i$-algebra structures.
Then it can be checked directly that f is also compatible
with the algebra structure on $\caO':=\coprod_{i \in D} \caO_i$.
But since the maps $\caO'(n) \rightarrow \caO(n)$ are
coequalizers in $\caC$ the claim follows.
\end{proof}

The first main result of this section is
\begin{thm} \label{alg-mod-str-1}
Let $\caO \in \Op(\caC)$ be cofibrant.
Then the category $\Alg(\caO)$ is a cofibrantly generated
$J$-semi model category over
$\caC$ with generating cofibrations $F_\caO I$ and generating
trivial cofibrations $F_\caO J$.
If $\caC$ is left proper (resp. right proper), then
$\Alg(\caO)$ is left proper relative to $\caC$ (resp. right proper).
If the monoid axiom holds in $\caC$, then $\Alg(\caO)$ is a
cofibrantly generated model category.
\end{thm}
We want to describe pushouts by free algebra maps.
The following definition has its origin in
\cite{hin2}{Definitions 3.3.1 and 3.3.2].

\begin{defi}
\begin{enumerate}
\item
A {\em doubly colored $\am$-tree} is the same as a doubly colored
$n$-tree except that instead of the labeling of the tails every
tail is marked by either $a$ or $m$.
\item A {\em proper} doubly colored $\am$-tree is
a doubly colored $\am$-tree such that every arrow starting from
an old vertex is either a tail or goes to a new vertex and
every vertex with only tails as outgoing arrows is new and
at least one of the outgoing tails is marked by $m$.
\end{enumerate}
\end{defi}
Note that in particular a proper doubly colored $\am$-tree
has no vertices of valency $0$.

Let $\caT_\am$ be the set of isomorphism classes of doubly
colored $\am$-trees and $\caT_\am^p$ the set of isomorphism classes
of proper doubly colored $\am$-trees.
For $T \in \caT_\am$ let $a(T)$ be the set of tails of $T$ marked
by $a$ and $m(T)$ the set of tails of $T$ marked by $m$.

Let $T \in \caT_\am^p$. Similarly as in the case of operads there
is the operation of changing a new vertex $v$ of $T$ into an old vertex
and also of changing a tail marked by $m$ into a tail marked by $a$.
Denote the resulting trees in $\caT_\am^p$
by $\ch_v(T)$ for $v \in V_\new(T)$
and by $\ch_t(T)$ for $t \in a(T)$.
For $\caO \in \Op(\caC)$ and $A \in \Alg(\caO)$
there is as in the operad case a concatenation map 
$$\conc_T^\caO(v):\;\; \caO(\val(v)) \otimes \!\!\!\!
\bigotimes_{v' \in V_\old(T)}
\!\!\!\! \caO(\val(v')) \otimes A^{\otimes (a(T))} \longrightarrow$$
$$\bigotimes_{v' \in V_\old(\ch_T(v))}
\!\!\!\!\!\!\!\! \caO(\val(v')) \otimes A^{\otimes (a(\ch_T(v)))}$$
induced by the operad maps of $\caO$ and the structure maps
of $A$.
There is also
a concatenation map
$$\conc_t^{\caO,A}(T):\; A \; \otimes \! \bigotimes_{v \in V(T)} 
\!\! \caO(\val(v))
\otimes A^{\otimes (a(T))} \longrightarrow
\!\!\!\! \bigotimes_{v \in V(\ch_T(t))} \!\!\!\!\!\! \caO(\val(v))
\otimes A^{\otimes (a(\ch_T(t)))}$$
induced by the structure maps of the algebra $A$.

\begin{prop} \label{alg-pushout}
Let $\caO \in \Op(\caC)$ and $f:\; X \rightarrow Y$ 
and $\varphi:\; X \rightarrow \caO^\sharp$ 
be maps in $\caC^\cN$. 
Let $\caO'$ be the pushout of $\caO$ by $Ff$ with attaching map 
the adjoint of $\varphi$.
Let $A$ be an $\caO'$-algebra and let $g:\; M \rightarrow N$ and
$\psi:\; M \rightarrow A^\sharp$ be maps in $\caC$. Let $B$ be
the pushout of $A$ as $\caO$-algebra by $F_\caO (g)$ with attaching map
the adjoint of $\psi$ and $B'$ the pushout of $A$
as $\caO'$-algebra by $F_{\caO'} (g)$.
Then the canonical map $h:\; B \rightarrow B'$
is given by a $\omega\times \omega \times (\omega+1)$-sequence
$B'=\colim_{(i,j,k)} B_{(i,j,k)}$,
where for $(i,j,k)$ a successor $B_{(i,j,k)}$ is a pushout of $B_{(i,j,k)-1}$
by the quotient of the map
$$\coprod
\!\left( \bigotimes_{v \in V_\old(T)
\setminus U_\old(T)} \!\!\!\!\!\!\!\!\!\! \caO(\val(v))\right)\otimes
A^{\otimes (a(T))} \otimes e^{\Box(U_\old(T))}
\Box g^{\Box (m(T))} \Box
\: \underset{v \in V_\new(T)}{\text{\huge $\Box$}}
f(\val(v)) \: \text{,}$$
where the coproduct is over all $T \in \caT_\am^p$ with
$\sharp V_\new(T)=i$, $\sharp m(T)=j$ and $u_\old(T)=k$,
with respect to the equivalence relation
which identifies for every isomorphism of directed graphs
$\varphi:\; T \rightarrow T'$, $T,T' \in \caT_\am^p$, which respects
the labeling of the tails and of the arrows which start at new
vertices, the summands corresponding to $T$ and $T'$
by a map which is described on the $\otimes$-part of the summands
involving vertices from $V_\old(T)\setminus U_\old(T)$ as in
Proposition \ref{free-op}.3 and on the other parts by the
identification of the indexing sets via $\varphi$.
The attaching map is induced on the different parts of the domain
of the above map 
by either the operation of removing
a vertex of valency one, by changing a new vertex into an old vertex or
by changing a tail labelled by $m$ into a tail labelled by $a$
and then by applying either a unit map, a map $\conc_T(v)$
or a map $\conc_T(t)$.
\end{prop}
\begin{proof}
We have to do the same steps as in the proof of
Proposition \ref{op-pushout}. Let $C$ be the colimit described
in the Proposition. The attaching maps are again well-defined
because the various concatenation processes commute with each other
and because of the symmetry properties of $\caO$ and the
equivalence relations appearing in previous steps.

We equip $C$ with an $\caO'$-algebra structure: 
Let us define the structure map $\caO'(n) \otimes A^{\otimes n}
\rightarrow A$. For $T \in \caT_\dc^p(n)$ let $S(T)$ be as in the proof
of Proposition \ref{op-pushout}.
For $T \in \caT_\am^p$ let
$$S^a(T):= \left( \bigotimes_{v \in V_\old(T)} \caO(\val(v)) \right)
\otimes A^{\otimes (a(T))} \otimes N^{\otimes (m(T))} \otimes
\!\!\!\! \bigotimes_{v \in V_\new(T)} \!\!\!\! Y(\val(v))\;\text{.}$$
Let $T \in \caT_\dc^p(n)$ and $T_i \in \caT_\am^p$, $i=1,\ldots,n$.
We obtain a tree $\widetilde{T} \in \caT_\am^p$ by glueing
$T_i$ to the tail of $T$ labelled by $i$ and then concatenating.
By applying operad and algebra structure maps we get a map
$S(T) \otimes S^a(T_1) \otimes \cdots \otimes S^a(T_n)
\rightarrow S^a(\widetilde{T})$.
It is then possible by similar considerations as in the proof
of Proposition \ref{op-pushout} to get from these maps the desired
structure map of $C$. It is easy to see that these structure maps
are associative and symmetric. Hence $C$ is an $\caO'$-algebra
which receives an $\caO$-algebra map from $B$ and
$\caO'$-algebra maps from $A$ and $F_{\caO'}(N)$ which are compatible
with each other in the obvious way.

We have to check that for an $\caO'$-algebra $D$ a map
$c:\;C \rightarrow D$ is the same as a map of $\caO'$-algebras
$a:\; A \rightarrow D$ and a map $n:\; N \rightarrow A^\sharp$ which are
compatible with each other.
We get the maps $a$ and $n$ from $c$ by the obvious compositions.
Given $a$ and $n$ we first obtain a map of $\caO$-algebras
$B \rightarrow D$. Moreover for any $T \in \caT_\am^p$ there
is a map $S^a(T) \rightarrow D$ by applying the 
$\caO'$-algebra structure maps of $D$. It is then easy to check
that these maps glue together to give the map $c$.
These processes are invers to each other.
\end{proof}

\begin{lem} \label{alg-pushout-cof}
Let the notation be as in the Proposition above.
If $\caO$ is cofibrant as an object in $\caC^{\Sigma,\bullet}$,
$A$ is cofibrant as an object in $\caC$, $f$ is a cofibration in
$\caC^\cN$ and $g$ is a cofibration in $\caC$ then
the map $h:\; B \rightarrow B'$ is a cofibration in $\caC$. If
$f$ or $g$ is a trivial cofibration then so is $h$.
If $f$ or $g$ is a trivial cofibration and $A$ is arbitrary,
then $h$ lies in $(\caC \otimes J) \cof$, hence is a weak equivalence
if the monoid axiom holds in $\caC$.
\end{lem}
\begin{proof}
Let $\sim$ be the equivalence relation on $\caT_\am^p$ which 
identifies $T$ and $T'$ in $\caT_\am^p$ if there is an
isomorphism of directed graphs $T \rightarrow T'$ 
which respects the labeling of the tails and
of the arrows starting at new vertices.
Let $C$ be an equivalence class of $\sim$ in $\caT_\am^p$.
We have to show that the appropriate quotient of the map
$$\coprod_{T \in C}
\!\left( \bigotimes_{v \in V_\old(T)
\setminus U_\old(T)} \!\!\!\!\!\!\!\!\!\! \caO(\val(v))\right)\otimes
A^{\otimes (a(T))} \otimes e^{\Box(U_\old(T))}
\Box g^{\Box (m(T))} \Box
\: \underset{v \in V_\new(T)}{\text{\huge $\Box$}}
f(\val(v)) $$
is a (trivial) cofibration in $\caC$ (or lies in $(\caC \otimes J) \cof$
under the assumptions of the last statement).
This is done as in
the proof of Lemma \ref{op-pushout-cof} by induction
on the depth of the trees in $C$. This time instead
of using Lemma \ref{semi-dir-cof} it is sufficient
to use Lemma \ref{rmod-tensor} applied to rings of
the form $\eins[\prod_{i=1}^k \Sigma_{n_i}]$.
\end{proof}
\begin{proof}[Proof of Theorem \ref{alg-mod-str-1}]
We apply Theorem \ref{mon-semi-mod-str} to the monad $\cT_\caO$ which
maps $X$ to $(F_\caO X)^\sharp$.
It is known that $\Alg(\caO)$ is cocomplete.
Since filtered colimits in $\Alg(\caO)$ are computed
in $\caC$ we are reduced to show that
the pushout of an $\caO$-algebra $A$ which is cofibrant 
as an object in $\caC$ by a map in $F_\caO I$ (resp.
in $F_\caO J$) is a cofibration (resp. trivial cofibration)
in $\caC$. Since $\caO$ is a retract of a cell operad
(i.e. a cell complex in $\Op(\caC)$)
such a pushout is a retract of a pushout of the same kind
with the additional hypothesis that $\caO$ is a cell
operad. So let $\caO$ be a cell operad.
Then the pushout in question is a transfinite composition
of maps $h$ as in Proposition \ref{alg-pushout}, hence
by Lemma \ref{alg-pushout-cof} it is a (trivial) cofibration.

It is clear that $\Alg(\caO)$ is right proper if $\caC$ is.
The pushout in $\Alg(\caO)$ by a cofibration whose domain
is cofibrant in $\caC$ is a retract of a transfinite
composition of pushouts by cofibrations in $\caC$,
hence if $\caC$ is left proper
weak equivalences are preserved by these pushouts,
so $\Alg(\caC)$ is also left proper.

The last statement follows again from Lemma \ref{alg-pushout-cof}.
\end{proof}

The second result concerning algebras is
\begin{thm} \label{alg-mod-str-2}
Let $\caO$ be an operad in $\caC$ which is cofibrant as an
object in $\caC^\Sigma$. Then $\Alg(\caO)$ is a cofibrantly generated
$J$-semi model category with generating cofibrations
$F_\caO I$ and generating trivial cofibrations $F_\caO J$.
If $\caC$ is right proper, so is $\Alg(\caO)$.
\end{thm}

The next result enables one to control
pushouts of cofibrant algebras by free algebra maps.

For an ordinal $\lambda$ denote by $S_\lambda$ the set
of all maps $f: \; \lambda \rightarrow \frac{1}{2}\cN$
such that $f(i)$ is $\ne 0$ only for finitely many $i < \lambda$,
if $f(i) \notin \cN$ then $i > 0$ and $f(i')=0$ for all $i'< i$
and if $\lambda$ is a successor then $f(\lambda -1)=0$.
For $f, f' \in S_\lambda$ say that $f<f'$ if there is
an $i< \lambda$ such that $f(i')=f'(i')$ for all $i'> i$
and $f(i)<f'(i)$. With this ordering $S_\lambda$ is
well-ordered.
For $i < \lambda$ denote by $f_i$ the element of $S_\lambda$ with
$f_i(i)=\frac{1}{2}$ and $f_i(i')=0$ for $i' \ne i$.
Set $S_{\lambda,+}:=S_\lambda \sqcup \{*\}$, where $*$ is
by definition smaller than any other element in $S_{\lambda,+}$.
Note that $f \in S_{\lambda,+}$ is a successor if and only if $f\ne *$
and $f(\lambda) \subset \cN$.
For $f \in S_{\lambda,+}$ a successor let 
$|f|:=\sum_{i < \lambda} f(i) \in \cN$ and
$\Sigma_f := \prod_{i < \lambda} \Sigma_{f(i)}$.

\begin{prop} \label{cell-alg-str}
Let $\caO \in Op(\caC)$ and $A=\colim_{i < \lambda} A_i$ be
a $F_\caO(\Mor(\caC))\cell$ $\caO$-algebra ($\Mor(\caC)$ is
the class of all morphisms in $\caC$) with $A_0=\caO(0)$,
where the transition maps $A_i \rightarrow A_{i+1}$ are
pushouts of free $\caO$-algebra maps on maps $g_i:\; K_i \rightarrow L_i$
in $\caC$ by maps adjoint to 
$\varphi_i: \; K_i \rightarrow A_i^\sharp$.
Then $A$ is a transfinite composition $A=\colim_{f \in S_{\lambda,+}}
A_f$ in $\caC$ such that
\begin{enumerate}
\item $A_*=0$ and $A_{f_i}=A_i$ for $i<\lambda$,
\item for $f \in S_\lambda$ such that for an $i_0 < \lambda$
we have $f(i_0) \notin \cN$, there is for all $m \in \cN$,
successors $l \in S_{\lambda,+}$ with $l < f$
and $n:=m + |l|$ a map
$$\Psi_{f,m,l}:\; \caO(n)\otimes_{(\Sigma_m \times
\Sigma_l)} \left( A_{i_0}^{\otimes m}
\otimes \bigotimes_{i<\lambda} L_i^{\otimes l(i)} \right)
\rightarrow A_f$$ compatible with the structure map
$\caO(n)\otimes_{\Sigma_n} A^{\otimes n} \rightarrow A$.
By applying permutations to $\caO(n)$ and the big bracket there
are similar maps for other orders of the factors in the big bracket.
These maps satisfy the following conditions:
\begin{enumerate}
\item They are compatible with the maps $L_i \rightarrow A_{i_0}$
for $i < i_0$. Moreover, if we replace a factor $L_{i_0}$ by
$K_{i_0}$ we can either go to $L_{i_0}$ or to $A_{i_0}$ and
apply suitable maps $\Psi$. Then the two compositions coincide.
\item They are associative in the following sense:
Let $f_1,\ldots,f_k \in S_\lambda$ be limit elements
with $f_i < f$, $i=1,\ldots,k$,
and let for each $f_i$ be given $m_i$, $l_i$ and $n_i$ 
satisfying the same conditions as $m$, $l$ and $n$
for $f$. Let $D_i$ be the domain of $\Psi_{f_i,m_i,l_i}$.
Then the two possible ways to get from
$$\caO(n) \otimes \left(\bigotimes_{i=1}^k D_i \right)
\otimes A_{i_0}^{\otimes m} \otimes \bigotimes_{i < \lambda}
L_i^{\otimes l(i)}$$
to $A_f$ given by either applying the $\Psi_{f_i,m_i,l_i}$ and
then $\Psi_{f,m+k,l}$ or by applying the obvious operad structure maps
and a suitable permutation of 
$\Psi_{f,m+\sum_{i=1}^k m_i,\; l + \sum_{i=1}^k l_i}$ coincide.
\end{enumerate}
\item For any successor
$f \in S_{\lambda,+}$ the map $A_{f-1} \rightarrow A_f$ is a pushout
by $$\caO(|f|)\otimes_{\Sigma_f}
\text{\LARGE $\Box$}_{i < \lambda} \; g_i^{\Box f(i)} \;\text{,}$$
where the attaching maps on the various parts of the domain of
this map are induced from the maps in (2) (see below).
\end{enumerate}
\end{prop}
\begin{proof}
The whole Proposition is shown by induction on $\lambda$,
so suppose that it is true for ordinals less than $\lambda$.
We construct the map in 2, prove its properties
and define the attaching map in 3 by transfinite
induction: Suppose $f \in S_{\lambda,+}$ is a successor, that
$A_{f'}$ is defined for $f' < f$ and that the map in
2 is defined for all limit elements $\tilde{f} \in S_\lambda$
with $\tilde{f} < f$.
Let $i_0 \in \lambda$ with $f(i_0) > 0$ and let 
$f'$ coincide with $f$ except that  $f'(i_0)=f(i_0)-1$.
The attaching map on the summand
$$S:=\caO(|f|) \otimes_{\Sigma_{f'}}
\left(\left(\bigotimes_{i < i_0} L_i^{\otimes f(i)} \right)
\otimes K_{i_0} \otimes L_{i_0}^{\otimes (f(i_0) -1)}
\otimes \bigotimes_{i_0 < i < \lambda} L_i^{\otimes f(i)}\right)$$
of the domain of 
$$\caO(|f|)\otimes_{\Sigma_f}
\text{\LARGE $\Box$}_{i < \lambda} \; g_i^{\Box f(i)}$$
is given as follows: Let $\tilde{f},l \in S_\lambda$ be defined
by $\tilde{f}(i_0)= f(i_0) - \frac{1}{2}$, $l(i_0)= f(i_0)-1$,
$\tilde{f}(i)=l(i)=0$ for $i< i_0$ and $\tilde{f}(i)=l(i)=f(i)$
for $i > i_0$.
Let $m:= 1+ \sum_{i < i_0} f(i)$.
There is a canonical map
$$S \rightarrow \caO(|f|)\otimes_{\Sigma_{f'}} 
\left(A_{i_0}^{\otimes (m-1)}
\otimes A_{i_0} \otimes L_{i_0}^{\otimes (f(i_0) -1)}
\otimes \bigotimes_{i_0 < i < \lambda} L_i^{\otimes f(i)}\right)$$
whose codomain maps naturally to the domain of
$\Psi_{\tilde{f},m,l}$.
So we get maps $S \rightarrow A_{\tilde{f}}
\rightarrow A_{f-1}$ the composition of which is the
attaching map on the summand $S$.
These maps glue together for various summands $S$:
There are two cases to distinguish. In the first one
the intersection of two summands contains $K_{i_0}$ twice. Then
the two maps on this intersection coincide
because of the symmetric group invariance.
In the second case the intersection $I$ contains $K_{i_0'}$ 
and $K_{i_0}$ with $i_0' < i_0$.
Let $\tilde{f}$ be as above and $\tilde{f}'$ be similarly defined
for $i_0'$. 
Now the two properties 2(a) of the maps $\Psi$
state that both maps $I \rightarrow A_f$ are equal the
map induced by first mapping both $K_{i_0'}$ and $K_{i_0}$
to $A_{i_0}$ and then applying a suitable map $\Psi$.

Now suppose $f \in S_\lambda$ is a limit element with
$f(i_0) \notin \cN$ for some $i_0 < \lambda$. 
Define $A_f$ as the colimit of the preceeding $A_{f'}$, $f' < f$.
Let $m$, $l$ and $n$ be as in 2.
We define $\Psi_{f,m,l}$ by induction on $m$
and on $S_{i_0}$ 
using the fact that $A_{i_0}=\colim_{f' \in S_{i_0}} A_{f'}$
by induction hypothesis for the induction on $\lambda$.
For abbreviation set $\caL:=\bigotimes_{i < \lambda} L_i^{\otimes l(i)}$.
Let $f' \in S_{i_0}$ be a successor and let a map
$$\psi_{f' -1}:\; \caO(n) \otimes \left(A_{i_0}^{\otimes (m-1)} \otimes
A_{f' -1} \otimes \caL
\right) \rightarrow A_f$$
be already defined. $A_{f'}$ is a pushout of $A_{f' -1}$ by
$$\varphi:\; \caO(|f'|)\otimes_{\Sigma_{f'}}
\text{\LARGE $\Box$}_{i < i_0} \; g_i^{\Box f'(i)} \;\text{.}$$
Let $C:= \bigotimes_{i < i_0} L_i^{\otimes f'(i)}$.
Then the codomain of $\varphi$ is $\caO(|f'|) \otimes_{\Sigma_{f'}} C$.
Moreover by induction hypothesis for the $m$-induction
there is a map $$\caO(n+|f'|-1)\otimes \left(A_{i_0}^{\otimes (m-1)} \otimes C
\otimes \caL
\right) \rightarrow A_f\;\text{,}$$
hence by plugging in $\caO(|f'|)$ into the $m$-th place of
$\caO(n)$ we get a map
$$\caO(n)\otimes \left(A_{i_0}^{\otimes (m-1)} \otimes \caO(|f'|) \otimes C
\otimes \caL
\right) \rightarrow A_f\;\text{.}$$
This map and $\psi_{f' -1}$ glue together to a map
$\psi_{f'}$: We have to show that they coincide after composition
on domains of the form
$$\caO(n)\otimes A_{i_0}^{\otimes (m-1)} \otimes \caO(|f'|) 
\otimes_{\Sigma_{f''}} S'
\otimes \caL$$
for $\caO(|f'|)\otimes_{\Sigma_{f''}} S'$ 
a summand of the domain of $\varphi$
containing $K_{i_0'}$ for some $i_0' < i_0$ (the definition of
$f''$ is similar to the one of $f'$).
To do this we can restrict for every $A_{i_0}$ to objects 
$\caO(|f_i'|) \otimes_{\Sigma_{f_i'}} C_i$, $i=1,\ldots,m-1$, for $C_i$
of the same shape as $C$ and $f_i' \in S_{i_0,+}$ successors.
Then the two possible ways to get from
$$\caO(n)\otimes \left(\bigotimes_{i=1}^{m-1} \caO(|f_i'|) 
\otimes_{\Sigma_{f_i'}} C_i \right) \otimes \caO(|f'|)
\otimes_{\Sigma_{f''}} S' \otimes \caL$$
to $A_f$ can be compared by mapping $K_{i_0'}$ to $A_{i_0'}$,
unwrapping the definitions of $A_f$ and $\Psi_{f' -1}$
and using associativity of $\caO$.
We arrive at a map $\caO(n) \otimes A_{i_0}^{\otimes m} \otimes \caL
\rightarrow A_f$. That it
factors through the $(\Sigma_m \times \Sigma_l)$-quotient follows
after replacing
$A_{i_0}^{\otimes m}$ by $\left(\bigoplus_{i=1}^k 
\caO(|f_i'|) \otimes_{\Sigma_{f_i'}} C_i \right)^{\otimes m}$
(the $C_i$ and $f_i'$ as above) in the domain of this map
since then the $(\Sigma_m \times \Sigma_l)$-relation is obviously
also valid in $A_f$.

Both properties 2(a) and (b) follow easily by the technique
of restricting any appearing $A_i$ by a factor 
$\caO(|f'|) \otimes_{\Sigma_{f'}} C$.

Now using the maps $\Psi$ and property 2(b) we can equip
$\widetilde{A}:=\colim_{f \in S_{\lambda,+}}$ with an $\caO$-algebra structure
(to do this accurately we have to enlarge $\lambda$ a bit
and the corresponding sequence by trivial pushouts).

We are left to prove the universal property for $\widetilde{A}$
by transfinite induction on $\lambda$.
So let it be true for ordinals less than $\lambda$.
If $\lambda$ is a limit ordinal or the successor
of a limit ordinal there is nothing to show.
Let $\lambda=\alpha + 2$, let $B$ be an $\caO$-algebra
and $A_{\alpha} \rightarrow B$ a map in $\Alg(\caO)$ and
$L_\alpha \rightarrow B^\sharp$ a map in $\caC$ such that these
two maps are compatible via the attaching map.
We define maps $A_f \rightarrow B$ by transfinite induction
on $S_{\lambda,+}$, starting with the given map on
$A_{f_\alpha}= A_\alpha$. So let $f_\alpha < f < \lambda$ be a successor.
Since for any $i \le \alpha$ there is a map $L_i \rightarrow B^\sharp$
we have a natural map $$\caO(|f|)\otimes_{\Sigma_f}
\bigotimes_{i < \lambda} \; L_i^{\otimes f(i)}
\rightarrow B$$ using the algebra structure maps of $B$.
We have to show that this is compatible via the attaching map
from the domain $D$ of $\caO(|f|)\otimes_{\Sigma_f}
\text{\LARGE $\Box$}_{i < \lambda} \; g_i^{\Box f(i)}$
to $A_{f-1}$ with the map $A_{f-1} \rightarrow B$ coming from
the induction hypothesis. We check this again on a summand
$S$ of $D$ containing some $K_{i_0}$. The attaching map
on $S$ is induced from $\Psi_{\tilde{f},m,l}$ as above. The
canonical map from the domain of $\Psi_{\tilde{f},m,l}$ to $B$
is compatible with $A_{\tilde{f}} \rightarrow B$ (as one checks
again by replacing any $A_{i_0}$ by essentially products of $L_i$'s
as above), which together with the fact that $L_{i_0} \rightarrow B$
and $A_{i_0+1} \rightarrow B$ coincide on $K_{i_0}$ implies the
compatibility.
By construction and the definition of the algebra structure on
$A_{\alpha + 1}$ the map $A_{\alpha +1} \rightarrow B$ just defined
is an $\caO$-algebra map.

If we have on the other hand a map of $\caO$-algebras
$A_{\alpha +1} \rightarrow B$ we can restrict it to get compatible maps
$A_\alpha \rightarrow B$ and $L_\alpha \rightarrow B^\sharp$.
These two assignments are inverse to each other.
\end{proof}

\begin{proof}[Proof of Theorem \ref{alg-mod-str-2}]
Let $\caO \in \Op(\caC)$ be cofibrant in $\caC^\Sigma$.
We have to show that the pushout of an $\caO$-algebra such that
the map from the initial $\caO$-algebra to $A$ is in $F_\caO I \cof$
by a map from $F_\caO I$ (resp. $F_\caO J$) is a cofibration
(resp. trivial cofibration) in $\caC$. We can assume that
$A$ is a $F_\caO I\cell$ $\caO$-algebra, since in the general
situation all maps we look at are retracts of corresponding
maps in this situation.
But if $A$ is a cell $\caO$-algebra our claim
immediately follows from Proposition \ref{cell-alg-str}
and Lemma \ref{rmod-tensor}.
\end{proof}

\section{Module structures} \label{module-structures}

In this section we want to show that if $\caC$
is simplicial $\Alg(\caO)$ is also a simplicial $J$-semi
model category in the cases when the assumptions
of Proposition \ref{alg-mod-str-1} or Proposition \ref{alg-mod-str-2}
are fulfilled. Also $\Op(\caC)$ is simplicial if $\caC$ is.

\begin{defi}
Let $\caD$ and $\caE$ be $J$-semi model categories (maybe over $\caC$)
and let $\caS$ be a model category.
Then a {\em Quillen bifunctor} $\caD \times \caS \rightarrow \caE$
is an adjunction of two variables $\caD \times \caS \rightarrow \caE$
such that for any cofibration $g:\; K \rightarrow L$
in  $\caS$ and fibration $p:\; Y \rightarrow Z$ in $\caE$,
the induced map
$$\Hom_{r,\Box}(g,p):\; \Hom_r(L,Y) \rightarrow
\Hom_r(L,Z)\times_{\Hom_r(K,Z)} \Hom_r(K,Y)$$
is a fibration in $\caD$ which is trivial if $g$ or $p$ is.

(See also \cite{hov-book}[Lemma 4.2.2].)
\end{defi}

It follows that for $f$ a cofibration in $\caD$ and $g$ a cofibration 
in $\caS$ both of which have cofibrant domains the pushout
$f \Box g$ is a cofibration in $\caE$ which is trivial if $f$
or $g$ is.

\begin{defi}
Let $\caD$ be a $J$-semi model category (maybe over $\caC$)
and let $\caS$ be a symmetric monoidal model category.
Then a {\em Quillen $\caS$-module structure} on $\caD$ is 
a $\caS$-module structure on $\caD$ such that the
action map $\otimes:\; \caD \times \caS \rightarrow \caD$ is a Quillen
bifunctor and the map $X \otimes (QS) \rightarrow X \otimes S \cong X$
is a weak equivalence for all cofibrant $X \in \caD$, where
$QS \rightarrow S$ is a cofibrant replacement.
\end{defi}

If $\caD$ has a Quillen $\caS$-module structure we say
that $\caD$ is an $\caS$-module.

Let now $\caS$ be
a symmetric monoidal model category where the tensor product
is the categorical product on $\caS$, so let us denote this by $\times$
(e.g. $\caS=\sset$).
Let be given a symmetric monoidal left Quillen functor
$\caS \rightarrow \caC$.

\begin{prop} \label{alg-enrich}
Let the situation be as above and assume that either $\eins$ 
is cofibrant in $\caS$ or that $\caC$ is left proper and the
maps in $I$ have cofibrant domains.
Let $\caO$ be an operad in $\caC$ which is either cofibrant
in $Op(\caC)$ or cofibrant as an object in $\caC^\Sigma$.
Then the $J$-semi model category (in the first case over $\caC$)
$\Alg(\caO)$ is naturally an $\caS$-module and the functor
$\caC \rightarrow \Alg(\caC)$ is an $\caS$-module homomorphism.
\end{prop}
\begin{proof}
Let $A \in \caC$ and $K \in \caS$. We denote by
$A^K$ the homomorphism object $\uhom(K,A) \in \caC$.
There is a map of operads $$\End^\Op(A) \rightarrow \End^\Op(A^K)\;\text{,}$$
which is described on as follows:
We give the maps 
$$\uhom(A^{\otimes n},A) \rightarrow \uhom((A^K)^{\otimes n},A^K)$$ 
on $T$-valued points ($T \in \caC$):
A map $T \otimes A^{\otimes n} \rightarrow A$ is sent
to the composition $$T \otimes (A^K)^{\otimes n} \rightarrow
T \otimes (A^{\otimes n})^{K^n} \rightarrow T \otimes (A^{\otimes n})^K
\rightarrow A^K\;\text{,}$$ where the second map is induced by the
diagonal $K \rightarrow K^n$.

Hence for objects $K \in \caS$ and $A \in \Alg(\caO)$
the object $(A^\sharp)^K$ has a natural structure of $\caO$-algebra
given by the composition
$\caO \rightarrow \End^\Op(A) \rightarrow \End^\Op(A^K)$.
We denote this $\caO$-algebra by $A^K$.

For a fixed $K \in \caS$ the functor $\Alg(\caO) \rightarrow \Alg(\caO)$,
$A \mapsto A^K$, has a left adjoint $A \mapsto A\otimes K$,
which is given for a free $\caO$-algebra $F_\caO(X)$, $X \in \caC$,
by $F_\caO(X)\otimes K = F_\caO(X \otimes K)$ and
which is defined in general by be requirement that
$\_ \otimes K$ respects coequalizers (note that every
$\caO$-algebra is a coequalizer of a diagram where only
free $\caO$-algebras appear).
So we have a functor $\Alg(\caO) \times \caS \rightarrow \Alg(\caO)$.

Let now $B \in \Alg(\caO)$ be fixed.
By a similar argument as above the functor $\caS^\op \rightarrow \Alg(\caO)$,
$K \mapsto B^K$, has a left adjoint $A \mapsto \uhom^\caS(A,B)$,
which sends a free $\caO$-algebra $F_\caO(X)$, $X \in \caC$,
to the image of $\uhom(X,B^\sharp)$ in $\caS$.

One checks that the functor $\Alg(\caO) \times \caS \rightarrow \Alg(\caO)$
we constructed defines an action of $\caS$ on $\Alg(\caO)$.

It remains to show that this functor is a Quillen bifunctor
and that the unit property is fulfilled.
So let $g:\; K \rightarrow L$ be a cofibration in $\caS$
and $p: \; Y \rightarrow Z$ a fibration in $\Alg(\caO)$.
We have to show that $\Hom_{\Box,r}(g,p)$ is a fibration in $\Alg(\caO)$,
i.e. lies in $F_\caO J \inj$.
By adjointness this means that $p$ has the right lifting
property with respect to the maps $(F_\caO f) \Box g=F_\caO(f \Box g)$
for all $f \in J$, which is by adjointness the case because
$f \Box g$ is a 
trivial cofibration. When $p$ or $f$ is trivial we want to show
that $\Hom_{\Box,r}(g,p)$ lies in $F_\caO I \inj$, so $p$ should
have the right lifting property with respect to the maps
$F_\caO(f \Box g)$ for all $f \in I$, which is again the case by adjointness.

If $\eins$ is cofibrant in $\caS$ we are ready.
In the other case the unit property follows by transfinite induction
from the explicit description of algebra pushouts, and hence
the structure of cell algebras, 
given in Proposition \ref{alg-pushout} and the structure of
cell algebras given in Proposition \ref{cell-alg-str}.
\end{proof}

In a similar manner one shows

\begin{prop} \label{op-enrich}
Let the situation be as before Proposition \ref{alg-enrich}
and assume that either $\eins$ 
is cofibrant in $\caS$ or that $\caC$ is left proper and the
maps in $I$ have cofibrant domains.
Then $\Op(\caC)$ is naturally an $\caS$-module and the functor
$\caC \rightarrow \Op(\caC)$ is an $\caS$-module homomorphism.
\end{prop}

\section{Modules} \label{modules}

Let $\caO \in \Op(\caC)$ and $A \in \Alg(\caO)$.
We denote the category of $A$-modules by $(\caO,A)\mmod$,
or $A\mmod$ if no confusion is likely.
Let $F_{(\caO,A)}:\; \caC \rightarrow
A\mmod$ (or $F_A$ for short) be the free $A$-module
functor. It is given by $M \mapsto U_\caO(A)\otimes M$,
where $U_\caO(A)$ is the universal enveloping algebra of 
the $\caO$-algebra $A$.
%Let us call an $\caO$-algebra $A$ {\em cofibrant}
%if the map from the initial $\caO$-algebra to $A$ lies in
%$F_\caO I \cof$.
Recall that $\Ass(\caC)$ denotes the category of associative unital
algebras in $\caC$, and let $F_\Ass$ be the free associative algebra functor
$\caC \rightarrow \Ass(\caC)$.

The main result of this section is
\begin{thm} \label{mod-mod-str}
Let $\caO \in \Op(\caC)$ and $A \in \Alg(\caO)$.
Let one of the following two conditions be satisfied:
\begin{enumerate}
\item $\caO$ is cofibrant as an object in $\caC^\Sigma$ and $A$
is a cofibrant $\caO$-algebra.
\item $\caO$ is cofibrant in $\Op(\caC)$ and $A$ is cofibrant
as an object in $\caC$.
\end{enumerate}
Then there is cofibrantly generated model structure on
$A\mmod$ with generating cofibrations $F_A I$ and generating
trivial cofibrations $F_A J$.
There is a right $\caC$-module structure on $A\mmod$.
\end{thm}

This theorem will follow from the fact that in each of the two
cases the enveloping algebra $U_\caO(A)$ is cofibrant in $\caC$,
since $A \mmod$ is canonically equivalent to
$U_\caO(A)\mmod$.

Note that there is a canonical surjection
from the tensor algebra to the universal enveloping algebra
$$T_\caO(A):=\coprod_{n \in \cN} \caO(n+1)\otimes_{\Sigma_n} A^{\otimes n}
\rightarrow U_\caO(A)\;\text{.}$$

\begin{prop}
Let $\caO \in \Op(\caC)$ and $f:\; X \rightarrow Y$ 
and $\varphi:\; X \rightarrow \caO^\sharp$ 
be maps in $\caC^\cN$. 
Let $\caO'$ be the pushout of $\caO$ by $f$ with attaching map 
the adjoint of $\varphi$.
Let $A$ be a $\caO'$-algebra.
Then $U_{\caO'}(A)$ is a pushout of $U_\caO(A)$
in $\Ass(\caC)$ by the map
$F_\Ass(\coprod_{n \in \cN} f(n)\otimes A^{\otimes (n-1)})$ with attaching map
the adjoint to the
composition $\coprod_{n \in \cN} X(n+1) \otimes A^{\otimes n}
\rightarrow \coprod_{n \in \cN} \caO(n+1) \otimes_{\Sigma_n} A^{\otimes n}
\rightarrow U_\caO(A)$.
\end{prop}
\begin{proof} (Compare to \cite{hin1}[6.8.1. Lemma.])
A $(\caO',A)$-module structure on a $(\caO,A)$-module
$M$ is given by maps $Y(n+1)\otimes A^{\otimes n} \otimes M \rightarrow M$
for $n \in \cN$ such that the compositions with $f(n+1) \otimes 
A^{\otimes n} \otimes M$ equals the composition
$X(n+1)\otimes A^{\otimes n} \otimes M \rightarrow
\caO(n+1)\otimes A^{\otimes n} \otimes M \rightarrow M$.
The same statement is true for a module structure under the
described pushout algebra on a $U_\caO(A)$-module.
\end{proof}

\begin{cor}
Let $\caO \in \Op(\caC)$ be cofibrant and let $A$ be an $\caO$-algebra
which is cofibrant as an abject in $\caC$. Then
$U_\caO(A)$ is cofibrant in $\Ass(\caC)$, in particular is cofibrant
as an object in $\caC$.
\end{cor}
Hence the first part of Theorem \ref{mod-mod-str} is proven.

\begin{cor} \label{univ-vergl-1}
Let $\caC$ be left proper,
let $\caO \in \Op(\caC)$ be cofibrant and let $A \rightarrow A'$
be a weak equivalence between $\caO$-algebras both of which are
cofibrant as objects in $\caC$. Then the map
$U_\caO(A) \rightarrow U_\caO(A')$ is a weak equivalence.
\end{cor}

We have an analoguous result to Proposition \ref{cell-alg-str}
for the enveloping algebra of a cell algebra.

\begin{prop} \label{cell-univ-alg-str}
Let $\caO \in Op(\caC)$ and $A=\colim_{i < \lambda} A_i$ be
a $F_\caO(\Mor(\caC))\cell$ $\caO$-algebra with $A_0=\caO(0)$,
where the transition maps $A_i \rightarrow A_{i+1}$ are
pushouts of free $\caO$-algebra maps on maps $g_i:\; K_i \rightarrow L_i$
in $\caC$ by maps adjoint to 
$\varphi_i: \; K_i \rightarrow A_i^\sharp$.
Then $U:=U_\caO(A)$ is a transfinite composition 
$U=\colim_{f \in S_{\lambda,+}} U_f$ in $\caC$ such that
\begin{enumerate}
\item $U_*=0$ and $U_{f_i}=U_\caO(A_i)$ for $i<\lambda$,
\item for $f \in S_\lambda$ such that for an $i_0 < \lambda$
we have $f(i_0) \notin \cN$, there is for all $m \in \cN$,
successors $l \in S_{\lambda,+}$ with $l < f$
and $n:=m+ |l|$ a map
$$\caO(n+1)\otimes_{(\Sigma_m \times
\Sigma_l)} \left( A_{i_0}^{\otimes m}
\otimes \bigotimes_{i<\lambda} L_i^{\otimes l(i)} \right)
\rightarrow U_f$$ compatible with the map
$\caO(n+1)\otimes_{\Sigma_n} A^{\otimes n} \rightarrow U$ and
\item
for any successor
$f \in S_{\lambda,+}$ the map $U_{f-1} \rightarrow U_f$ is a pushout
by $$\caO(|f|+1)\otimes_{\Sigma_f}
\text{\LARGE $\Box$}_{i < \lambda} \; g_i^{\Box f(i)} \;\text{,}$$
where the attaching maps on the various parts of the domain of
this map are induced from the maps in (2).
\end{enumerate}
\end{prop}
\begin{proof}
This Proposition is proven in essentially the same way as
Proposition \ref{cell-alg-str} except that this time we
have to define associative algebra structures on
the $U_{f_i}$ and to verify the universal property stating
the equivalence of module categories.
For the associative algebra structure one uses the same
formulas as for the tensor algebra and checks that they are compatible
with the attaching maps.
For the universal property one uses the fact that an $A$-module $M$ is
given by maps $$\caO(|f|+1) \otimes_{\Sigma_f} 
\left(\bigotimes_{i < \lambda} L_i^{\otimes f(i)} \right) \otimes M
\rightarrow M$$ which are compatible in various ways the explicit
formulation of which we leave to the reader.
\end{proof}

\begin{cor}
For $\caO$ an operad in $\caC$ which is cofibrant in $\caC^\Sigma$
and $A$ a cofibrant $\caO$-algebra the enveloping algebra
$U_\caO(A)$ is cofibrant as an object in $\caC$.
\end{cor}

Hence also the second part of Theorem \ref{mod-mod-str} is proven.

%\begin{cor} \label{univ-vergl-2}
%Let $\caC$ be left proper,
%let $\caO \in \Op(\caC)$ be cofibrant in $\caC^\Sigma$ and
%$A \rightarrow A'$ be a map between cofibrant $\caO$-algebras
%which is a weak equivalence.
%Then the map $U_\caO(A) \rightarrow U_\caO(A')$ is a weak equivalence.
%\end{cor}
%\begin{proof}
%Let $Q\caO \rightarrow \caO$ be a cofibrant replacement of $\caO$.
%Then it follows from Proposition \ref{cell-univ-alg-str}
%that the vertical maps in the commutative square
%$$\xymatrix{U_{Q\caO}(A) \ar[r] \ar[d] & U_{Q\caO}(A') \ar[d] \\
%U_\caO(A) \ar[r] & U_\caO(A')}$$
%are weak equivalences, and the upper horizontal map
%is a weak equivalence by Corollary \ref{univ-vergl-1}, hence
%the claim follows.
%\end{proof}

\begin{cor} \label{univ-vergl-2}
Let $\caC$ be left proper,
let $f:\; \caO \rightarrow \caO'$ be a weak equivalence between
operads in $\caC$ both of which are cofibrant as objects in
$\caC^\Sigma$ and let $A$ be a cofibrant $\caO$-algebra.
Let $A'$ be the pushforward of $A$ with respect to $f$.
Then the induced maps $A \rightarrow A'$ and
$U_\caO(A) \rightarrow U_{\caO'}(A')$ are weak equivalences.
\end{cor}

\begin{defi} Let $\caC$ be left proper and let $\eins$ and the
domains of the maps in $I$ be cofibrant in $\caC$.
\begin{enumerate}
\item
For $\caO \in \Op(\caC)$ define the {\em derived category
of $\caO$-algebras} $D\Alg(\caO)$ to be $\ho \Alg(Q\caO)$, where
$Q\caO \rightarrow \caO$ is a cofibrant repalcement in $\Op(\caC)$.
Define the {\em derived $2$-category of $\caO$-algebras}
$D^{\le 2}\Alg(\caO)$ to be $\ho^{\le 2} \Alg(Q\caO)$.
\item
For $\caO \in \Op(\caC)$ and $A \in \Alg(\caO)$ define
the {\em derived category of $A$-modules $D(A\mmod)$} to be
$\ho (QA\mmod)$, where $QA \rightarrow A$ is a cofibrant replacement of $A$
in $\Alg(Q\caO)$ with $Q\caO \rightarrow \caO$ a cofibrant
replacement in $\Op(\caC)$.
\end{enumerate}
\end{defi}

Note that these definitions do not depend
(up to equivalence up to unique isomorphism
or up to equivalence up to isomorphism, which is itself
defined up to unique isomorphism in the case
of $D^{\le 2}\Alg(\caO)$) on the choices
by Corollary \ref{univ-vergl-2} and \cite[Theorem 2.4]{hov-mon},
that if $\caO \in \Op(\caC)$ is cofibrant in $\caC^\Sigma$
there is a canonical equivalence $D\Alg(\caO) \sim \ho \Alg(\caO)$
and that for a cofibrant $\caO \in \Op(\caC)$ 
and $A \in \Alg(\caO)$ which is cofibrant in $\caC$ there
is a canonical equivalence $D(A\mmod) \sim \ho (A \mmod)$.

\section{Functoriality}

In this section let $\caC$ be left proper 
and let $\eins$ and the domains of the maps in $I$ be cofibrant in $\caC$.
\begin{prop}
\begin{enumerate}
\item
There is a well defined $2$-functor
$$\ho^{\le 2} \Op(\caC) \rightarrow \Cat \;\text{,}$$
$$ \caO \mapsto D\Alg(\caO)$$
such that for any cofibrant operad $\caO$ in $\caC$ there
is a canonical equivalence $D\Alg(\caO) \sim \ho \Alg(\caO)$ and
every functor in the image of this $2$-functor
has a right adjoint.
\item For $\caO \in \Op(\caC)$ there is a well defined
$2$-functor $$D^{\le 2} \Alg(\caO) \rightarrow \Cat \;\text{,}$$
$$A \mapsto D(A\mmod)$$
such that for any cofibrant $A \in \Alg(Q\caO)$
($Q\caO \rightarrow \caO$ a cofibrant replacement) there
is a canonical equivalence $D(A\mmod) \sim \ho(A\mmod)$ and
every functor in the image of this $2$-functor has
a right adjoint.
\end{enumerate}
\end{prop}

\begin{rem}
The $2$-functor in the second part of the Proposition should be
well defined for an object $\caO \in \ho^{\le 3} \Op(\caC)$
and should depend on $\caO$ functorially.
\end{rem}

\begin{proof}
We prove the first part of the Proposition, the second one is similar.
Let $\caO, \caO' \in \Op(\caC)_{cf}$, $f,g \in \Hom(\caO,\caO')$ and
$\varphi$ a $2$-morphism from $f$ to $g$ in $\ho^{\le 2} \Op(\caC)$.
First of all it is clear that the pushforward functor
$f_*:\;\Alg(\caO) \rightarrow \Alg(\caO')$ is a left Quillen functor
between $J$-semi model categories by the definition of the
$J$-semi model structures.
We have to show that $\varphi$ induces a natural isomorphism
between $f_*$ and $g_*$ on the level of homotopy categories.
So let $\caO^\bullet$ be a cosimplicial frame on $\caO$.
$\varphi$ can be represented by a chain of $1$-simplices in
$\Hom(\caO^\bullet,\caO')$, and a homotopy between two representing
chains by a chain of $2$-simplices.
So we can assume that $\varphi$ is a $1$-simplex,
i.e. $\varphi \in \Hom(\caO^1,\caO')$.
We have maps $\caO \sqcup \caO \overset{i_0 \sqcup i_1}{\longrightarrow}
\caO^1 \overset{p}{\rightarrow} \caO$,
and $\ho \Alg(\caO^1) \rightarrow \ho \Alg(\caO)$ is an equivalence.
Hence for $A \in \ho \Alg(\caO)$ there is a unique isomorphism
$\varphi'(A):\; i_{0*}(A) \rightarrow i_{1*}(A)$ with $p_*(\varphi'(A))=\id$.
Then the $\varphi(\varphi'(A))$ define a natural isomorphism
between $(\varphi \circ i_0)_*$ and $(\varphi \circ i_1)_*$.
Now if we have a homotopy $\Phi \in \Hom(\caO^2,\caO')$,
the three natural transformations which are defined by the three
$1$-simplices of $\Phi$ are compatible, since on a given object
they are the images in
$\ho \Alg(\caO')$ of three compatible isomorphisms between the
three possible images of $A$ in $\ho \Alg(\caO^2)$.
\end{proof}

Let $f:\; \caO \rightarrow \caO'$ be a map of operads in
$\caC$ and let $A \in D^{\le 2} \Alg(\caO)$.
Then there is an adjunction
$$\xymatrix{D(A\mmod) \ar@<.6ex>[r] & D(f_*A \mmod) 
\ar@<.6ex>[l]}\;\text{.}$$
It follows that for $B \in D^{\le 2} \Alg(\caO')$ there
is also an adjunction
$$\xymatrix{D(f^*B \mmod) \ar@<.6ex>[r] & D(B \mmod) 
\ar@<.6ex>[l]}\;\text{.}$$
Of course for $A$ and $B$ as above and a map $f_*A \rightarrow B$
there is a similar adjunction.

Now let $\caD$ be a second left proper
symmetric monoidal cofibrantly generated
model category with suitable smallness assumptions on the
domains of the generating cofibrations and trivial cofibrations
(depending on which definition of $J$-semi model category one takes)
and with a cofibrant unit.
Let $L:\; \caC \rightarrow \caD$ be a symmetric monoidal
left Quillen functor with right adjoint $R$.
For objects $X,Y \in \caD$ there is
always a natural map $$R(X) \otimes R(Y) \rightarrow R(X \otimes Y)$$
adjoint to the map $$F(R(X) \otimes R(Y)) \cong FR(X) \otimes FR(Y)
\rightarrow X \otimes Y$$ which respects the associativity and
commutativity isomorphisms (so $R$ is a {\em pseudo symmetric monoidal}
functor). It follows that $L$ can be lifted to preserve
operad, algebra and module structures.

Hence there is induced a pair of adjoint functors
$$\xymatrix{\Op(\caC) \ar@<.6ex>[r]^-{L_\Op} & 
\Op(\caD) \ar@<.6ex>[l]^-{R_\Op}}\;\text{,}$$
which is a Quillen adjunction between $J$-semi model categories by
the definition of the model structures.

For $\caO \in \Op(\caC)$ there is induced a pair of adjoint functors
$$\xymatrix{\Alg(\caO) \ar@<.6ex>[r]^-{L_\caO} & 
\Alg(L_\Op(\caO)) \ar@<.6ex>[l]^-{R_\caO}}\;\text{,}$$
which is a Quillen adjunction between $J$-semi model categories
in the cases where $\caO$ is either cofibrant in $\Op(\caC)$
or cofibrant as an object in $\caC^\Sigma$.

So for $\caO \in \Op(\caC)$, $\caO' \in \Op(\caD)$ and
$f:\; L_\Op(\caO) \rightarrow \caO'$ a map there are induced
adjunctions
$$\xymatrix{D\Alg(\caO) \ar@<.6ex>[r] & 
D\Alg(\caO') \ar@<.6ex>[l]}\;\text{ and}$$
$$\xymatrix{D^{\le 2}\Alg(\caO) \ar@<.6ex>[r]^-{\Psi} & 
D^{\le 2}\Alg(\caO') \ar@<.6ex>[l]}\;\text{.}$$
Now let $A \in D^{\le 2}\Alg(\caO)$,
$B \in D^{\le 2}\Alg(\caO')$ and $g:\; \Psi(A) \rightarrow B$
be a map. Then there is induced an adjunction
$$\xymatrix{D(A\mmod) \ar@<.6ex>[r] & 
D(B \mmod) \ar@<.6ex>[l]}\;\text{.}$$

All the adjunctions are compatible (in an appropriate weak 
categorical sense) with compositions of the
maps which induce these adjunctions.

%There are three kinds of functoriality: We can change
%the underlying model category 
%$\caC$, we can change the operad 
%$\caO \in \Op(\caC)$ and we can change the algebra
%$A \in \Alg(\caO)$.

\section{$E_\infty$-Algebras} \label{einfty-alg}

Let $\caN$ be the operad in $\caC$ whose algebras are
just the commutative unital algebras in $\caC$, i.e.
$\caN(n)=\eins$ for $n \in \cN$, and let $\caP$ be the operad
whose algebras are objects in $\caC$ pointed by $\eins$,
i.e. $\caP(n)=\eins$ for $n=0,1$, $\caP(n)=0$ otherwise.
There is an obvious map $\caP \rightarrow \caN$.

\begin{defi}
\begin{enumerate}
\item An {\em $E_\infty$-operad} in $\caC$ is an operad $\caO$ in $\caC$
which is cofibrant as an object in $\caC^\Sigma$ together
with a map $\caO \rightarrow \caN$ which is a weak equivalence.
\item A {\em pointed} $E_\infty$-operad in $\caC$ is an
$E_\infty$-operad $\caO$ in $\caC$ together with a map
$\caP \rightarrow \caO$ such that the composition with the
map $\caO \rightarrow \caN$ is the canonical map $\caP \rightarrow \caN$.
\item A {\em unital} $E_\infty$-operad in $\caC$
is a pointed $E_\infty$-operad in $\caC$ such that the map
$\caP(0) \rightarrow \caO(0)$ is an isomorphism
(this is the same as an $E_\infty$-operad $\caO$ in $\caC$
such that the map $\caO(0) \rightarrow \caN(0)$ is an isomorphism).
\end{enumerate}
\end{defi}
The unit $\eins$ is an $\caN$-algebra, hence it is an algebra
for any $E_\infty$-operad.

We first want to show that under suitable conditions unital
$E_\infty$-operads always exist.

For $\caO \in \Op(\caC)$ let us denote by $\caO_{\le 1}$
the operad with $\caO_{\le 1} (0)=\caO(0)$,
$\caO_{\le 1} (1)=\caO(1)$ and $\caO_{\le 1} (n)=0$ for $n > 1$.
There is a canonical map $\caO_{\le 1} \rightarrow \caO$ in $\Op(\caC)$.
If $\caO$ is an $E_\infty$-operad there is also a map
$\caO_{\le 1} \rightarrow \caP$ in $\Op(\caC)$,
and we denote by $\widetilde{\caO}$ the pushout of $\caO$
with respect to this map.

\begin{lem} \label{unit-lem}
Let $\caO$ be an $E_\infty$-operad which admits a pointing.
\begin{enumerate}
\item Then there is a canonical equivalence $\Alg^u(\caO) \sim
\Alg(\widetilde{\caO})$, in particular an $\caO$-algebra is unital
if and only if it comes from an $\widetilde{\caO}$-algebra.
\item
Assume that $\caC$ is left proper, that $\eins$ is cofibrant
in $\caC$ and that $\caO$ is cofibrant in $\Op(\caC)$.
Then $\widetilde{\caO}$ is a unital $E_\infty$-operad
in $\caC$.
\end{enumerate}
\end{lem}
\begin{proof}
By Lemma \ref{alg-mor}(1)
an $\widetilde{\caO}$-algebra $A$ is the same as an $\caO$-algebra $A$
together with a map $\eins \rightarrow A$ such that the structure map
$\caO(0) \rightarrow A$ is the composition $\caO(0) \rightarrow \eins
\rightarrow A$.
Hence a unital $\caO$-algebra comes from an $\widetilde{\caO}$-algebra.
On the other hand
if $A$ is an $\widetilde{\caO}$-algebra we have to show that
the induced pointing $\eins \rightarrow A$ is a map of algebras.
This follows easily from the fact that the map $\caO(0)$ has
a right inverse (a pointing of $\caO$).
For the first part of the Lemma it remains to prove
that an $\caO$-algebra morphism between $\widetilde{\caO}$-algebras
is in fact an $\widetilde{\caO}$-algebra morphism,
which follows from Lemma \ref{alg-mor}(2).

Consider the commutative square
$$\xymatrix{\caO(0) \ar[r] \ar[d] & F_\caO(\eins) \ar[d] \\
\eins \ar[r] & F_{\widetilde{\caO}} (\eins)}$$
of $\caO$-algebras
and let $P$ be the pushout of the left upper triangle of
the square. We want to show that the canonical map $P \rightarrow
F_{\widetilde{\caO}}(\eins)$ is an isomorphism.
By the first part of the Lemma $P$ is an $\widetilde{\caO}$-algebra.
Now again by the first part of the Lemma it is easily seen
that $P$ has the same universal property as $\widetilde{\caO}$-algebra
as $F_{\widetilde{\caO}}(\eins)$.

So the above square is a pushout square in $\Alg(\caO)$, 
and hence by left properness
of $\Alg(\caO)$ over $\caC$ (Theorem \ref{alg-mod-str-1})
the right vertical arrow is a weak equivalence.
This implies that $\caO \rightarrow \widetilde{\caO}$
is a weak equivalence.
It remains to prove that $\widetilde{\caO}$ is cofibrant
as object in $\caC^\Sigma$, which follows from Corollary \ref{unit-cof}.
\end{proof}

Let us call a vertex $v \in V(T)$ of a tree $T \in \caT$
a {\em no-tail} vertex if one cannot reach a tail from $v$.
Let us call $T$ {\em $0$-special} if the only no-tail vertices
of $T$ are vertices of valency $0$. A proper $0$-special doubly colored
tree is a doubly colored tree which is $0$-special such that
any vertex of valency $0$ is old. Let $\widetilde{\caT}_\dc^p(n)$
be the set of isomorphism classes of such trees with $n$ tails.

\begin{lem} \label{unit-pushout}
Let $\caO=\colim_{i < \lambda} \caO_i$ be an operad in $\caC$
such that the transition maps $\caO_i \rightarrow \caO_{i+1}$
are pushouts of free operad maps on maps $g_i:\; K_i \rightarrow L_i$
in $\caC^\cN$ and such that $\caO_0$ is the initial operad.
Let $E \in \caC$ and let $\caO(0) \rightarrow E$
be a morphism in $\caC$. Let $\widetilde{E}$ be the operad with
$\widetilde{E}(0)=E$, $\widetilde{E}(1)=\eins$ and
$\widetilde{E}(n)=0$ for $n>1$. Let the squares
$$\xymatrix{\caO_{i,\le 1} \ar[r] \ar[d] & \caO_i \ar[d] \\
\widetilde{E} \ar[r] & \widetilde{\caO}_i}$$ be pushout squares
in $\Op(\caC)$, where either $i < \lambda$ or $i$ is the blanket.
Then $\widetilde{\caO}=\colim_{i<\lambda} \widetilde{\caO}_i$,
and every map $\widetilde{\caO}_i \rightarrow \widetilde{\caO}_{i+1}$
is a $\omega \times (\omega+1)$-sequence in $\caC^\cN$ as in
Proposition \ref{op-pushout}, where for $j < \omega$
$\caO_{(i,j)}(n)$ is a pushout of $\caO_{(i,j)-1}(n)$ in $\caC$ 
by the quotient of the map
$$\coprod \left( \bigotimes_{v \in V_\old(T)
\setminus U_\old(T)} \!\!\!\!\!\!\! 
\widetilde{\caO}_i (\val(v))\right)\otimes
e^{\Box(U_\old(T))}\Box \: 
\underset{v \in V_\new(T)}{\text{\huge $\Box$}}
g_i(\val(v)) \;\text{,}$$
where the coproduct is over all $T \in \widetilde{\caT}_\dc^p(n)$
with $\sharp V_\new(T)=i$ and $u_\old(T)=j$,
with respect to an equivalence relation analoguous
to the one in Proposition \ref{op-pushout}.
In particular we have $\widetilde{\caO}_i(0)=E$ for all $i<\lambda$
or $i$ the blanket.
\end{lem}
\begin{cor} \label{unit-cof}
Let the notation be as in the Lemma above
and assume that the maps $g_i$ are cofibrations in $\caC^\cN$
and that $E$ is cofibrant in $\caC$.
Then the operad $\widetilde{\caO}$ is cofibrant in $\caC^{\Sigma,\bullet}$.
\end{cor}
\begin{proof}
The proof is along the same lines as the proof of
Lemma \ref{op-pushout-cof}.
\end{proof}

For the rest of this section let us fix a pointed
$E_\infty$-operad $\caO$ in $\caC$.
An $\caO$-algebra $A$ is naturally pointed, i.e.
there is a canonical map $\eins \rightarrow A$,
but note that this need not be a map of algebras.
If it is, we say that $A$ is a {\em unital} $\caO$-algebra.
Let us denote the category of unital $\caO$-algebras
by $\Alg^u(\caO)$. This is just the category of
objects in $\Alg(\caO)$ under $\eins$.
If $\caO$ is unital, then every $\caO$-algebra is unital.

\begin{lem}
If $\eins$ is cofibrant in $\caC$ and $\caO$ is cofibrant
in $\Op(\caC)$ there is a $J$-semi model structure
on $\Alg^u(\caO)$ over $\caC$.
\end{lem}
\begin{proof}
In any $J$-semi model category $\caD$ over $\caC$ the category of objects
under an object from $\caD$ which becomes cofibrant in $\caC$
is again a $J$-semi model category over $\caC$.
\end{proof}

\begin{lem} \label{ealg-univ-vergl}
Assume that $\caC$ is left proper and
that the domains of the maps in $I$ are cofibrant.
Let $A \in \Alg(\caO)$ be cofibrant. Then the canonical
map of $A$-modules $U_\caO(A) \rightarrow A$ adjoint to
the pointing $\eins \rightarrow A$ is a weak equivalence.
\end{lem}
\begin{proof}
We can assume that $A$ is a cell $\caO$-algebra.
It is easy to see that the map $U_\caO(A) \rightarrow A$
is compatible with the descriptions of $A$ and $U_\caO(A)$
in Proposition \ref{cell-alg-str} and 
Proposition \ref{cell-univ-alg-str} as
transfinite compositions, and the map $\psi$ from
the map of part (3) of Proposition \ref{cell-univ-alg-str} to
the map of part (3) of Proposition \ref{cell-alg-str} 
is induced by the map $\caO(|f|+1)=\caO(|f|+1) \otimes \eins^{\otimes |f|}
\otimes \eins
\rightarrow \caO(|f|+1) \otimes \caO(1)^{\otimes |f|}\otimes 
\caO(0) \rightarrow \caO(|f|)$, which itself is induced
by the unit, the pointing and a structure map of $\caO$.
Since $\caO$ is an $E_\infty$-operad this is a weak equivalence,
hence since the domains of the maps in $I$ are cofibrant
$\psi$ is a weak equivalence. Now the claim follows
by transfinite induction and left properness of $\caC$.
\end{proof}
\begin{cor}
Assume that $\caC$ is left proper,
that the domains of the maps in $I$ are cofibrant
and that $\caO$ is cofibrant in $\Op(\caC)$.
Let $A \in \Alg(\caO)$ be cofibrant as object in $\caC$. Then the canonical
map of $A$-modules $U_\caO(A) \rightarrow A$ adjoint to
the pointing $\eins \rightarrow A$ is a weak equivalence.
\end{cor}
\begin{proof}
Let $QA \rightarrow A$ be a cofibrant replacement.
Then in the commutative square
$$\xymatrix{U_\caO(QA) \ar[r] \ar[d] & U_\caO(A) \ar[d] \\
QA \ar[r] & A}$$ the horizonrtal maps are weak equivalences
(the upper one by Corollary \ref{univ-vergl-1}) and the
left vertical arrow is a weak equivalence by the Lemma above,
hence the right vertical map is also a weak equivalence.
\end{proof}
\begin{cor} \label{cof-alg-univ-equiv}
Assume that $\caC$ is left proper and
that the domains of the maps in $I$ are cofibrant.
Let $A \rightarrow A'$ be a weak equivalence between
cofibrant $\caO$-algebras. Then the map
$U_\caO(A) \rightarrow U_\caO(A')$ is also a weak equivalence.
\end{cor}
\begin{proof}
This follows immediately from Lemma \ref{ealg-univ-vergl}.
\end{proof}

This Corollary has the consequence that under the assumptions
of the Corollary there is a canonical equivalence
$D(A\mmod) \sim \ho A\mmod$ for a cofibrant
$\caO$-algebra $A$.

\section{$\cS$-Modules and Algebras} \label{s-modules}

In this section we generalize the theories developed
in \cite{ekmm} and \cite{km}.

\begin{defi}
(I) A {\em symmetric monoidal category with pseudo-unit} is
a category $\caD$ together with
\begin{itemize}
\item
a functor $\_ \boxtimes \_ : \; \caD \times \caD \rightarrow \caD$,
\item
natural isomorphisms 
$(X \boxtimes Y) \boxtimes Z \rightarrow X \boxtimes (Y \boxtimes Z)$ 
and $X \boxtimes Y \rightarrow Y \boxtimes X$ which
satisfy the usual equations and
\item an object $\eins \in\caD$ with morphisms
$\eins \boxtimes X \rightarrow X$ (and hence morphisms
$X \boxtimes \eins \rightarrow X$ induced by the symmetry isomorphisms)
such that the diagram
$$\xymatrix{\eins \boxtimes (X \boxtimes Y) \ar[d] \ar[r] & X \boxtimes Y 
\\ (\eins \boxtimes X) \boxtimes Y \ar[ru]}
$$ commutes and such that the two possible maps
$\eins \boxtimes \eins \rightarrow \eins$ agree.
\end{itemize}
(II) A {\em symmetric monoidal functor} between symmetric
monoidal categories with pseudo-unit $\caD$ and $\caD'$ is
a functor $F:\; \caD \rightarrow \caD'$ together with
natural isomorphisms $F(X) \boxtimes F(Y) \rightarrow F(X \boxtimes Y)$
compatible with the associativity and commutativity isomorphisms and
with a map $F(\eins_\caD) \rightarrow \eins_{\caD'}$ compatible
with the unit maps.
\end{defi}
\begin{defi}
(I) A {\em symmetric monoidal model category with weak unit}
is a model category $\caD$ which is a symmetric monoidal category
with pseudo-unit such that the functor
$\caD \times \caD \rightarrow \caD$ has the
structure of a Quillen bifunctor (\cite[p. 108]{hov-book})
and such that the composition
$Q\eins \boxtimes X \rightarrow \eins \boxtimes X \rightarrow X$ is
a weak equivalence for all cofibrant $X \in \caD$, where
$Q\eins \rightarrow \eins$ is a cofibrant replacement.

(II) A {\em symmetric monoidal Quillen functor} between symmetric
monoidal model categories with weak unit $\caD$ and $\caD'$
is a left Quillen functor $\caD \rightarrow \caD'$
which is a symmetric monoidal functor between symmetric
monoidal categories with pseudo-unit such that the composition
$F(Q \eins_\caD) \rightarrow F(\eins_\caD) \rightarrow \eins_{\caD'}$
is a weak equivalence.
\end{defi}

The homotopy category of a symmetric monoidal model category
with weak unit is a closed symmetric monoidal category.

Let us assume now that $\caC$ is either simplicial
(i.e. there is a symmetric monoidal left Quillen functor
$\sset \rightarrow \caC$) or that there is a symmetric monoidal
left Quillen functor $\comp_{\ge 0}(\Ab) \rightarrow \caC$, where
$\comp_{\ge 0}(\Ab)$ is endowed with the projective model structure.
In both cases we denote by $\caL$ the image of the linear
isometries operad in $\Op(\caC)$ via either the simplicial
complex functor or the simplicial complex functor followed
by the normalized chain complex functor.
Clearly $\caL$ is a unital $E_\infty$-operad.
Let $\cS:=\caL(1)$. $\cS$ is a ring with unit in $\caC$
which is cofibrant as an object in $\caC$.

As in \cite{ekmm} or \cite{km} we define a tensor
product on $\cS\mmod$ by
$$M\boxtimes N:=\caL(2)\otimes_{\cS \otimes \cS} M \otimes N \;\text{.}$$
\cite[Theorem V.1.5]{km} and \cite[Lemma V.1.6]{km} also work in
our context, hence $\cS\mmod$ is a symmetric monoidal category with
pseudo-unit.

There is an internal Hom in $\cS\mmod$ given by
$$\uhom^\boxtimes(M,N) := \uhom_\cS(\caL(2)\otimes_\cS M,N)\;\text{,}$$
where, when forming $\caL(2) \otimes_\cS M$, $\cS$ acts on $\caL(2)$
through $\cS =\eins\otimes \cS \rightarrow \cS \otimes \cS$,
when forming $\uhom_\cS$, $\cS$ acts on $\caL(2)\otimes_\cS M$ via its
left action on $\caL(2)$ and the left action of $\cS$ on
$\uhom^\boxtimes (M,N)$ is induced through the right action
of $\cS$ on $\caL(2)$ through $\cS= \cS \otimes \eins \rightarrow 
\cS \otimes \cS$.

There is an augmentation $\cS \rightarrow \eins$ which is
a map of algebras with unit.

\begin{prop}
The category $\cS\mmod$ is a cofibrantly generated symmetric
mono-idal model category with weak unit with
generating cofibrations $\cS \otimes I$ and generating trivial
cofibrations $\cS \otimes J$.
The functor $\caC \rightarrow \cS\mmod$,
$X \mapsto \cS \otimes X$, is a Quillen equivalence,
and its left inverse, the functor $\cS\mmod \rightarrow \caC$,
$M \mapsto M \otimes_\cS \eins$, is a symmetric monoidal
Quillen equivalence.
Moreover there is a closed action of $\caC$ on $\cS\mmod$.
\end{prop}
\begin{proof}
That $R\mmod$ is a cofibrantly generated model category together
with a closed action of $\caC$ on it is true for any associative
unital ring $R$ in $\caC$ which is
cofibrant as an object in $\caC$.

Let $f$ and $g$ be cofibrations in $\caC$.
The $\boxtimes$-pushout product of $\cS \otimes f$ and $\cS \otimes g$
is isomorphic to $\caL(2) \otimes (f \Box g)$. As a left $\cS$-module
$\caL(2)$ is (non canonically) isomorphic to $\cS$, hence
$\caL(2) \otimes (f \Box g)$ is a cofibration $\cS\mmod$,
and it is trivial if one of $f$ or $g$ is trivial.
%To construct the tensor equivalence use the fact
%that for cofibrant $X,Y \in \caC$ there is a weak equivalence
%$\cS \otimes X \rightarrow X$ in $\cS\mmod$ induced by
%the augmentation $\cS \rightarrow \eins$ and the diagram
%$$(\cS \otimes X) \boxtimes (\cS \otimes Y) = \caL(2) \otimes X \otimes Y
%\rightarrow X \otimes Y \leftarrow \caS \otimes (X \otimes Y)\;\text{.}$$
To show that for a cofibrant $\cS$-module $M$ the map
$Q\eins \boxtimes M \rightarrow M$ is a weak equivalence we
can assume that $M$ is a cell $\cS$-module and we can take
$Q\eins = \cS$. Then $M$ is a transfinite composition where
the transition maps are pushouts of maps $f:\; \cS \otimes K \rightarrow
\cS \otimes L$, where $K \rightarrow L$ is a cofibration in 
$\caC$ with $K$ cofibrant. But the composition 
$\cS \boxtimes \cS \rightarrow \eins \boxtimes \cS
\rightarrow \cS$ is a weak equivalence between cofibrant objects
in $\cS\mmod$, hence the composition $\cS \boxtimes f \rightarrow \eins
\boxtimes f \rightarrow f$ is a weak equivalence between
cofibrations in $\cS\mmod$. So by transfinite induction the composition
$\cS \boxtimes M \rightarrow \eins \boxtimes M \rightarrow M$
is a weak equivalence between cofibrant objects in $\cS\mmod$.
\end{proof}

Note that in the simplicial case $\eins \boxtimes \cS$ is
cofibrant in $\caC$, hence for cofibrant $M$ both maps
$\cS \boxtimes M \rightarrow \eins \boxtimes M \rightarrow M$
are weak equivalences.

Let $\cS\mmod^u$ be the category of unital $\cS$-modules,
i.e. the objects in $\cS\mmod$ under $\eins \in \cS\mmod$.
For $M \in \cS\mmod^u$ and $N \in \cS\mmod$ there
are the products $M \vartriangleleft N$ and $N \vartriangleright M$,
and for $M,N \in \cS\mmod^u$ there is the product $M \boxdot N$.
These products are defined as in \cite[Definition V.2.1]{km}
and \cite[Definition V.2.6]{km}.

$\cS\mmod^u$ is a symmetric monoidal category with $\boxdot$
as tensor product.

Analoguous to \cite[Theorem V.3.1]{km} and
\cite[Theorem V.3.3]{km} we have

\begin{prop} \label{alg-char}
\begin{itemize}
\item
$\Alg(\caL)$ is naturally equivalent to the category of
commutative rings with unit in $\cS\mmod^u$.
Hence for $A,B \in \Alg(\caL)$ there is a natural isomorphism
$A \sqcup B \cong A \boxdot B$.
\item
For $A \in \Alg(\caL)$ an $A$-module $M$ is the same as
an $\cS$-module $M$ together with a
map $A \vartriangleleft M \rightarrow M$ satisfying the usual identities.
\end{itemize}
\end{prop}

For $A \in \Alg(\caL)$ let $\Comm(A)$ be the category of commutative
unital $A$-algebras in $\cS\mmod^u$, i.e. the objects in
$\Alg(\caL)$ under $A$. In particular we have 
$\Alg(\caL) \sim \Comm(\eins_\cS)=: \Comm_\caC$, where
we denote by $\eins_\cS$ the algebra $\eins$ in $\cS\mmod^u$.

For the rest of the section let us make the following
\begin{assump} \label{add-assump}
The model category $\caC$ is left proper
and $\eins$ and the domains of the maps in $I$ are cofibrant in $\caC$.
\end{assump}

\begin{cor}
\begin{itemize}
\item
$Comm_\caC$ is a cofibrantly generated $J$-semi
model category.
\item 
For any cofibrant $A \in \Comm_\caC$
the category $\Comm(A)$ is also a cofibrantly generated
$J$-semi model category.
\item  
If $A \rightarrow B$ is a weak equivalence 
between cofibrant $A,B \in \Comm_\caC$,
then the induced functor $\Comm(A) \rightarrow \Comm(B)$
is a Quillen equivalence.
\end{itemize}
\end{cor}
\begin{proof}
Follows from Theorem \ref{alg-mod-str-2}.
\end{proof}

\begin{defi}
For $A \in \Comm_\caC$ let $D\Comm(A)$ be
$\ho \Comm(QA)$ for $QA \rightarrow A$ a cofibrant
replacement in $\Comm_\caC$, and let
$D^{\le 2}\Comm(A):=\ho^{\le 2} \Comm(QA)$.
The $2$-functor $\Comm_\caC \rightarrow \Cat$,
$A \mapsto D\Comm(A)$, descents
to a $2$-functor $D^{\le 2}\Comm_\caC \rightarrow \Cat$,
$A \mapsto D\Comm(A)$.
\end{defi}

Let $A \in \Comm_\caC$ and $M,N \in A\mmod$.
As in \cite[Definition V.5.1]{km} or
\cite[Remark V.5.2]{km} we define the tensor product
$M \boxtimes_A N$ as the coequalizer in the diagram
$$\xymatrix{(M \vartriangleright A) \boxtimes N \cong
M \boxtimes (A \vartriangleleft N) \ar@<.6ex>[r]^-{m \boxtimes \id}
\ar@<-.6ex>[r]_-{\id \boxtimes m} &
M \boxtimes N \ar[r] & M \boxtimes_A N} \;\text{ or}$$
$$\xymatrix {M \boxtimes A \boxtimes N 
\ar@<.6ex>[r]^-{m \boxtimes \id} \ar@<-.6ex>[r]_-{\id \boxtimes m} &
M \boxtimes N \ar[r] & M \boxtimes_A N}\;\text{.}$$

With this product the category $A\mmod$ has the structure
of a symmetric monoidal category with pseudo-unit, where
the pseudo-unit is $A$.
As for $\cS$-modules one can define products
$\vartriangleleft_A$, $\vartriangleright_A$ and
$\boxdot_A$. There is also an analogue of Proposition 
\ref{alg-char} for $A$-algebras and modules over $A$-algebras.

The free $A$-module functor $\cS\mmod \rightarrow A\mmod$
is given by $M \mapsto A \vartriangleleft M$. More
generally for $A \rightarrow B$ a map in $\Comm_\caC$ the
pushforward of modules is given by $M \mapsto B \vartriangleleft_A M$.
In particular there is a canonical isomorphism of
$A$-modules $U_\caL(A) \cong A \vartriangleleft \cS$.

\begin{lem} \label{boxtimes-base-change}
Let $A \rightarrow B$ be a map in $Comm_\caC$, let
$M,N \in A\mmod$ and $P \in B\mmod$. Then there are canonical ismorphisms
$$M \boxtimes_A P \cong (B \vartriangleleft_A M) \boxtimes_B P
\;\text{ and}$$
$$(B \vartriangleleft_A M) \boxtimes_B (B \vartriangleleft_A N)
\cong B \vartriangleleft_A (M \boxtimes_A N) \;\text{.}$$
\end{lem}
\begin{proof}
Similar to the proof of \cite[Proposition V.5.8]{km}.
\end{proof}

For $M,N \in A\mmod$ define 
the internel Hom $\uhom_A^\boxtimes(M,N)$
in $A\mmod$ as the equalizer
$$\xymatrix{\uhom_A^\boxtimes(M,N) \ar[r] & \uhom^\boxtimes (M,N)
\ar@<.6ex>[r] \ar@<-.6ex>[r] & \uhom^\boxtimes(A \vartriangleleft M,N)}$$
like in \cite[Definition V.6.1]{km}.

\begin{prop}
For a cofibrant $A \in \Comm_\caC$ the category $A\mmod$ is a
cofibrantly generated symmetric monoidal model
category with weak unit with generating cofibrations
$A \vartriangleleft (\cS \otimes I)$ and generating trivial
cofibrations $A \vartriangleleft (\cS \otimes J)$.
If $f:\; A \rightarrow B$ is a map in $\Comm_\caC$
between cofibrant algebras
the pushforward $f_*:\; A\mmod \rightarrow B\mmod$
is a symmetric monoidal Quillen functor
which is a Quillen equivalence if $f$ is a weak equivalence.
\end{prop}
\begin{proof}
$A\mmod$ is a cofibrantly generated model category by
Theorem \ref{mod-mod-str}(1).
Let $f$ and $g$ be cofibrations in $\caC$.
By Lemma \ref{boxtimes-base-change}
the $\boxtimes_A$-pushout product of the maps $A \vartriangleleft
(\cS \otimes f)$ and $A \vartriangleleft (\cS \otimes g)$
is given by $A \vartriangleleft (\caL(2) \otimes (f \Box g))$,
hence since $\caL(2)\cong \cS$ as $\cS$-modules this
is a cofibration in $A\mmod$, and it is trivial if one of $f$ or
$g$ is trivial.

Note that $A \vartriangleleft \cS$ is cofibrant in $A\mmod$
and that the map $A \vartriangleleft \cS \cong
U_\caL(A) \rightarrow A$ is a weak equivalence by
Lemma \ref{ealg-univ-vergl}.
So we have to show
that for cofibrant $M \in A\mmod$ the map $(A \vartriangleleft \cS)
\boxtimes_A M \rightarrow M$ is a weak equivalence, which follows
from the fact that the maps of the form
$(A \vartriangleleft \cS) \boxtimes_A
(A \vartriangleleft (\cS\otimes f)) \rightarrow A \vartriangleleft
(\cS \otimes f)$ for cofibrations $f \in \caC$
with cofibrant domain are weak
equivalences between cofibrations in $A\mmod$.
The first part of the last statement follows
from Lemmas \ref{boxtimes-base-change} and \ref{ealg-univ-vergl},
and the second part by Corollary \ref{cof-alg-univ-equiv}.
\end{proof}

For $A \in D^{\le 2}\Comm_\caC$ there is unambiguously defined
the closed symmetric monoidal category $D(A\mmod)$
with tensor product denoted by $\otimes_A$.
The assignment $A \mapsto D(A\mmod)$ is a $2$-functor
$D^{\le 2}\Comm_\caC \rightarrow \Cat^{\text{sm}}$, where
$\Cat^{\text{sm}}$ is the $2$-category of symmetric monoidal
categories, such that the image functors of all maps in
$D^{\le 2}\Comm_\caC$ have right adjoints.

%Let $A \in \Comm_\caC$.
%There is a well defined functor
%$$\sqcup_A^\cL:\; \Comm(A) \times \Comm(A) \rightarrow \Comm(A)\;\text{,}$$
%which maps $A$-algebras $B$ and $C$ to the homotopy pushout 
%$B \sqcup_A^\cL C$ of $B$ and $C$ over $A$.

Let $$\xymatrix{B \ar[r]^-{g'} \ar@{=>}[dr]^\varphi & B' \\ 
A \ar[r]^-g \ar[u]_-f & 
A' \ar[u]_-{f'}}$$ be a commutative square in $D^{\le 2}\Comm_\caC$.
Let $M \in D(B\mmod)$. Then we have a base change morphism
$$g_*f^*M \rightarrow {f'}^*g'_* M$$ defined to be the adjoint of
the natural map $f^*M \rightarrow f^*{g'}^*g'_* M
\overset{\varphi}{\cong} g^*{f'}^* g'_* M$ or equivalently
of the map $f'_*g_*f^* M \overset{\varphi}{\cong} g'_*f_* f^* M
\rightarrow g'_*M$.

The base change morphism is natural with respect to
composition of commutative squares.

The following statement is trivial in the context of usual
commutative algebras, but is a rather strong structure result
in our context.

\begin{prop} \label{base-change-iso}
Let the notation be as above.
If the square is a homotopy pushout, then
the base change morphism $g_*f^*M \rightarrow {f'}^*g'_* M$
is an isomorphism.
\end{prop}

The proof will be given in the next section.

Let $A \rightarrow B$ be a map in $D^{\le 2}\Comm_\caC$.
Let $M \in D(A\mmod)$ and $N \in D(B\mmod)$.
There is a projection morphism
$$M \otimes_A f^*N \rightarrow f^*(f_*M \otimes_B N)$$
adjoint to the natural map
$f_*(M \otimes_A f^*N)=f_*M \otimes_B f_*f^*N \rightarrow f_*M \otimes_B N$.
Note that for $B$-modules $M',N'$ there is a natural map
$f^* M' \otimes_A f^* N' \rightarrow f^*(M' \otimes_B N')$,
and the projection morphism is equivalently described as the
composition $M \otimes_A f^*N \rightarrow f^* f_*M \otimes_A f^*N
\rightarrow f^*(f_*M \otimes_B N)$.

\begin{prop} \label{proj-mor-iso}
Let the notation be as above. Then the projection morphism
$M \otimes_A f^*N \rightarrow f^*(f_*M \otimes_B N)$ is
an isomorphism.
\end{prop}

We give the proof in the next section.

Let a square in $D^{\le 2} \Comm_\caC$ be
given as above and let $M \in D(B\mmod)$, $N \in D(A'\mmod)$
and $P \in D(A\mmod)$.
Set $M':=f^*M$, $N':=g^*N$, $\widetilde{M}:=g'_* M$,
$\widetilde{N}:=f'_* N$ and $\widetilde{P}:=
g'_*f_*P \cong f'_*g_* P$.

\begin{lem} \label{base-proj-vertr}
Let the notation be as above. Then the diagram
$$\xymatrix@C.7cm{(M' \otimes_A P) \otimes_A N' \ar[r] &
g^*(g_*(M' \otimes_A P) \otimes_{A'} N) \ar[r] & g^*{f'}^*((\widetilde{M}
\otimes_{B'} \widetilde{P}) \otimes_{B'} \widetilde{N}) \\
M' \otimes_A (P \otimes_A N') \ar[u] \ar[r] & 
f^*(M \otimes_B f_*(P \otimes_A N')) \ar[r] & f^*{g'}^*(\widetilde{M}
\otimes_{B'} (\widetilde{P} \otimes_{B'} \widetilde{N})) \ar[u]}\;\text{,}$$
where in the first two horizontal maps the projection morphism is
applied and in the second two an adjunction and the base change
morphism, commutes.
\end{lem}
\begin{proof}
Let $F:=g^*{f'}^* \overset{\varphi}{\cong} f^*{g'}^*$.
One checks that both compositions are equal to the composition
$M' \otimes_A P \otimes_A N' \rightarrow F^* \widetilde{M} \otimes_A
F^* \widetilde{P} \otimes_A F^* \widetilde{N} 
\rightarrow F^*(\widetilde{M} \otimes_{B'} \widetilde{P}
\otimes_{B'} \widetilde{N})$, where the first arrow is a tensor
product of obvious objectwise morphisms.
\end{proof}

Let $A \in D^{\le 2}\Comm_\caC$.
We can use the two Propositions above to give
the natural functor $M:\; D\Comm(A) \rightarrow D(A\mmod)$
a symmetric monoidal structure with respect to the coproduct on $D\Comm(A)$
and the tensor product $\otimes_A$ on $D(A\mmod)$:
We use the fact that $D\Comm(A)$ is equivalent
to the $1$-truncation of $A \downarrow D^{\le 2}\Comm_\caC$.
So let $B \leftarrow A \rightarrow C$ be a triangle
in $D^{\le 2}\Comm_\caC$ and complete it by a homotopy pushout to
a square with upper right corner $B \sqcup_A C$.
First we apply the base change isomorphism to the unit $\eins_B$
in $D(B\mmod)$, which says that 
there is a natural isomorphism
$$(C \rightarrow B \sqcup_A C)^* (\eins_{B \sqcup_A C}) \cong
(A \rightarrow C)_*(M(B))\;\text{.}$$
Applying $(A \rightarrow C)^*$ to the left hand side of this isomorphism
we get $M(B \sqcup_A C)$, applying this map to the right
hand side we get $M(B) \otimes_A M(C)$ by the projection
formula. This establishes the isomorphism
$M(B) \otimes_A M(C) \cong M(B \sqcup_A C)$. That this isomorphism
respects the commutativity isomorphisms follows from Lemma
\ref{base-proj-vertr} with $P=\eins_A$.
That it respects the associativity isomorphisms 
for objects $f:\; A \rightarrow B$, $h:\; A \rightarrow C$
and $g: \; A \rightarrow A'$ in $A \downarrow \ho^{\le 2} \caD$
also follows from Lemma \ref{base-proj-vertr}
with $M=\eins_B$, $N=\eins_{A'}$ and $P=h^*\eins_C$
and a diagram chase.

%Let $A \in \Comm_\caC$ be cofibrant.
%To prove the two Propositions above we need a modified
%model structure on $A\mmod$, which exists if the
%assumption below is satisfied.
%We will give another proof of the Propositions
%in the next section, which is a bit more technical, but
%which does not rely on this assumption.

%Let $\underline{f}=(f_1,\ldots,f_k)$ 
%with maps $f_i$ in $I$ and let $\underline{n}=(n_1,\ldots,n_k)$
%with $n_i \in \cN$, $i=1 \ldots,k$. Let $m=\sum_{i=1}^k n_i$
%and let $H < \prod_{i=1}^k \Sigma_{n_i}$ be a subgroup.
%Consider the map $F_{\underline{f},\underline{n}}:=
%A\vartriangleleft (\caL(m)\otimes_H ...$
%Tgeht wohl doch nicht, da man Faserungen nicht kontrollieren kann...

\section{Proofs}

%Let $A \in D\Comm_\caC$ and
%$B,C \in D\Comm(A)$. Let $B \sqcup_A C$ be the
%homotopy pushout of $B$ and $C$ over $A$.
%So $\sqcup_A$ is a symmetric monoidal structure
%on $\Comm(A)$.
%Note that there is a natural functor
%$M:\; D\Comm(A) \rightarrow D(A\mmod)$ which
%assigns to an $A$-algebra the underlying $A$-module.
%In the following we want to show that this
%functor is naturally symmetric monoidal,
%i.e. that there exist natural isomorphisms
%R$$M(B \sqcup_A C) \cong M(B) \otimes_A M(C)\;\text{.}$$
%Note that this already follows in the case
%$A=\eins$ from the Proposition above.

In this section we give the proofs of Propositions \ref{base-change-iso}
and \ref{proj-mor-iso}.
Assume throughout that Assumption \ref{add-assump}
is fulfilled.

We need the concept of
operads in $A\mmod$ for $A \in \Comm_\caC$.
We also give the definition of a pointed operad, because
it is needed in the Appendix.  In the context of symmetric
monoidal categories with pseudo-unit a pointed operad is
not just an operad $\caO$ together with a pointing of $\caO(0)$,
the domains of the structure maps also have to be adjusted
(see below).

So let us fix $A \in \Comm_\caC$.
%In the following we write $\boxtimes$, $\vartriangleleft$,
%$\vartriangleright$ and $\boxdot$ for
%$\boxtimes_A$, $\vartriangleleft_A$,
%$\vartriangleright_A$ and $\boxdot_A$.
Let $A\mmod^u$ be the category of pointed $A$-modules,
i.e. the category of objects in $A\mmod$ under $A$.
For $M$ a pointed or unpointed $A$-module
and $N$ a pointed or unpointed $A$-module
let $M \circledast N$ be either $M \boxtimes_A N$,
$M \vartriangleleft_A N$, $M \vartriangleright_A N$
or $M \boxdot_A N$, depending on
whether $M$ and $N$ are unpointed, $M$ is pointed and
$N$ is unpointed, $M$ is unpointed and $N$ is pointed
or $M$ and $N$ are pointed.
$M \circledast N$ is an object in $A\mmod$ unless
both $M$ and $N$ are pointed in which case it
is an object in $A\mmod^u$.
Note that for $M_1,\ldots, M_n$ $A$-modules each
of them either pointed or unpointed the product
$M_1 \circledast \cdots \circledast M_n$ is well defined,
despite the fact that for different bracketings
of this expression the symbols for which $\circledast$
actually stands can be different.

\begin{defi}
An {\em operad} $\caO$ in $A\mmod$ is an object
$\caO(n) \in (A\mmod)[\Sigma_n]$
for each $n \in \cN$, where $\caO(1)$ is pointed, together with maps
$$\caO(m) \circledast \caO(n_1) \circledast \cdots \circledast
\caO(n_m) \rightarrow \caO(n) \;\text{,}$$
where $m,n_1\ldots,n_m \in \cN$ and $n=\sum_{i=1}^m n_i$,
such that the usual diagrams for these structure maps
commute.
A {\em pointed operad} in $A\mmod$ is the same as above with the exception
that $\caO(0)$ is also pointed.
\end{defi}
Let $\Op(A\mmod)$ be the category of operads in $A\mmod$ and
$\Op^p(A\mmod)$ the category of pointed operads
in $A\mmod$. A pointed operad $\caO$ in $A\mmod$ is called {\em unital}
if the pointing $A \rightarrow \caO(0)$ is an isomorphism.
Let $\Op^u(A\mmod)$ be the category of unital operads in $A\mmod$.

Let $(A\mmod)^{\Sigma,\bullet \bullet}$ be the category of
collections of objects $\caO(n) \in (A\mmod)[\Sigma_n]$,
which are pointed for $n=0,1$ and unpointed otherwise.
As for ordinary operads we have free (pointed) operad functors $F$
starting from the categories $(A\mmod)^\cN$,
$(A\mmod)^\Sigma$, $(A\mmod)^{\Sigma,\bullet \bullet}$ in the pointed case
and various other pointed versions of these categories
to $\Op(A\mmod)$ or $Op^p(A\mmod)$.
Note that if $A$ is cofibrant all these source categories
of the functors $F$ are model categories.

\begin{thm} \label{op-mod-str-modi}
Let $A$ be cofibrant in $\Comm_\caC$.
Then the category $\Op(A\mmod)$
(resp. $\Op^p(A\mmod)$) is a cofibrantly generated
$J$-semi model category over 
$(A\mmod)^{\Sigma, \bullet}$ (resp. 
over $(A\mmod)^{\Sigma,\bullet \bullet}$)
with generating cofibrations $FF_AI$ and generating trivial cofibrations
$FF_AJ$.
If $\caC$ is left proper, then 
$\Op(A\mmod)$ (resp. $\Op^p(A\mmod)$)
is left proper relative to 
$(A\mmod)^{\Sigma,\bullet \bullet}$ (resp. relative to
$(A\mmod)^{\Sigma, \bullet}$).
If $\caC$ is right proper, so are $\Op(A\mmod)$
and $\Op^p(A\mmod)$.
\end{thm}

Let $f$ be a map in $A\mmod$ or $A\mmod^u$ and
let $g$ be a map in $A\mmod$ or $A\mmod^u$.
Let $f \boxstar g$ be the pushout product of $f$ and $g$
with respect to the product $\circledast$.
$f \boxstar g$ is a map in $A\mmod$ unless both
$f$ and $g$ are maps in $A\mmod^u$ in which case 
$f \boxstar g$ is a map in $A\mmod^u$.

Note that if $A$ is cofibrant the category $A\mmod^u$ has
a natural model structure as category of objects under
$A$ in the model category $A\mmod$.
Note however that $A\mmod^u$ is {\em not} symmetric monoidal
(with potential tensor product $\boxdot_A$), since
this product is not closed.

\begin{lem} \label{unit-pushouts}
Let $A$ be cofibrant in $\Comm_\caC$, let $f$
be a cofibration in $A\mmod$ or $A\mmod^u$,
let $g$ be a cofibration in $A\mmod$ or $A\mmod^u$
let $M$ be cofibrant in $A\mmod$ or $A\mmod^u$ and
let $N$ be cofibrant in $A\mmod$ or $A\mmod^u$.
Then 
\begin{itemize}
\item the pushout product $f \boxstar g$ is a cofibration
in $A\mmod$ or $A\mmod^u$ which is trivial if $f$ or $g$ is,
\item the product $M \circledast f$ is a cofibration in
$A\mmod$ or $A\mmod^u$ which is trivial if $f$ is and
\item the product $M \circledast N$ is cofibrant in
$A\mmod$ or $A\mmod^u$.
\end{itemize}
There is also a version of this statement
when the map or object in $A\mmod$ has a right action
of a discrete group $G$ and the other
map or object is in $A\mmod^u$ (resp. when both
maps or objects are in $A\mmod$ and have actions of
discrete groups $G$ and $G'$). The resulting
map or object is then a cofibration or cofibrant object
in $(A\mmod)[G]$ (resp. $(A\mmod)[G \times G']$).
\end{lem}
Note that in a symmetric monoidal category
cases 2 and 3 would be special cases of case 1.
\begin{proof}
It suffices to show this for relative cell complexes
$f$ and $g$ and cell complexes
$M$ and $N$, for which it follows for the first case
by writing the pushout
product of a $\lambda$-sequence and a $\mu$-sequence
as a $\lambda\times \mu$-sequence.
Let $M \in A\mmod^u$.
Then if $A \rightarrow M$ is a $\lambda$-sequence, $M$ itself
is a $(1+\lambda)$-sequence in $A\mmod$. One concludes now
by writing the
products in cases 2 and 3 again as appropriate sequences.
The cases with group actions work in the same way.
\end{proof}

We remark now that there are versions of
Propositions \ref{free-op} and \ref{op-pushout}
for $\Op^p(A\mmod)$ where all tensor products are
repaced by $\circledast$-products and all pushout products
by the $\circledast$-pushout product $\boxstar$.
There is also a version of Lemma \ref{op-pushout-cof},
from which Theorem \ref{op-mod-str-modi} follows
in the same way as Theorem \ref{op-mod-str}.

\begin{defi}
Let $\caO \in \Op(A\mmod)$ (resp. $\caO \in \Op^p(A\mmod)$). 
\begin{enumerate}
\item An {\em $\caO$-algebra}
is an object $B \in A\mmod$
(resp. $B \in A\mmod^u$) together with maps
$$\caO(n) \circledast B^{\circledast n} \rightarrow A$$
satisfying the usual identities. The category of $\caO$-algebras
is denoted by $\Alg(\caO)$.
\item Let $B \in \Alg(\caO)$. A {\em $B$-module} is an
object $M \in A\mmod$ together with maps
$$\caO(n+1) \circledast B^{\circledast n}
\circledast M \rightarrow M$$
satisfying the usual identities. The category of $B$-modules
is denoted by $B\mmod$.
\end{enumerate}
\end{defi}

Let $\caO \in \Op^{(p)}(A\mmod)$. The free $\caO$-algebra
functor $F_\caO:\; A\mmod \rightarrow \Alg(\caO)$
is given by $$F_\caO(M)= \coprod_{n \ge 0} \caO(n)
\circledast_{\Sigma_n} M^{\boxtimes_A n}\;\text{.}$$
In the pointed case $F_\caO$ factors through $A\mmod^u$.

As in section \ref{algebras} one shows the

\begin{thm} \label{alg-mod-str-modi}
Let $A$ be cofibrant in $\Comm_\caC$ and let $\caO \in \Op(A\mmod)$
(resp. $\caO \in \Op^p(A\mmod)$).
\begin{enumerate}
\item If $\caO$ is cofibrant
the category $\Alg(\caO)$ is a cofibrantly generated
$J$-semi model category over
$A\mmod$ with generating cofibrations $F_\caO F_A I$ and generating
trivial cofibrations $F_\caO F_A J$.
If $\caC$ is left proper (resp. right proper), then
$\Alg(\caO)$ is left proper relative to $A\mmod$ (resp. right proper).
\item
Let $\caO$ be cofibrant as an
object in $(A\mmod)^{\Sigma, \bullet}$
(resp. in $(A\mmod)^{\Sigma, \bullet \bullet}$).
Then $\Alg(\caO)$ is a cofibrantly generated
$J$-semi model category with generating cofibrations
$F_\caO F_A I$ and generating trivial cofibrations $F_\caO F_A J$.
If $\caC$ is right proper, so is $\Alg(\caO)$.
\end{enumerate}
\end{thm}

Let $\caN_A \in \Op(A\mmod)$
(resp. $\caN_A^u \in \Op^u(A\mmod)$) be the operad with 
$\caN_A(n)=A$ (resp. $\caN_A^u(n)=A$)
for $n \in \cN$ and the natural structure maps.
Note that both categories $\Alg^{(u)}(\caN_A)$ are {\em not} equivalent to the
category $\Comm(A)$ but there are functors
$$C_{\caN_A}^{(u)}:\; \Alg(\caN_A^{(u)}) \rightarrow \Comm(A)\;\text{,}$$
which are defined to be the left adjoints of the pullback
functors $\Comm(A) \rightarrow \Alg(\caN_A^{(u)})$. These adjoints
exist since they exist on free algebras and every
algebra is a coequalizer of two maps between free algebras
(as is always the case for algebras over a monad).

\vskip.2cm

Let $\caO \in \Op^{(p)}(A\mmod)$ and $B \in \Alg(\caO)$.
As for ordinary algebras one defines the universal
enveloping algebra $U_\caO(B)$ as the
quotient of the tensor algebra
$$\coprod_{n \ge 0} \caO(n+1) \circledast_{\Sigma_n} B^{\circledast n}$$
by the usual relations.
$U_\caO(B)$ is an associative unital algebra in $A\mmod$,
hence it is an $A_\infty$-algebra in $\caC$
(i.e. an algebra over the operad $\caL$
considered as a {\em non-$\Sigma$} operad), which also
has a universal enveloping algebra $U_\caL(U_\caO(B)) \in \Ass(\caC)$.
One has canonical equivalences $$B\mmod \sim U_\caO(B)\mmod
\sim U_\caL(U_\caO(B))\mmod \; \text{.}$$
Let $F_B:\; A\mmod \rightarrow B\mmod$ be the free $B$-module functor.

\vskip.2cm

As in section \ref{modules} one shows the

\begin{thm} \label{mod-mod-str-modi}
Let $A$ be cofibrant in $\Comm_\caC$,
let $\caO \in \Op(A\mmod)$ (resp. $\caO \in \Op^p(A\mmod)$)
and $B \in \Alg(\caO)$.
Let one of the following two conditions be satisfied:
\begin{enumerate}
\item $\caO$ is cofibrant as an object in $(A\mmod)^{\Sigma, \bullet}$
(resp. in $(A\mmod)^{\Sigma,\bullet \bullet}$) and $B$
is a cofibrant $\caO$-algebra.
\item $\caO$ is cofibrant in $\Op(A\mmod)$ (resp. $\Op^p(A\mmod)$)
and $A$ is cofibrant as an object in $A\mmod$ (resp. in $A\mmod^u$).
\end{enumerate}
Then there is cofibrantly generated model structure on
$B\mmod$ with generating cofibrations $F_B F_A I$ and generating
trivial cofibrations $F_B F_A J$.
\end{thm}

\begin{defi}
An {\em $E_\infty$-operad} (resp. pointed $E_\infty$-operad)
in $A\mmod$ is an object 
$\caO \in \Op(A\mmod)$ (resp. $\caO \in \Op^p(A\mmod)$)
which is cofibrant as an object in $(A\mmod)^{\Sigma, \bullet}$
(resp. in $(A\mmod)^{\Sigma, \bullet \bullet}$)
together with a map $\caO \rightarrow \caN_A$ which is a weak equivalence.
A pointed $E_\infty$-operad
$\caO$ is called {\em unital} if it is unital as an object in 
$\Op^p(A\mmod)$.
\end{defi}

For $\caO$ a pointed $E_\infty$-operad in $A\mmod$ let us define
the operad $\widetilde{\caO}$ in the same way as in section \ref{einfty-alg}.
Then we have analogues of Lemmas \ref{unit-lem}
and \ref{unit-pushout} and Corollary \ref{unit-cof}.
So we are able to construct a unital $E_\infty$-operad
in $A\mmod$ by first taking a cofibrant resolution $\caO \rightarrow
\caN_A$ in $\Op(A\mmod)$ and then forming $\widetilde{\caO}$.
This will be relevant in the Appendix.

Let $B \in \Alg(\caO)$ be cofibrant.
As in Lemma \ref{ealg-univ-vergl} one can show that the map
$U_\caO(B) \rightarrow B$ adjoint to the pointing
$A \rightarrow B$ is a weak equivalence.

\vskip.2cm

For the rest of this section let us fix an unpointed
$E_\infty$-operad
$\caO$ in $A\mmod$ (we could also take a pointed one).
Let $\pi$ be the map $\caO \rightarrow \caN_A$.

\begin{lem} \label{comm-alg-equiv}
Let $A$ be cofibrant in $\Comm_\caC$.
Then the composition $$\xymatrix{\Alg(\caO) \ar[r]^-{\pi_*} &
\Alg(\caN_A) \ar[r]^-{C_{\caN_A}} & \Comm(A)}$$
is a Quillen equivalence.
\end{lem}
\begin{proof}
This follows from the fact that for a cofibrant $A$-module $M$
the map $$\caO(n)\circledast_{\Sigma_n} M^{\boxtimes_A n}
\rightarrow M^{\boxtimes_A n}/\Sigma_n$$ is a weak equivalence.
\end{proof}

\begin{lem} \label{mod-equiv}
Let $A$ be cofibrant in $\Comm_\caC$ and let $B \in \Alg(\caO)$
be cofibrant. 
Then the functor $$B\mmod \rightarrow (C_{\caN_A} \circ \pi)_*(B)\mmod$$
is a Quillen equivalence.
\end{lem}
\begin{proof}
This follows from the fact that the map
$U_\caL(U_\caO(B)) \rightarrow U_\caL((C_{\caN_A} \circ \pi)_*(B))$
is a weak equivalence, which follows itself from the
description of these algebras in terms of transfinite compositions
as in Propositions \ref{cell-alg-str} and \ref{cell-univ-alg-str}.
\end{proof}

\begin{lem} \label{univ-cof}
Let $A$ be cofibrant in $\Comm_\caC$ and
let $B \in \Alg(\caO)$ be cofibrant. Then $U_\caO(B)$
is cofibrant as object in $A\mmod^u$.
\end{lem}
\begin{proof}
Follows by the description of $U_\caO(B)$ as
in Proposition \ref{cell-univ-alg-str}.
\end{proof}

\begin{cor} \label{cof-mod}
Let $A$ be cofibrant in $\Comm_\caC$. Then
for cofibrant $B \in \Alg(\caO)$ and cofibrant $M \in B\mmod$
the underlying $A$-module $M$ is cofibrant in $A\mmod$.
\end{cor}
\begin{proof}
Follows from Lemmas \ref{univ-cof} and \ref{unit-pushouts}
and transfinite induction.
\end{proof}

\begin{lem} \label{s-prod-str}
Let $\mu$ and $\lambda$ be ordinals and let $S_{\mu,+}$ and
$S_{\lambda,+}$ be as in Proposition \ref{cell-alg-str}. Then
there is a (necessarily unique) isomorphism
$$\varphi:\; S_{\lambda+ \mu,+} \cong S_{\mu,+} \times S_{\lambda,+}$$
of well-ordered sets.
\end{lem}
\begin{proof}
There is a natural inclusion $S_{\lambda,+} \hookrightarrow
S_{\lambda + \mu, +}$, and $\varphi$ maps its image to $\{*\} \times
S_{\lambda,+}$ in the natural way.
Now let $f \in S_{\lambda + \mu}$ with $f(i) \notin \cN$
for some $i \in \mu$. There is a segment $M_f \subset S_{\lambda + \mu}$
starting at $f$ which is isomorphic to $S_{\lambda, +}$ as a well-ordered
set. Via this identification
$S_\lambda$ corresponds to all $f' \in S_{\lambda+ \mu}$
with $f'|_\mu = f|_\mu + \frac{1}{2}$.
Then $\varphi$ maps $M_f$ to $\{f|_\mu \} \times S_{\lambda, +}$ if
$i > 0$ and to $\{f|_\mu - \frac{1}{2}\} \times S_{\lambda,+}$
if $i=0$. It is easy to see that this way $\varphi$ is well-defined,
bijective and order-preserving.
\end{proof}
\begin{rem}
If $f \in S_{\lambda,+}$ and $g\in S_{\mu,+}$ are successors,
then $\varphi$ maps $(f \sqcup g) -1$ to $(g-1,f-1)$.
\end{rem}

\begin{proof}[Proof of Proposition \ref{base-change-iso}]
By Lemmas \ref{comm-alg-equiv} and \ref{mod-equiv}
we can work in $\Alg(\caO)$.
So let $B,C \in \Alg(\caO)$ be cofibrant.
Let us denote the coproduct in $\Alg(\caO)$ by $\sqcup_A$.
We have to prove the base change isomorphism for the diagram
$$\xymatrix{B \ar[r]^{g'} & B \sqcup_A C \\ A \ar[u]_f \ar[r]^g & C
\ar[u]_{f'}}\;\text{.}$$
Let $M \in B\mmod$ be cofibrant.
Then $f^*M$ is cofibrant in $A\mmod$ by Corollary \ref{cof-mod}.
Hence the base change morphism is represented by the
morphism of $U_\caO(C)$-modules $U_\caO(C) \vartriangleleft_A M
\rightarrow U_\caO(B \sqcup_A C) \vartriangleleft_{U_\caO(B)} M$
which is adjoint to the map
$M \cong A \vartriangleleft_A M \rightarrow 
U_\caO(B \sqcup_A C) \vartriangleleft_{U_\caO(B)} M$.
We can assume that $M$ is a cell module.
Then by transfinite induction we are reduced to the
following statement:
Let $K \in A\mmod$ be cofibrant.
Then the map $U_\caO(C) \vartriangleleft_A (U_\caO(B) \vartriangleleft_A K)
\rightarrow U_\caO(B \sqcup_A C) \vartriangleleft_A K$ is
a weak equivalence.
By Lemma \ref{unit-pushouts} this follows if we show that the
map of $B$-modules
$\psi:\; U_\caO(B) \boxdot_A U_\caO(C) \rightarrow U_\caO(B \sqcup C)$
(where we exchanged the roles of $B$ and $C$) is a weak equivalence.
It suffices to prove this for cell algebras $B$ and $C$.
So let $B=\colim_{i<\lambda} B_i$, where
the transition maps are given by pushouts by
cofibrations $g_i: K_i \rightarrow L_i$ in $A\mmod$ with cofibrant domain
as in Proposition \ref{cell-alg-str}.
Similarly let $C=\colim_{i< \mu} C_i$, where the transition maps
are given by pushouts by cofibrations 
$h_i: M_i \rightarrow N_i$ in $A\mmod$ with cofibrant domain.
Then the map $0 \rightarrow U_\caO(B \sqcup_A C)$
is described as in Proposition \ref{cell-univ-alg-str}
by a $S_{\lambda+\mu,+}$-sequence (1).
Since the maps $0 \rightarrow U_\caO(B)$ resp. 
$0 \rightarrow U_\caO(C)$ are
$S_{\lambda,+}$- resp. $S_{\mu,+}$-sequences, the
map $0 \rightarrow U_\caO(B) \boxdot_A U_\caO(C)$ is a 
$S_{\mu,+}\times S_{\lambda,+}$-sequence (2)
by Lemma \ref{cell-pushout-str} (this also holds in the
case of a symmetric monoidal category with pseudo-unit).
Let $\alpha:\; S_{\mu,+}\times S_{\lambda,+} \rightarrow S_{\lambda+\mu,+}$
be the isomorphism of well-ordered sets of Lemma \ref{s-prod-str}.
Let $f \in S_{\lambda,+}$ and $f' \in S_{\mu,+}$ be successors.
Then $\alpha$ identifies $(f \sqcup f')-1$ and $(f'-1,f-1)$,
and the relevant pushouts in the sequences (1) and (2) are by maps
$$\caO(|f \sqcup f'|+1) \circledast_{\Sigma_{f \sqcup f'}}
\text{\LARGE $\Box_*$}_{i< \lambda} g_i^{\Box_* f(i)} \Box_*
\text{\LARGE $\Box_*$}_{i< \mu} h_i^{\Box_* f'(i)} \; \text{ and}$$
$$\caO(|f|+1) \circledast \caO(|f'|+1) 
\circledast_{\Sigma_f \times \Sigma_{f'}}
\text{\LARGE $\Box_*$}_{i< \lambda} g_i^{\Box_* f(i)} \Box_*
\text{\LARGE $\Box_*$}_{i< \mu} h_i^{\Box_* f'(i)} \; \text{.}$$
It is easy to see by transfinite induction
that the map $\psi$ is compatible with sequences
(1) and (2) via the identification $\alpha$ on the indexing sets
and with the above pushouts by the map induced by
the tensor multiplication map
$\caO(|f|+1) \circledast \caO(|f'|+1) \rightarrow
\caO(|f \cup f'|+1)$ which inserts the second object into the last
slot of the first object. This map is a weak equivalence because
$\caO$ is an $E_\infty$-operad, hence the claim follows
by transfinite induction.
\end{proof}

\begin{proof}[Proof of Proposition \ref{proj-mor-iso}]
By Lemmas \ref{comm-alg-equiv} and \ref{mod-equiv}
we can assume that we have a cofibrant $\widetilde{B} \in \Alg(\caO)$,
a cofibrant $\widetilde{N} \in \widetilde{B}\mmod$
and a cofibrant $M \in A\mmod$ and prove the projection isomorphism
for $M$ and the image $N$ of $\widetilde{N}$ in $B\mmod$,
where $B$ is the image of $\widetilde{B}$ in $\Comm(A)$.
Since $\widetilde{N}$ is cofibrant as $A$-module by Corollary
\ref{cof-mod}
the projection morphism is represented by the composition
$$M \boxtimes_A \widetilde{N} \rightarrow M \boxtimes_A N
\cong (B \vartriangleleft_A M) \boxtimes_B N\;\text{,}$$
where the isomorphism at the second place is from
Lemma \ref{boxtimes-base-change}.
So we have to show that the first map is a weak equivalence.
We can assume that $\widetilde{N}$
is a cell module. Then by transfinite induction one
is left to show that for a cofibrant $A$-module $K$
the map $M \boxtimes_A (U_\caO(\widetilde{B}) \vartriangleleft_A K)
\rightarrow M \boxtimes_A (B \vartriangleleft_A K)$ is a weak
equivalence. But this map is the map from the free $\widetilde{B}$-module
on $M \boxtimes_A K$ to the free $B$-module on $M \boxtimes_A K$,
which is a weak equivalence by Lemma \ref{mod-equiv}. Hence we are finished.
\end{proof}

\section{Appendix}

Assume that Assumption \ref{add-assump}
is fulfilled.

In this section we give an alternative definition
of a product on the derived category of
modules over an algebra in $D^{\le 2}\Comm_\caC:=D^{\le 2}\Alg(\caN)$
without using the special properties of
the linear isometries operad.
Unfortunately it seems to be rather ugly (or difficult)
to construct associativity and commutativity
isomorphisms, and we did not try hard to do this!
Note that $D^{\le 2}\Comm_\caC$ is the same up to canonical equivalence
as the category denoted with the same symbol in section \ref{s-modules}.
If $\caO$ is a unital $E_\infty$-operad and $A \in D^{\le 2}\Comm_\caC$,
then there is a representative $\widetilde{A} \in \ho^{\le 2} \Alg(\caO)$
which is well defined up to an isomorphism which itself is well defined
up to a unique $2$-isomorphism. There is a similar statement
for a lift of $A$ into $\Alg(\caO)$.

Let us first treat the case where $\caC$ is simplicial, since
it is a bit nicer.
Let $\caO$ be a pointed $E_\infty$-operad in $\sset$
and denote by $\caO$ also its image in $\Op(\caC)$.
In $\sset$ the diagonal $\bigtriangleup:\;\caO \rightarrow \caO \times \caO$
is a map of operads, hence we also have a map of operads
$\caO \rightarrow \caO \otimes \caO$ in $\Op(\caC)$.

We will define a tensor product on $\ho A\mmod$ for
a cofibrant $\caO$-algebra $A$.

First note that for $\caO$-algebras $A$ and $B$ the
tensor product $A \otimes B$ is a $\caO \otimes \caO$-algebra,
hence also a $\caO$-algebra via $\bigtriangleup$.
Also for an $A$-module $M$ and a $B$-module $N$ the
tensor product $M \otimes N$ has a natural structure
of an $A\otimes B$-module.
If $A,B$ are unital there are induced maps in $\Alg^u(\caO)$
$A=A\otimes \eins \rightarrow A \otimes B$ and
$B = \eins \otimes B \rightarrow A \otimes B$.

\begin{prop} \label{pushout-prod-vergl}
Assume that $\caO$ is either unital or cofibrant in $\Op(\caC)$.
Let $A,B \in \Alg^u(\caO)$ be cofibrant.
Then the canonical map $A \sqcup B \rightarrow A \otimes B$
in $\Alg^u(\caO)$ induced by the maps
$A \rightarrow A\otimes B$ and $B \rightarrow A \otimes B$
is a weak equivalence.
\end{prop}
\begin{proof}
This proof is very similar to a part of the proof of Proposition
\ref{base-change-iso}.
By Lemma \ref{unit-lem} we are reduced to the case
where $\caO$ is unital.
It suffices to prove the claim for cell algebras $A$ and $B$.
So let $A=\colim_{i<\lambda} A_i$, where
the transition maps are given by pushouts by
maps $g_i: K_i \rightarrow L_i$ as in Proposition \ref{cell-alg-str}.
Similarly let $B=\colim_{i< \mu} B_i$, where the transition maps
are given by pushouts by maps $h_i: M_i \rightarrow N_i$.
Then the map $0 \rightarrow A \sqcup B$
is described by Proposition \ref{cell-alg-str}
by a $S_{\lambda+\mu,+}$-sequence (1).
Since the maps $0 \rightarrow A$ resp. $0 \rightarrow B$ are
$S_{\lambda,+}$- resp. $S_{\mu,+}$-sequences, the
map $0 \rightarrow A\otimes B$ is a 
$S_{\mu,+}\times S_{\lambda,+}$-sequence (2).
Let $\alpha:\;S_{\lambda+\mu,+} \rightarrow S_{\lambda,+}\times S_{\mu,+}$
be the isomorphism of well-ordered sets of Lemma \ref{s-prod-str}.
Let $f \in S_{\lambda,+}$ and $f' \in S_{\mu,+}$ be successors.
Then $\alpha$ identifies $(f \sqcup f')-1$ and $(f'-1,f-1)$.
The relevant pushouts in the sequences (1) and (2) are by maps
$$\caO(|f \sqcup f'|) \otimes_{\Sigma_{f \sqcup f'}}
\text{\LARGE $\Box$}_{i< \lambda} g_i^{\Box f(i)} \Box
\text{\LARGE $\Box$}_{i< \mu} h_i^{\Box f'(i)} \; \text{ and}$$
$$\caO(|f|) \otimes \caO(|f'|) \otimes _{\Sigma_f \times \Sigma_{f'}}
\text{\LARGE $\Box$}_{i< \lambda} g_i^{\Box f(i)} \Box
\text{\LARGE $\Box$}_{i< \mu} h_i^{\Box f'(i)} \; \text{,}$$
and again one shows by transfinite induction that
the map 
$\psi:\;A \sqcup B \rightarrow A\otimes B$ is compatible with sequences
(1) and (2) via the identification $\alpha$ on the indexing sets
and with the above pushouts by the map induced by
$$\xymatrix{\caO(|f| + |f'|) \ar[r]^-{\bigtriangleup} &
\caO(|f| + |f'|) \otimes \caO(|f| + |f'|) \ar[r]^-{\beta \otimes \gamma}
& \caO(|f|) \otimes
\caO(|f'|)} \;\text{,}$$
where $\beta$ inserts the pointing $\eins \rightarrow \caO(0)$
into the last $|f'|$ slots of $\caO(|f|+|f'|)$ and $\gamma$
inserts the pointing into the first $|f|$ slots.
This map is a weak equivalence since $\caO$ is an $E_\infty$-operad,
so our claim follows by transfinite induction and the assumptions.
\end{proof}

Assume that $\caO$ is either unital or cofibrant in $\Op(\caC)$.
For any cofibrant $\caO$-algebra $A$ let $Q_A$ denote a cofibrant
replacement functor in $A\mmod$.
Let $A \in \Alg^u(\caO)$ be cofibrant.
Then the map $A \sqcup A \rightarrow A \otimes A$ is a weak
equivalence. 
Now define a functor $$T:\; A\mmod \times A\mmod \rightarrow A\mmod
\; \text{ by}$$
$$T(M,N):=(A \sqcup A \rightarrow A)_* (Q_{(A\sqcup A)} 
(Q_A M \otimes Q_A N))\;\text{.}$$
It is clear that $T$ descents to a functor
$$T:\; D(A\mmod) \times D(A\mmod) \rightarrow D(A\mmod)\;\text{.}$$
We will see that this functor is naturally isomorphic
to the tensor product defined in section \ref{s-modules}.

\vskip.2cm

Now we skip the restriction of $\caC$ being simplicial.
Let $\caO$ be a unital $E_\infty$-operad in $\caC$ which always
exists by Lemma \ref{unit-lem}.
Then the operad $\caO \otimes \caO$ is also a unital $E_\infty$-operad.
Let $A,B \in \Alg(\caO \otimes \caO)$. Let
$\pi_1: \; \caO \otimes \caO \rightarrow \caO \otimes \caN \cong \caO$
and $\pi_2:\; \caO \otimes \caO \rightarrow \caN \otimes \caO \cong \caO$
be the two projections and define
$A_i:= \pi_{i,*}A$, $B_i:=\pi_{i,*}B$, $i=1,2$.
Note that $\pi_1$ and $\pi_2$ are weak equivalences.
There are maps 
$$A_1 \otimes \eins \rightarrow A_1 \otimes B_2\;\text{ and}$$ 
$$\eins \otimes B_2 \rightarrow A_1 \otimes B_2$$ of
$\caO \otimes \caO$-algebras and natural isomorphisms
of $\caO \otimes \caO$-algebras $A_1\otimes \eins \cong \pi_1^*A_1$
and $\eins \otimes B_2 \cong \pi_2^* B_2$, which are
on the underlying objects in $\caC$ the isomorphisms
$A_1\otimes \eins \cong A_1$
and $\eins \otimes B_2 \cong B_2$.
Using the adjunction units $A \rightarrow \pi_1^*A_1$
and $B \rightarrow \pi_2^*B_2$ we finally get
maps $A \rightarrow A_1\otimes B_2$ and
$B \rightarrow A_1 \otimes B_2$, hence
a map $$A \sqcup B \rightarrow A_1\otimes B_2$$
of $\caO \otimes \caO$-algebras.

\begin{prop}
Let $A,B \in \Alg(\caO \otimes \caO)$ be cofibrant.
Then the map $A \sqcup B \rightarrow A_1 \otimes B_2$ constructed
above is a weak equivalence.
\end{prop}
\begin{proof}
The proof of this Proposition is exactly the same as the one
for Proposition \ref{pushout-prod-vergl}, except that this
time the relevant pushouts in the sequences (1) and (2) are by maps
$$(\caO(|f \sqcup f'|)\otimes \caO(|f \sqcup f'|))
\otimes_{\Sigma_{f \sqcup f'}}
\text{\LARGE $\Box$}_{i< \lambda} g_i^{\Box f(i)} \Box
\text{\LARGE $\Box$}_{i< \mu} h_i^{\Box f'(i)} \; \text{ and}$$
$$\caO(|f|) \otimes \caO(|f'|) \otimes _{\Sigma_f \times \Sigma_{f'}}
\text{\LARGE $\Box$}_{i< \lambda} g_i^{\Box f(i)} \Box
\text{\LARGE $\Box$}_{i< \mu} h_i^{\Box f'(i)} \; \text{.}$$
The map $A \sqcup B \rightarrow A_1 \otimes B_2$ is again compatible
with these pushouts by the map induced by
$$\xymatrix{
\caO(|f| + |f'|) \otimes \caO(|f| + |f'|) \ar[r]^-{\beta \otimes \gamma}
& \caO(|f|) \otimes
\caO(|f'|)} \;\text{,}$$
where $\beta$ inserts the pointing $\eins \rightarrow \caO(0)$
into the last $|f'|$ slots of $\caO(|f|+|f'|)$ and $\gamma$
inserts the pointing into the first $|f|$ slots.
This map is again a weak equivalence since $\caO$ is an $E_\infty$-operad,
so we are done.
\end{proof}

Let $D\Comm_\caC:=D\Alg(\caN)$.

\begin{cor}
The natural functor $M:\;D\Comm_\caC \rightarrow \ho \caC$
has a natural symmetric monoidal structure with respect to the
coproduct on $D\Comm_\caC$ and the tensor product on $\ho \caC$.
\end{cor}

If $\cS$-modules are available in $\caC$
it is clear that this symmetric monoidal structure
is naturally isomorphic to the one constructed at the end
of section \ref{s-modules}.

Let now $A \in \Alg(\caO \otimes \caO)$ be cofibrant. 
Note that for $M,N \in A\mmod$ the tensor product
$\pi_{1,*} M \otimes \pi_{2,*} N$ is an $A_1 \otimes A_2$-module,
hence also an $A \sqcup A$-module.
Consider the functor
$$T:\; A\mmod \times A\mmod \rightarrow A\mmod\;\text{,}$$
$$(M,N) \mapsto (A \sqcup A \rightarrow A)_*
(Q_{A \sqcup A}(\pi_{1,*}(Q_A M) \otimes 
\pi_{2,*}(Q_A N)))\;\text{.}$$
It is again clear that $T$ descents to a functor
$$T:\; D(A\mmod) \times D(A\mmod) \rightarrow D(A\mmod)\;\text{.}$$

To see that this functor is
isomorphic to the previous functor $T$ in the simplicial case
one takes the previous $\caO$ to be $\caO \otimes \caO$
and looks at the map of $\caO \otimes \caO$-algebras
(obtained via the diagonal) $A \otimes A \rightarrow
(A_1 \otimes \eins)\otimes (\eins \otimes A_2)$.
The last algebra is isomorphic to the $\caO \otimes \caO$-algebra
$A_1 \otimes A_2$.
Hence for $A$-modules $M$ and $N$ we get a map of
$A \otimes A$-modules $M\otimes N \rightarrow M_1 \otimes N_2$
which is a weak equivalence. From this one gets the
natural isomorphism we wanted to construct.

It remains to show that in the cases $\caC$
receives a symmetric monoidal left Quillen functor
from $\sset$ or $\comp_{\ge 0}(\Ab)$ the functor $T$ is isomorphic
to the tensor product $\otimes_A$ defined in section \ref{s-modules}.

To do this let $\caO$ be a unital $E_\infty$-operad in
$\cS\mmod=\eins_\cS\mmod$ and let $\overline{\caO}:=\caO\otimes_\cS \eins$ be
its image in $\Op(\caC)$.
The operad $\caO \circledast \caO$ (which is defined componentwise)
is also a unital $E_\infty$-operad whose image in $\Op(\caC)$
is $\overline{\caO} \otimes \overline{\caO}$.
Then by the above procedure one can define a tensor product on
$\ho (A\mmod)$ for a cofibrant $\caO \circledast \caO$-algebra $A$,
and it is easy to see that this coincides (after the appropriate
identifications) with the product $T$ defined above on
$\ho (\overline{A}\mmod)$ ($\overline{A}$ is
the image of $A$ in $\Alg(\overline{\caO} \otimes \overline{\caO})$)
on the one hand and with the product $\boxtimes_{A'}$
on $\ho (A'\mmod)$, where $A'$ is the image of $A$ in $\Comm_\caC$,
on the other hand.

\end{document}